\documentclass[10pt,a4paper]{amsart}
\usepackage{amsfonts,layout}
\usepackage{bm,epsfig}

\newcount\minute
\newcount\hour
\def\currenttime{%
	\minute\time
	\hour\minute
	\divide\hour60
	\the\hour:\multiply\hour60\advance\minute-\hour\the\minute}
\def\draftnote{{\it \today \quad  \currenttime \hfill  tex-file :   \jobname}}

\catcode`@=11
\@addtoreset{equation}{section}

\catcode`@=12
\newtheorem{Theorem}{Theorem}[section]

\newtheorem{Proposition}{Proposition}[section]
\newtheorem{Lemma}{Lemma}[section]

\newtheorem{Hyp.}{Hyp.}[section]

\begin{document}

\title[]{Null controllability of parabolic equations with interior degeneracy and one-sided control}

\author{P. Cannarsa} 
\address{Dipartimento di Matematica, Universit\`a di Roma "Tor Vergata",
Via della Ricerca Scientifica, 00133 Roma, Italy}
\email{cannarsa@mat.uniroma2.it}

\author{R. Ferretti} 
\address{Dipartimento di Matematica e Fisica, Universit\`a Roma Tre, L.go S. Leonardo Murialdo, 1, 00146 Roma, Italy} \email{ferretti@mat.uniroma3.it}

\author{P. Martinez} 
\address{Institut de Math\'ematiques de Toulouse, UMR CNRS 5219, Universit\'e Paul Sabatier Toulouse III \\ 118 route de Narbonne, 31 062 Toulouse Cedex 4, France} \email{Patrick.Martinez@math.univ-toulouse.fr}

\thanks{This research was partly supported by Istituto Nazionale di Alta Matematica through the European Research Group GDRE CONEDP}

\begin{abstract}
For $\alpha\in (0,2)$ we study the null controllability of the parabolic operator
$$Pu= u_t - (\vert x \vert ^\alpha u_x)_x\qquad (1<x<1),$$
which degenerates at the interior point $x=0$, for locally distributed controls acting only one side of the origin (that is, on some interval $(a,b)$ with $0<a<b<1$). 
Our main results guarantees that $P$ is null controllable if and only if it is weakly degenerate, that is, $\alpha \in (0,1)$. So, in order to steer the system to zero, one needs controls to act on both sides of the point of degeneracy in the strongly degenerate case $\alpha\in [1,2)$.

Our approach is based on spectral analysis and the moment method. Indeed, we completely describe the eigenvalues and  eigenfunctions of the associated stationary operator in terms of Bessel functions and their zeroes for both weakly and strongly degenerate problems. Hence, we  obtain lower $L^2$ bounds for the eigenfunctions on the control region in the case $\alpha \in [0,1)$ and deduce the lack of observability in the case of $\alpha \in [1,2)$. We also provide numerical evidence to illustrate our theoretical results.
\end{abstract}
\maketitle

\section{Introduction}

\subsection{Presentation of the problem and main results} 
Degenerate parabolic equations have received increasing attention in recent years because of  their connections with several applied domains such as climate science~(\cite{bu68,bu69}, \cite{se69}, \cite{ghi76}, \cite{di97}),  populations genetics~(\cite{E-1976, ET-1993}), vision~(\cite{Manfredini04}), and mathematical finance~(\cite{BlaSch73})---just to mention a few. Indeed, in all these fields, one is naturally led to consider parabolic problems where  the diffusion coefficients lose uniform ellipticity. Different situations may occur: degeneracy (of uniform ellipticity) may take place at the boundary or in the interior of the space domain. Moreover, the equation may be degenerate on a small set or even on the whole domain. 

From the point of view of control theory, interesting phenomena have been pointed out for degenerate parabolic equations. We refer the reader to \cite{sicon2008,memoire} and the discussion below for the study  of null controllability of  boundary-degenerate parabolic operators. A similar analysis for interior-degenerate equations, associated with certain classes of hypoelliptic diffusion operators, was developed in \cite{Karine-Grushin1, Karine-Grushin3} (see also \cite{Karine-Grushin2}) for Grushin type structures, and in \cite{Karine-Grushin4} for the Heisenberg operator.

In this paper, intrigued by numerical tests (section \ref{sec-num}), we investigate the null controllability of a degenerate parabolic operator in one space dimension, which  degenerates at a single point inside the space domain, under the action of a locally distributed control supported only on one side of the  domain with respect to the point of degeneracy. In formulas, we consider the   problem
\begin{equation}
\label{eq-control}
\begin{cases}
u_t - (\vert x \vert ^\alpha u_x)_x = h(x,t) \chi _{(a,b)} (x) , \quad x\in (-1,1) \\
u(-1,t)=0=u(1,t) , \\
u(x,0)=u_0 (x) ,
\end{cases}
\end{equation}
assuming either $0<a<b<1$ or $-1<a<b<0$. Observe that this is the most general form a locally distributed control can be reduced to for the above operator.

In brief, we will  prove that null controllability:
\begin{itemize}
\item  fails for $\alpha \in [1,2)$,
\item  holds true when $\alpha \in [0,1)$.
\end{itemize}
Consequently, the controllability properties of the above operator, when degeneracy occurs inside the domain, are quite different from those of boundary degeneracy described in \cite{sicon2008, memoire}. In particular, in order for \eqref{eq-control} to be null controllable when $\alpha \in [1,2)$, the control support must lie on both sides of the degeneracy point.  We now proceed to describe our results more precisely and comment on the literature on this subject.



The case of boundary degeneracy  was addressed in \cite{CMV4, sicon2008} for equations in one space dimension, that is, for the problem
\begin{equation}
\label{eq-control-typ}
\begin{cases}
u_t - (x ^\alpha u_x)_x = h(x,t) \chi _{(a,b)} (x) , \quad x\in (0,1) \\
\begin{cases}
u(0,t)=0 =u(1,t) & \text{ if } \alpha \in [0,1), \\
(x^\alpha u_x)(0,t) =0 = u(1,t)& \text{ if } \alpha >1, 
\end{cases}
\\
u(x,0)=u_0 (x).
\end{cases}
\end{equation}
For the above equation, by deriving Carleman estimates with degeneracy-adapted weights, it was proved that null controllability holds if and only if $\alpha \in [0,2)$. 
These results were later extended in various directions---yet still limited to boundary degeneracy, see \cite{memoire} and the references therein.


The analysis of problem \eqref{eq-control} that we develop in the present paper, based on a detailed study of the associated spectral problem, allows us to discover some interesting properties, both positive and negative from the point of view of null controllability. More precisely,  we obtain:
\begin{itemize}
\item {\em Negative results for $\alpha \in [1,2)$.} This  fact is a little surprising if compared with the behaviour of  the usual problem \eqref{eq-control-typ}. The negative result we prove means that, when $\alpha \in [1,2)$, the degeneracy is too strong to allow the control to act on the other side of the domain with respect to the point of degeneracy. However, null controllability still holds true for those initial condition that are supported in the same region as the control.

\item {\em Positive results for $\alpha \in [0,1)$.} The proof of the fact that the control is sufficiently strong to cross  the degeneracy point does require to use fine properties of Bessel functions. We also give a sharp estimate of the blow-up rate of the null controllability cost as $\alpha \to 1^-$.
\end{itemize}

Degenerate parabolic equations with one (or more) degeneracy point inside the domain have also been studied by the flatness method developed by Martin-Rosier-Rouchon in \cite{flatness1, flatness2, flatness3} (see also Moyano \cite{moyano} for some strongly degenerate equations). More specifically, one can use the null controllability result with boundary control derived in \cite{flatness3} to construct a locally distributed control which steers the initial datum to $0$ for $\alpha \in [0,1)$. On the other hand, neither our analysis of the cost in the weakly degnerate case, nor our negative result for the strongly degenerate case seem to be attainable by the flatness approach.

Parabolic equations with interior degeneracy were also considered by Fragnelli-Mugnai in \cite{Fra-Mug2,Fra-Mug1}, where positive null controllability results were obtained for a general class of coefficients. Their approach,  based on Carleman estimates,  gives the controllability result when the control region is on both sides of the space domain with respect to the degeneracy point.
Indeed, as our negative result shows, strongly degenerate problems ($\alpha \in [1,2)$) fail to be null controllable otherwise. On the other hand, for weakly degenerate problems,  Carleman estimates  do not seem to lead to null controllability results with the same generality as we obtain in this paper for problem \eqref{eq-control} (see Proposition \ref{prop-NCstr-w}).

Our method is based on a careful analysis of the spectral problem associated with \eqref{eq-control-typ}. As is known from the work of Kamke \cite{Kamke}, the eigenvalues of problem \eqref{eq-control} are related to the zeros of Bessel functions. Indeed this fact has been recently used in a clever way for a boundary control problem by Gueye \cite{Gueye}. When degeneracy occurs inside, the problem looks like a simplified version of the one studied in Zhang-Zuazua \cite{XuZu} in the case of a 1-D fluid-structure model: we solve the problem on both sides of the degeneracy, and we study the transmission conditions.
Once the spectral problem is solved, negative results come quite immediately when $\alpha \in [1,2)$. For the positive part, we combine the moment method with general results obtained in \cite{cost-weak} concerning the existence of biorthogonal families under general gap conditions on the square roots of the eigenvalues. Then,  we complete the analysis with some $L^2$ lower bounds for the eigenfunctions on the control region. 

The theoretical results of our foundings are completed with a final numerical work which concludes the paper (see section \ref{sec-num}). Starting from the pioneering works of J.-L. Lions (see, e.g., \cite{Lions}), numerical approximation of controllability problems for parabolic equations has become an established matter. Among the rich literature, we quote here the basic results in \cite{LopZua}, devoted to the nondegenerate heat equation with boundary control, along with the more recent studies in \cite{BHLR1} and the general framework provided in \cite{LabTre}.

Most of the literature is concerned with semi-discretized problems and aims at constructing an approximation of the stabilizing control. We will rather use here a fully discrete approximation, and avoid the problem of convergence of approximate stabilizing controls to a limit solution, which is known to be a very ill-conditioned problem (see \cite{BHLR2} for a study of the full discretization, as well as a review of the relevant literature). Therefore, the last section should be understood as a numerical illustration of the theoretical results, while a rigorous study of numerical approximations for this problem will be postponed to future works.


\subsection{Plan of the paper} \hfill

\begin{itemize}

\item In section \ref{s-results-1D}, we state our results, distinguishing the weakly degenerate case ($\alpha \in [0,1)$) from the strongly degenerate one ($\alpha \in [1,2)$) (the eigenvalues and the eigenfunctions are not the same according to the case we are considering).

\item In section \ref{sec-wp}, we prove the well-posedness of the problem in both cases.

\item In section \ref{sec-strong}, we solve the spectral problem in the strongly degenerate case ($\alpha \in [1,2)$).

\item In section \ref{sec-NCstrong}, we prove our negative and positive null controllability result for the strongly degenerate case.

\item In section \ref{sec-SLweak}, we solve the spectral problem associated with the weakly degenerate case ($\alpha \in [0,1)$).

\item In section \ref{sec-NCweak}, we prove our positive null controllability result for the weakly degenerate case by combining gap estimates for the eigenvalues with a classical property for the $L^2$ norm of the eigenfunctions in the control region  (in the spirit of Lagnese \cite{Lagnese} for a nondegenerate setting).

\item In section \ref{sec-num}, we provide numerical examples illustrating the above positive and negative controllability results.
\end{itemize}


\section{Main results}
\label{s-results-1D}


\subsection{The strongly degenerate case: $\alpha \in [1,2)$} \hfill

\subsubsection{Functional setting and well-posedness when $\alpha \in [1,2)$}  \hfill

We consider
$$ X = L^2 (-1,1) $$
endowed with the natural scalar product
$$ \forall f,g \in   X, \quad \langle f , g \rangle = \int _{-1} ^1 fg.$$

For $1 \leq \alpha < 2$, we define
\begin{multline}
\label{*H^1_a-2}
 H^1_{\alpha} (-1,1):=\Bigl \{ u \in L^2 (-1,1) \ \mid \  u \text{ locally absolutely continuous in } (0,1]
\\ 
\text{ and in } [-1,0), \int _{-1} ^1 \vert x \vert ^\alpha u_x ^2  \, dx < \infty
\text{ and } u(-1)=0 =u(1) \Bigr \} .
 \qquad \qquad
\end{multline}
$H^1_{\alpha} (-1,1)$ is endowed with the natural scalar product
$$ \forall f,g \in  H^1_{\alpha} (-1,1), \quad (f,g) = 
\int _{-1} ^1  \bigl( \vert x \vert ^\alpha f'(x) g'(x) + f(x) g(x) \bigr) \, dx .$$

Next, consider 
\begin{equation*}
H^2_{\alpha} (-1,1):=   \Bigl \{ u \in H^1_{\alpha }(-1,1)  \ \mid \int _{-1} ^1 \vert (\vert x \vert ^\alpha u'(x) ) ' \vert ^2 \, dx < \infty  \Bigr \} ,
\end{equation*}
and the operator $A:D(A)\subset X \to X$ will be  defined by
\begin{equation*}
D(A) :=  H^2_{\alpha} (-1,1) \quad \text{ and } \quad
\forall u \in D(A), \quad    Au:= (\vert x \vert ^\alpha  u_x)_x .
\end{equation*}

Then the following results hold:

\begin{Proposition}
\label{Prop-A-s}
Given $\alpha \in [1,2)$, we have:

a) $H^1_{\alpha} (-1,1)$ is a Hilbert space.

b) $A: D(A) \subset X \to X$ is a self-adjoint negative operator with dense domain.

\end{Proposition}

Hence, $A$ is the infinitesimal generator of an analytic semigroup of contractions $e^{tA}$ on $X$. 
Given a source term $h$ in  $L^2((-1,1)\times (0,T))$ and an initial condition $v_0 \in  X$, consider the problem
\begin{equation}
\label{eq-v-h-s}
\begin{cases}
v_t - (\vert x \vert ^\alpha v_x)_x = h(x,t), \\
v (-1,t)=0 = v(1,t) , \\
v(x,0)=v_0(x)  .
\end{cases}
\end{equation}
 The function $v \in \mathcal C ^0 ([0,T]; X)  \cap L^2(0,T; H^1_{\alpha} (-1,1))$
given by the variation of constant formula
$$ v(\cdot ,t) = e^{tA} v_0 + \int _0 ^t e^{(t-s)A} h(\cdot, s) \, ds
$$ 
is called the mild solution of \eqref{eq-v-h-s}. We say that a function 
$$v \in
\mathcal C ^0 ([0,T]; H^1_{\alpha} (-1,1))  \cap   H^1(0,T; L^2(-1,1)) \cap    L^2 (0,T; D(A)) $$
is a strict solution of \eqref{eq-v-h-s} if $v$ satisfies $v_t - (\vert x \vert ^\alpha v_x)_x = h(x,t)$ almost everywhere
in $(-1,1 )\times (0,T)$, and the initial and boundary conditions  for all $t\in [0,T]$ and all $x\in [-1,1]$.  And one can prove the existence and the uniqueness of the strict solution.

\begin{Proposition}
\label{prop-wp}
If $v_0 \in H^1 _{\alpha}(-1,1)$, then the mild solution of \eqref{eq-v-h-s} is the unique strict solution of \eqref{eq-v-h-s}.
\end{Proposition}


\subsubsection{The eigenvalue problem when $\alpha \in [1,2)$}  \hfill

The knowledge of the eigenvalues and associated eigenfunctions of the degenerate diffusion operator $ u \mapsto - (\vert x \vert^\alpha u')'$, i.e., the nontrivial solutions $(\lambda, \Phi)$ of
\begin{equation}
\label{*pbm-vp}
\begin{cases}
 - (\vert x \vert ^\alpha \Phi '(x))' =\lambda \Phi (x)  \qquad x\in (-1,1),\\
\Phi(-1)=0 =\Phi(1),
\end{cases}
\end{equation}
 will be essential for our purposes.
It is well-known that Bessel functions play an important role in this problem, see, e.g., Kamke \cite{Kamke}. For $\alpha \in [1,2)$, let
$$ \nu _\alpha := \frac{ \vert \alpha -1 \vert }{2-\alpha} = \frac{ \alpha -1 }{2-\alpha}, 
\qquad \kappa _\alpha:= \frac{2-\alpha}{2}.$$
Given $\nu$, we denote by $J_\nu$ the Bessel function of the first kind of order $\nu$ (see sections \ref{Bessel-intro}-\ref{Bessel-zeros}) and by $j_{\nu,1}< j_{\nu,2} < \dots < j_{\nu,n} <\dots$ the positive zeros of $J_\nu$.

When $\alpha \in [1,2)$, we have the following description of the spectrum of the associated operator:

\begin{Proposition}
\label{*prop-vp}
The admissible eigenvalues $\lambda$ for problem \eqref{*pbm-vp} are given by 
\begin{equation}
\label{*vp}
\forall n \geq 1, \qquad \lambda_{\alpha, n} = \kappa _\alpha ^2 j_{\nu _\alpha ,n}^2  .
\end{equation}
The associated eigenspace is of dimension 2, and an orthonormal basis (in  $L^2(-1,1))$) is given by the following eigenfunctions
\begin{equation}
\label{*fp-d}
\tilde \Phi ^{(r)} _{\alpha, n}(x) :=
 \begin{cases}
\frac{\sqrt{2 \kappa _\alpha }}{ \vert J'_{\nu_\alpha} (j_{\nu_\alpha,n} ) \vert} x^{(1-\alpha)/2} J_{\nu _\alpha} (j_{\nu_\alpha,n} x ^{\kappa_\alpha}) \quad & \text{ if } x \in (0,1)\\
0 \quad & \text{ if } x \in (-1,0) 
\end{cases} ,
\end{equation}
and
\begin{equation}
\label{*fp-g}
\tilde  \Phi ^{(l)} _{\alpha, n}(x) :=
 \begin{cases}
0 \quad & \text{ if } x \in (0,1)\\
\frac{\sqrt{2 \kappa _\alpha }}{ \vert J'_{\nu_\alpha} (j_{\nu_\alpha,n} ) \vert} \vert x \vert ^{(1-\alpha)/2} J_{\nu _\alpha} (j_{\nu_\alpha,n} \vert x \vert  ^{\kappa_\alpha})\quad & \text{ if } x \in (-1,0) 
\end{cases} .
\end{equation}
Moreover $\{ \tilde \Phi ^{(r)} _{\alpha, n}, \tilde \Phi ^{(l)} _{\alpha, n}, n\geq 1 \}$
forms an orthonormal basis of $L^2(-1,1)$.
\end{Proposition}


\subsubsection{Null controllability when $\alpha \in [1,2)$} \hfill

The following controllability result is a direct consequence of the above proposition.

\begin{Proposition}
\label{prop-NCstr}
Assume that $\alpha \in [1,2)$ and let $0 < a < b < 1$. Then
null controllability fails, and the initial conditions that can be steered to $0$ in time $T$ are exactly those which are supported in $[0,1)$.
\end{Proposition}


\subsection{The weakly degenerate case: $\alpha \in [0,1)$} \hfill


\subsubsection{Functional setting and well-posedness when $\alpha \in [0,1)$}  \hfill

For  $0 \leq \alpha < 1$, we consider
\begin{multline}
\label{*H^1_a-w}
 H^1_{\alpha} (-1,1):=\Bigl \{ u \in L^2 (-1,1) \ \mid \  u \text{ absolutely continuous in } [-1,1],
\\ 
\int _{-1} ^1 \vert x \vert ^{\alpha} u_x ^2 \, dx < \infty 
\text{ and } u(-1)=0 =u(1) \Bigr \} .
 \qquad \qquad
\end{multline}
$H^1_{\alpha} (-1,1)$ is endowed with the natural scalar product
$$ \forall f,g \in  H^1_{\alpha} (-1,1), \quad (f,g) = 
\int _{-1} ^1 \bigl( \vert x \vert ^\alpha f'(x) g'(x) + f(x) g(x) \bigr) \, dx .$$

Next, consider 
\begin{equation*}
H^2_{\alpha} (-1,1):=   \Bigl \{ u \in H^1_{\alpha }(-1,1)  \ \mid \int _{-1} ^1 \vert (\vert x \vert ^\alpha u'(x) ) ' \vert ^2 \, dx < \infty  \Bigr \} ,
\end{equation*}
and the operator $A:D(A)\subset X \to X$ will be  defined in a similar way by
\begin{equation*}
D(A) :=  H^2_{\alpha} (-1,1) \quad \text{ and } \quad
\forall u \in D(A), \quad    Au:= (\vert x \vert ^\alpha  u_x)_x .
\end{equation*}

Then the following results hold:

\begin{Proposition}
\label{Prop-A-w}

Given $\alpha \in [0,1)$, we have the following:

a) $H^1_{\alpha} (-1,1)$ is a Hilbert space;

b) $A: D(A) \subset X \to X$ is a self-adjoint negative operator with dense domain.

\end{Proposition}

Hence, once again, $A$ is the infinitesimal generator of an analytic semigroup of contractions $e^{tA}$ on $X$. 
Given a source term $h$ in  $L^2((-1,1)\times (0,T))$ and an initial condition $v_0 \in  X$, consider the problem
\begin{equation}
\label{eq-v-h-w}
\begin{cases}
v_t - (\vert x \vert ^\alpha v_x)_x = h(x,t), \\
v (-1,t)=0 = v(1,t) , \\
v(x,0)=v_0(x)  .
\end{cases}
\end{equation}
 The function $v \in \mathcal C ^0 ([0,T]; X)  \cap L^2(0,T; H^1_{\alpha} (-1,1))$
given by the variation of constant formula
$$ v(\cdot ,t) = e^{tA} v_0 + \int _0 ^t e^{(t-s)A} h(\cdot, s) \, ds
$$ 
is called the mild solution of \eqref{eq-v-h-w}. We say that a function 
$$v \in
\mathcal C ^0 ([0,T]; H^1_{\alpha} (-1,1))  \cap   H^1(0,T; L^2(-1,1)) \cap    L^2 (0,T; D(A)) $$
is a strict solution of \eqref{eq-v-h-w} if $v$ satisfies $v_t - (\vert x \vert ^\alpha v_x)_x = h(x,t)$ almost everywhere
in $(-1,1 )\times (0,T)$, and the initial and boundary conditions  are fulfilled for all $t\in [0,T]$ and all $x\in [-1,1]$. And once again one can prove the existence and the uniqueness of the strict solution.

\begin{Proposition}
If $v_0 \in H^1 _{\alpha}(-1,1)$, then the mild solution of \eqref{eq-v-h-w} is the unique strict solution of \eqref{eq-v-h-w}.
\end{Proposition}


\subsubsection{Eigenvalues and eigenfunctions when $\alpha \in [0,1)$} \hfill

Once again, the knowledge of the eigenvalues and associated eigenfunctions of the degenerate diffusion operator $ u \mapsto - (\vert x \vert^\alpha u')'$, i.e. the solutions $(\lambda, \Phi)$ of
\begin{equation}
\label{*pbm-vp-w}
\begin{cases}
 - (\vert x \vert ^\alpha \Phi '(x))' =\lambda \Phi (x)  \qquad x\in (-1,1),\\
\Phi(-1)=0 =\Phi(1)
\end{cases}
\end{equation}
 will be essential for our purposes.
When $\alpha \in [0,1)$, let
$$ \nu _\alpha := \frac{ \vert \alpha -1 \vert }{2-\alpha} = \frac{ 1-\alpha }{2-\alpha}, 
\qquad \kappa _\alpha:= \frac{2-\alpha}{2}.$$
Now we will need the zeros of the Bessel function $J_{\nu_\alpha}$, and also the zeros of the Bessel function of negative order $J_{-\nu_\alpha}$ (see subsection \ref{Bessel-zeros}).

We prove the following description for \eqref{*pbm-vp-w}:

\begin{Proposition}
\label{*prop-vp-w}
When $\alpha \in [0,1)$, we have exactly two sub-families of eigenvalues and associated eigenfunctions for problem \eqref{*pbm-vp-w}, that is:

\begin{itemize}
\item the eigenvalues of the form
$\kappa _\alpha ^2 j_{\nu_\alpha,n}^2$, associated with the odd functions
\begin{equation}
\label{Phi-alpha-n-w-i}
\Phi _{\alpha,n} ^{(o)} (x)=
\begin{cases}
 x^{\frac{1-\alpha}{2}}  J_{\nu _\alpha} (j_{\nu_\alpha,n} x^{\kappa _\alpha}) \quad & \text{ if } x \in (0,1) 
 \\
 -   \vert x \vert ^{\frac{1-\alpha}{2}}  J_{\nu _\alpha} (j_{\nu_\alpha,n} \vert x \vert ^{\kappa _\alpha}) \quad & \text{ if } x \in (-1,0) 
 \end{cases},
\end{equation}

\item the eigenvalues of the form
$\kappa _\alpha ^2 j_{-\nu_\alpha,n}^2$, associated with the even functions
\begin{equation}
\label{Phi-alpha-n-w-p}
\Phi _{\alpha,n} ^{(e)} (x) = 
\begin{cases}
 x ^{\frac{1-\alpha}{2}} J_{-\nu _\alpha} (j_{-\nu_\alpha,n} x^{\kappa _\alpha}) \quad & \text{ if } x \in (0,1) 
 \\
\vert x \vert ^{\frac{1-\alpha}{2}} J_{-\nu _\alpha} (j_{-\nu_\alpha,n} \vert x \vert ^{\kappa _\alpha}) \quad & \text{ if } x \in (-1,0) 
 \end{cases} .
\end{equation}
\end{itemize}
Moreover, the family $\{ \Phi _{\alpha,n} ^{(o)}, \Phi _{\alpha,n} ^{(e)}, n\geq 1 \}$ forms an orthogonal basis of $L^2(-1,1)$.
\end{Proposition}


\subsubsection{Null controllability when $\alpha \in [0,1)$} \hfill

\begin{Proposition} 
\label{prop-NCstr-w}
Assume that $\alpha \in [0,1)$ and that $0<a<b<1$. Then null controllability holds: given $u_0 \in L^2 (-1,1)$, there exists a control $h$ that drives the solution $u$ to $0$ in time $T$.
\end{Proposition}


\subsubsection{Blow-up of the control cost as $\alpha \to 1^-$} \hfill

Given $\alpha \in [0,1)$, $T>0$ and $u_0 \in L^2 (-1,1)$, consider 
$\mathcal{U}^{ad} (\alpha,T; u_0) $ the set of admissible controls:
$$ \mathcal{U}^{ad} (\alpha,T; u_0) := \Bigl \{h \in L^2 ((a,b)\times (0,T))\  | \ u^{(h)}(T)=0 \Bigr \} .$$
Since null controllability holds if and only if $\alpha <1$, it is natural to expect that the null controllability cost 
\begin{equation}
\label{def-cost}
C_{NC} (\alpha, T) := \sup _{\Vert u_0 \Vert _{L^2 (-1,1)}} \inf _{h \in \mathcal{U}^{ad} (\alpha,T; u_0) } \Vert h \Vert _{L^2((a,b)\times (0,T))} 
\end{equation}
blows up when $\alpha \to 1^-$. This is the object of the following result:

\begin{Theorem}
\label{prop-cost}
a) Estimate from above: there exists some $C>0$ independent of $\alpha \in [0,1)$ and of $T>0$ such that
\begin{equation}
\label{1310-12}
C_{NC} (\alpha,T) \leq \frac{C}{(1-\alpha)^2} e^{-T/C} e^{C/T} ;
\end{equation}

b) Estimate from below: there exists some $C'>0$ independent of $\alpha \in [0,1)$ and of $T>0$ such that
\begin{equation}
\label{1310-11}
C_{NC} (\alpha,T)  \geq \frac{C'}{(1-\alpha) \sqrt{T}} e^{-T/C'} .
\end{equation}
\end{Theorem}

Note that Theorem \ref{prop-cost} proves that the null controllability cost blows up when $\alpha \to 1^-$, and when $T\to 0 ^+$. Moreover:
\begin{itemize}
\item with respect to $\alpha$: we have a good estimate of the behavior when $\alpha \to 1^-$,
but the upper and lower estimates are not of the same order;
\item with respect to $T$: in \cite{CMV-cost-loc} we prove that the blow-up of the null controllability cost is of the order $e^{C/T}$; here we only obtain a weak blow-up estimate, of the order $1/\sqrt{T}$; we conjecture that the blow-up rate is of the form $e^{C/T}$.
\end{itemize}
It would be interesting to have better blow-up estimates, with respect to $\alpha \to 1^-$ and  $T\to 0^+$.


\section{Well-posedness: proof of Propositions \ref{Prop-A-s} and \ref{Prop-A-w}}
\label{sec-wp}


\subsection{The strongly degenerate case} \hfill

\subsubsection{Proof of Proposition \ref{Prop-A-s}, part a)} \hfill

Let us verify that $H^1 _\alpha (-1,1)$, defined in \eqref{*H^1_a-2}, is complete.
Take $(f_j)_j$ a Cauchy sequence in $H^1 _\alpha (-1,1)$. Then
$(f_j)_j$ is a Cauchy sequence in $L^2 _\alpha (-1,1)$ and in 
$W^{1,2}(\varepsilon, 1)$ and $W^{1,2}(-1,-\varepsilon)$ for all fixed positive $\varepsilon$.
Thus $(f_j)_j$ converges to some limit $f$, which has to be in $L^2 _\alpha (-1,1)$ and in $W^{1,2}(\varepsilon, 1)$ and $W^{1,2}(-1,-\varepsilon)$ for all $\varepsilon \in (0,1)$. Hence $f$ is locally absolutely continuous on $(0,1]$ and $[-1,0)$. Moreover, 
$$ \int _{\varepsilon} ^1 \vert x \vert ^\alpha f'(x) ^2 \, dx 
\leq \liminf _j \int _{\varepsilon} ^1 \vert x \vert ^\alpha f'_j (x) ^2 \, dx,$$
which implies that $f\in H^1 _\alpha (-1,1)$. 
Finally, $f_j \to f$ in 
$H^1 _\alpha (-1,1)$: indeed, since $(f_j)_j$ is a Cauchy sequence in $H^1 _\alpha (-1,1)$, we have
$$ \forall \varepsilon >0, \exists n_0: \quad  i,j \geq n_0 \quad \implies \quad
\Vert f_i - f_j \Vert _{H^1 _\alpha (-1,1)} \leq \varepsilon
.$$
If $(f_j)_j$ does not converge to $f$ in $H^1 _\alpha (-1,1)$, then there exists $\varepsilon _0 >0$ and a subsequence $(f_{j'})_{j'}$ such that
$$ \Vert f_{j'} - f \Vert _{H^1 _\alpha (-1,1)} \geq \varepsilon _0 .$$
Using the Cantor diagonal process, one can find a subsequence $(f_{j''})_{j''}$
such that 
$$ f' _{j''} \to f' \quad a.e. (-1,1) .$$
Then, with the Cauchy criterion applied with $\varepsilon < \varepsilon _0$, we have 
$$ \Vert f_{j''} - f_{i''} \Vert _{H^1 _\alpha (-1,1)} \leq \varepsilon $$
if $i'', j'' $ are large enough. So, using Fatou's lemma, we obtain
$$  \Vert f_{j''} - f \Vert _{H^1 _\alpha (-1,1)} \leq \varepsilon < \varepsilon _0, $$
which is a contradiction. Hence $(f_j)_j$ converges to $f$ in $H^1 _\alpha (-1,1)$. This concludes the proof of Proposition \ref{Prop-A-s}, part a). \qed


\subsubsection{Integration by parts} \hfill

Let us prove the following integration by parts formula:

\begin{Lemma}
\label{lem-IPP-s}
\begin{equation}
\label{IPP-s}
\forall f,g \in H^2 _\alpha (-1,1), \quad
\int _{-1} ^1 (\vert x \vert ^\alpha f'(x) )' \, g(x) \, dx 
= - \int _{-1} ^1 \vert x \vert ^\alpha \, f'(x) \, g '(x) \, dx .
\end{equation}
\end{Lemma}

\noindent {\it Proof of Lemma \ref{lem-IPP-s}.} 
If $f \in H^2 _\alpha (-1,1)$, then 
$$ F(x):= \vert x \vert ^\alpha f'(x) \in H^1 (-1,1) .$$
Take $g \in H^2 _\alpha (-1,1)$, and $\varepsilon \in (0,1)$.
Decompose
$$ \int _{-1} ^1 F'(x) g(x) \, dx  = 
\int _{-1} ^{-\varepsilon} F'(x) g(x) \, dx
+ \int _{-\varepsilon} ^\varepsilon F'(x) g(x) \, dx
+ \int _{\varepsilon} ^1  F'(x) g(x) \, dx .$$
Then, since $g \in H^2 _\alpha (-1,1) \subset H^1 _0 (-1,1)$, the usual integration by parts formula gives
$$ 
\int _{-1} ^{-\varepsilon} F'(x) g(x) \, dx
= [F(x) g(x) ] _{-1} ^{-\varepsilon} - \int _{-1} ^{-\varepsilon} F(x) g'(x) \, dx ,$$
and 
$$ 
\int _{\varepsilon} ^1 F'(x) g(x) \, dx
= [F(x) g(x) ] _{\varepsilon} ^1 - \int _{\varepsilon} ^1 F(x) g'(x) \, dx .$$
Now, since  $F$ and $g'$ belong to $L^2 (-1,1)$, we have
$$ \int _{-1} ^{-\varepsilon} F(x) g'(x) \, dx +
\int _{\varepsilon} ^1 F(x) g'(x) \, dx \to \int _{-1} ^1 F(x) g'(x) \, dx
\quad \text{ as } \varepsilon \to 0 ,$$
and since $F'$ and $g$ belong to $L^2 (-1,1)$, we have
$$  \int _{-\varepsilon} ^\varepsilon F'(x) g(x) \, dx \to 0 
\quad \text{ as } \varepsilon \to 0 .$$
It remains to study the boundary terms:
first, because of Dirichlet boundary conditions at $\pm1$, we have
$$ [F(x) g(x) ] _{-1} ^{-\varepsilon} = F(-\varepsilon) g(-\varepsilon)
\quad \text{ and } \quad [F(x) g(x) ] _{\varepsilon} ^1
= - F(\varepsilon) g(\varepsilon) .$$
Now we note that 
\begin{multline*}
 \forall x \in (0,1), \quad (F(x)g(x))' = F'(x) g(x) + F(x) g'(x) 
\\
= F'(x) g(x) + (x^{\alpha /2} f'(x)) (x^{\alpha /2} g'(x)),
\end{multline*}
and since $F'$, $g$, $x^{\alpha /2} f'$, $x^{\alpha /2} g'$ belong to $L^2(0,1)$, 
we see that $(Fg)' \in L^1(0,1)$. Hence $Fg$ is absolutely continuous on $(0,1]$. Therefore it has a limit as $x\to 0$, there exists $L$ such that
$$ F(x) g(x) \to L \text{ as } x \to 0 ^+.$$
In fact, we claim that $L=0$. Indeed:
\begin{itemize}
\item first, the function $x\mapsto x^\alpha f'(x)$ belongs to $H^1(-1,1)$, hence it has a limit as $x\to 0^+$:
$$ x^\alpha f'(x) \to \ell \quad \text{ as } x \to 0 ^+ ;$$
\item if $\ell \neq 0$,  
$$ x^{\alpha /2} f'(x) \sim \frac{\ell}{x^{\alpha /2}} \quad \text{ as } x \to 0 ^+ ;$$
but since $\alpha \geq 1$, we have that $\frac{\ell}{x^{\alpha /2}} \notin L^2 (0,1)$; so $\ell =0$;
\item then 
$$ \forall x \in (0,1), \quad x^\alpha f'(x) = \int _0 ^x (s^\alpha f'(s))' \, ds $$
and using the Cauchy-Schwarz inequality, one has
$$ \forall x\in (0,1), \quad \vert x ^\alpha f'(x) \vert \leq C \sqrt{ x } ;$$
\item finally, 
$$ \forall x\in (0,1), \quad \vert x ^\alpha f'(x) g(x) \vert \leq C \sqrt{ x } \vert g(x) \vert,$$
hence
$$\forall x\in (0,1), \quad \vert F(x) g(x) \vert \leq C \sqrt{ x } \vert g(x) \vert;$$
if $L\neq 0$, then if $x$ is sufficiently close to $0$ we have
$$ \vert g(x) \vert \geq \frac{CL}{2 \sqrt{x} },$$
which is in contradiction with $g\in L^2(-1,1)$. Hence $L=0$.
\end{itemize}
This implies that
$$ [F(x) g(x) ] _{\varepsilon} ^1
=  - F(\varepsilon) g(\varepsilon) 
\to 0 \quad \text{ as } \varepsilon \to 0^+ .$$
Same property holds on $[-1,0)$, hence 
$$ [F(x) g(x) ] _{-1} ^{\varepsilon}  \to 0 \quad \text{ as } \varepsilon \to 0^+ .$$
This concludes the proof of Lemma \ref{lem-IPP-s}. \qed


\subsubsection{Proof of Proposition \ref{Prop-A-s}, part b) } \hfill

First we note that it is clear that $D(A)$ is dense, since it contains
all the functions of class $C^\infty$, compactly supported in $(-1,1)$.

Next we derive from Lemma \ref{lem-IPP-s} that 
$$ \forall f \in D(A), \quad 
\langle Af , f \rangle = \int _{-1} ^1 ( \vert x \vert ^\alpha f'(x))' f(x) \, dx 
= - \int _{-1} ^1  \vert x \vert ^\alpha f'(x) ^2 \, dx \leq 0 ,$$
hence $A$ is dissipative.

In order to show that $A$ is symmetric, we apply Lemma \ref{lem-IPP-s} twice to obtain that
\begin{multline*}
\forall f,g \in D(A), \quad 
\langle Af , g \rangle 
= \int _{-1} ^1 ( \vert x \vert ^\alpha f'(x))' g(x) \, dx 
= - \int _{-1} ^1  \vert x \vert ^\alpha f'(x) g'(x) \, dx 
\\
= - \int _{-1} ^1  (\vert x \vert ^\alpha g'(x)) f'(x) \, dx 
= \int _{-1} ^1 ( \vert x \vert ^\alpha g'(x))' f(x) \, dx 
= \langle f , Ag \rangle .
\end{multline*}

Finally, we check that $I-A$ is surjective. Let $f\in X$. Then, by Riesz theorem,
there exists one and only one $u\in H^1 _\alpha (-1,1)$ such that 
$$\forall v \in H^1 _\alpha (-1,1), \quad  
\int _{-1} ^1 \bigl( uv+ \vert x \vert ^\alpha u' v' \bigr) = \int _{-1} ^1 fv .$$
In particular the above relation holds true for all $v$ of class $C^\infty$, compactly supported in $(-1,1)$. Hence, $x\mapsto \vert x \vert ^\alpha u'$
has a weak derivative given by
$$ \Bigl( \vert x \vert ^\alpha u' \Bigr) '= - (f-u) ;$$
since $f-u \in L^2 (-1,1)$, we obtain that $( \vert x \vert ^\alpha u')'
\in L^2 (-1,1)$. Hence $u\in H^2 _\alpha (-1,1)$. Then $(I-A)u=f$. So the operator $I-A$ is surjective. This concludes the proof of Proposition \ref{Prop-A-s}, part b). \qed


\subsection{The weakly degenerate case} \hfill

\subsubsection{Proof of Proposition \ref{Prop-A-w}, part a)} \hfill

We have to check that $H^1 _\alpha (-1,1)$  defined by \eqref{*H^1_a-w} is complete. Take a Cauchy sequence $(f_j)_j$ in $H^1 _\alpha (-1,1)$.
Then $(f_j)_j$ is also a Cauchy sequence in $L^2 (-1,1)$,
$(\vert x \vert ^{\alpha /2} f'_j)_j$ is a Cauchy sequence in $L^2 (-1,1)$,
and since
$$ \forall f \in H^1 _\alpha (-1,1), \quad 
\int _{-1} ^1 \vert f'(x) \vert \, dx \leq 
\Bigl(\int _{-1} ^1 \frac{1}{\vert x \vert ^{\alpha }} \, dx \Bigr) ^{1/2} 
\Bigl( \int _{-1} ^1 \vert x \vert ^{\alpha } \vert f'(x) \vert ^2 \, dx\Bigr) ^{1/2} ,$$
we see that  $(f_j)_j$ is also a Cauchy sequence in $W^{1,1} (-1,1)$.
Hence there exists $f\in L^2 (-1,1) \cap W^{1,1} (-1,1)$ such that
$$ f_j \to f \text{ in } L^2(-1,1), \quad f'_j \to f' \text{ in } L^1(-1,1).$$
It remains to check that $\vert x \vert ^{\alpha /2} f' \in L^2 (-1,1)$.
This comes from the following remark:
up to a sequence, we have
$$ f'_j \to f' \quad \text{ a.e. on } (-1,1) ,$$
hence
$$ \vert x \vert ^{\alpha } \vert f'_j (x) \vert ^2
\to \vert x \vert ^{\alpha } \vert f'(x) \vert ^2 \quad \text{ a.e. on } (-1,1) ,$$
and the Fatou Lemma implies that
$$ \int _{-1} ^1 \vert x \vert ^{\alpha } \vert f'(x) \vert ^2 \, dx
\leq \liminf_j \int _{-1} ^1 \vert x \vert ^{\alpha } \vert f'_j(x) \vert ^2 \, dx ,$$
hence $\vert x \vert ^{\alpha /2} f' \in L^2 (-1,1)$. 
Finally, $f_j \to f$ in $H^1 _\alpha (-1,1)$, using the same kind of arguments as in 
Proposition \ref{Prop-A-s}, part a). This concludes the proof of Proposition \ref{Prop-A-w}, part a).
\qed


\subsubsection{Integration by parts} \hfill

Let us prove the following integration by parts formula:

\begin{Lemma}
\label{lem-IPP-w}
\begin{equation}
\label{IPP-w}
\forall f,g \in H^2 _\alpha (-1,1), \quad
\int _{-1} ^1 (\vert x \vert ^\alpha f'(x) )' \, g(x) \, dx 
= - \int _{-1} ^1 \vert x \vert ^\alpha \, f'(x) \, g '(x) \, dx .
\end{equation}
\end{Lemma}

\noindent {\it Proof of Lemma \ref{lem-IPP-w}.} 
If $f \in H^2 _\alpha (-1,1)$, then 
$$ F(x):= \vert x \vert ^\alpha f'(x) \in H^1 (-1,1) .$$
Take $g \in H^2 _\alpha (-1,1)$, and $\varepsilon \in (0,1)$.
Decompose
$$ \int _{-1} ^1 F'(x) g(x) \, dx  = 
\int _{-1} ^{-\varepsilon} F'(x) g(x) \, dx
+ \int _{-\varepsilon} ^\varepsilon F'(x) g(x) \, dx
+ \int _{\varepsilon} ^1  F'(x) g(x) \, dx .$$
Then, since $g \in H^2 _\alpha (-1,1) \subset H^1 _0 (-1,1)$, the usual integration by parts formula gives
$$ 
\int _{-1} ^{-\varepsilon} F'(x) g(x) \, dx
= [F(x) g(x) ] _{-1} ^{-\varepsilon} - \int _{-1} ^{-\varepsilon} F(x) g'(x) \, dx ,$$
and 
$$ 
\int _{\varepsilon} ^1 F'(x) g(x) \, dx
= [F(x) g(x) ] _{\varepsilon} ^1 - \int _{\varepsilon} ^1 F(x) g'(x) \, dx .$$
Now, since  $F$ and $g'$ belong to $L^2 (-1,1)$, we have
$$ \int _{-1} ^{-\varepsilon} F(x) g'(x) \, dx +
\int _{\varepsilon} ^1 F(x) g'(x) \, dx \to \int _{-1} ^1 F(x) g'(x) \, dx
\quad \text{ as } \varepsilon \to 0 ,$$
and since $F'$ and $g$ belong to $L^2 (-1,1)$, we have
$$  \int _{-\varepsilon} ^\varepsilon F'(x) g(x) \, dx \to 0 
\quad \text{ as } \varepsilon \to 0 .$$
It remains to study the boundary terms:
first, because of Dirichlet boundary conditions at $\pm1$, we have
$$ [F(x) g(x) ] _{-1} ^{-\varepsilon} = F(-\varepsilon) g(-\varepsilon)
\quad \text{ and } \quad [F(x) g(x) ] _{\varepsilon} ^1
= - F(\varepsilon) g(\varepsilon) .$$
Now we note that $G$ and $g$ are absolutely continuous on $[-1,1]$, hence
$Fg$ is also absolutely continuous on $[-1,1]$, hence, and therefore it has a limit as $x\to 0$: there exists $L$ such that
$$ F(x) g(x) \to L \text{ as } x \to 0 ,$$
which implies that
$$ [F(x) g(x) ] _{-1} ^{-\varepsilon} 
+ [F(x) g(x) ] _{\varepsilon} ^
= F(-\varepsilon) g(-\varepsilon) - F(\varepsilon) g(\varepsilon)
\to 0 \text{ as } \varepsilon \to  0 .$$


\subsubsection{Proof of Proposition \ref{Prop-A-s}, part b) } \hfill

It is similar to the proof of Proposition \ref{Prop-A-s}, part b),
using the integration by parts formula given by Lemma \ref{lem-IPP-w}. \qed


\section{The Sturm-Liouville problem in the strongly degenerate case} 
\label{sec-strong}

The goal of this section is to prove Proposition \ref{*prop-vp}:
we study the spectral problem \eqref{*pbm-vp} and the properties of the eigenvalues and eigenfunctions when $\alpha \in [1,2)$.

First, one can observe that if $\lambda$ is an eigenvalue, then $\lambda >0$: indeed, multiplying \eqref{*pbm-vp} by $\Phi $ and integrating by parts, then
$$ \lambda \int _{-1} ^1 \Phi ^2 = \int _{-1} ^1 \vert x \vert ^\alpha \Phi _x ^2 ,$$
which implies first $\lambda \geq 0$, and next that $\Phi =0$ if $\lambda =0$.

Now we make the following observation: if $(\lambda,\Phi )$ solves \eqref{*pbm-vp},
then

\begin{equation}
\label{*pbm-vp-d}
\begin{cases}
 - (x ^\alpha \Phi '(x))' =\lambda \Phi (x)  \qquad x\in (0,1),\\
\Phi (1) = 0, 
\end{cases} ,
\end{equation}
and
\begin{equation}
\label{*pbm-vp-g}
\begin{cases}
 - (\vert x \vert ^\alpha \Phi '(x))' =\lambda \Phi (x)  \qquad x\in (-1,0),\\
\Phi (-1)=0
\end{cases} .
\end{equation}
In the following, we study \eqref{*pbm-vp-d} and \eqref{*pbm-vp-g}, and then we will be able to solve \eqref{*pbm-vp}.


\subsection{The study of \eqref{*pbm-vp-d}} \hfill

\subsubsection{The link with the Bessel's equation} \hfill

There is a change a variables that allows one to transform the eigenvalue problem \eqref{*pbm-vp-d} into a
differential Bessel's equation (see in particular Kamke \cite[section  2.162, equation (Ia), p. 440]{Kamke}, Gueye \cite{Gueye}) and \cite{CMV-cost-loc}:
assume that $\Phi$ is a solution of \eqref{*pbm-vp-d} associated to the eigenvalue $\lambda$; then one easily checks that
the function $\Psi$ defined by
\begin{equation}
\label{*eq-lien}
\Phi  (x) =: x^{\frac{1-\alpha}{2}} \Psi \Bigl(\frac{2}{2-\alpha} \sqrt{\lambda} x^{\frac{2-\alpha}{2}} \Bigr)
\end{equation}
is solution of the following boundary problem:
\begin{equation}
\label{*pb-bessel}
\begin{cases}
y^2 \Psi ''(y) + y \Psi  '(y) + (y^2 - (\frac{\alpha -1}{2-\alpha}) ^2) \Psi (y) = 0, \quad y\in (0, \frac{2 }{2-\alpha}\sqrt{\lambda} ), \\
\Psi \Bigl(\frac{2}{2-\alpha} \sqrt{\lambda} \Bigr) = 0 .
\end{cases}
\end{equation}
This is exactly the Bessel's equation of order $\frac{\alpha -1}{2-\alpha}$.
%
%
For reader convenience, we recall here the definitions concerning Bessel's equation and functions together with some useful properties of these functions and of their zeros. {\it Throughout this section, we assume that $\nu \in \mathbb R_+$}.

 
\subsubsection{Bessel's equation and Bessel's functions of order $\nu$} \hfill
\label{Bessel-intro}

The Bessel's functions of order $\nu$ are the solutions of the following differential equation (see \cite[section 3.1, eq. (1), p. 38]{Watson} or \cite[eq (5.1.1), p. 98]{Lebedev}): 
\begin{equation}
\label{*eq-bessel-ordre-nu}
y^2 \Psi ''(y) + y \Psi '(y) +(y^2-\nu^2) \Psi (y)=0, \qquad y\in (0,+\infty).
\end{equation}
The above equation is called {\it Bessel's equation for functions of order $\nu$}. 
Of course the fundamental theory of ordinary differential equations says that the solutions of \eqref{*eq-bessel-ordre-nu} generate a vector space $S_\nu$ of dimension 2. In the following we recall what can be chosen 
as a basis of $S _\nu$.


\subsubsection{Fundamental solutions of Bessel's equation when $\nu \notin \Bbb N$} \hfill
\label{Bessel-base-pas entier}

Assume that $\nu \notin \Bbb N$. When looking for solutions of  \eqref{*eq-bessel-ordre-nu} of the form of series of ascending powers of $x$, one can  construct two series that are solutions:
$$ \sum_{m \geq 0}  \frac{(-1)^m}{m! \ \Gamma (\nu+m+1) } \left( \frac{y}{2}\right) ^{\nu+2m}
\ \text{ and } \ 
\sum_{m \geq 0}  \frac{(-1)^m}{m! \ \Gamma (-\nu+m+1) } \left( \frac{y}{2}\right) ^{-\nu+2m},$$
where $\Gamma$ is the Gamma function
(see \cite[section 3.1, p. 40]{Watson}).
 The first of the two series converges for all values of $y$ and defines the so-called Bessel function of order $\nu$ and {\it of the first kind} which is denoted by $J_\nu$:
\begin{equation}
\label{*def-Jnu}
 J_\nu (y)
:= \sum_{m = 0} ^\infty  \frac{(-1)^m}{m! \ \Gamma (m+\nu+1) } \left( \frac{y}{2}\right) ^{2m+\nu}
=\sum_{m = 0} ^\infty  c_{\nu,m} ^+ y^{2m+\nu} , \qquad y \geq 0,
\end{equation}
 (see \cite[section 3.1, (8), p. 40]{Watson} or \cite[eq. (5.3.2), p. 102]{Lebedev}).
The second series converges for all positive values of $y$ and is $J_{-\nu}$:
\begin{equation}
\label{*def-J-nu}
 J_{-\nu} (y)
:= \sum_{m = 0} ^\infty  \frac{(-1)^m}{m! \ \Gamma (m-\nu+1) } \left( \frac{y}{2}\right) ^{2m-\nu}
=\sum_{m = 0} ^\infty  c_{\nu,m} ^- y^{2m-\nu} , \qquad y > 0 .
\end{equation}

When $\nu \not \in  \mathbb N$, the two functions $J_\nu$ and $J_{-\nu}$ are linearly independent and therefore the pair $(J_\nu,J_{-\nu})$ forms a fundamental system of solutions of  \eqref{*eq-bessel-ordre-nu},
(see  \cite[section 3.12, eq. (2), p. 43]{Watson}).


\subsubsection{Fundamental solutions of Bessel's equation when $\nu =n \in \Bbb N$} \hfill
\label{*sub-bessel-entier}

Assume that $\nu =n \in \Bbb N$. When looking for solutions of  \eqref{*eq-bessel-ordre-nu} of the form of series of ascending powers of $y$, one sees that $J_n$ and $J_{-n}$ are still solutions of \eqref{*eq-bessel-ordre-nu}, 
where $J_n$ is still given by \eqref{*def-Jnu} and $J_{-n}$ is given by \eqref{*def-J-nu}; when $\nu =n \in \Bbb N$, $J_{-n}$ can be written
\begin{equation}
\label{*def-J-n}
J_{-n} (y) = \sum_{m \geq n}  \frac{(-1)^m}{m! \ \Gamma (m-n+1) } \left( \frac{y}{2}\right) ^{-n+2m} .
\end{equation}
However now $J_{-n} (y) = (-1)^n J_{n} (y)$, hence $J_n$ and $J_{-n}$ are linearly dependent, (see \cite[section 3.12, p. 43]{Watson} or \cite[eq. (5.4.10), p. 105]{Lebedev}). The determination of a fundamental system of solutions in this case requires further investigation. In this purpose, one introduces the Bessel's functions of order $\nu$ and {\it of the second kind}:
among the several definitions of Bessel's functions of second order, we recall here the definition by Weber. The Bessel's functions of order $\nu$ and {\it of second kind} are denoted by $Y_\nu$ and defined by (see \cite[section 3.54, eq. (1)-(2), p. 64]{Watson} or \cite[eq. (5.4.5)-(5.4.6), p. 104]{Lebedev}):
\begin{equation}
\label{def-Ynu}
\begin{cases}
\forall \nu \not \in \mathbb N, &\qquad 
Y_\nu(y) := 
\displaystyle{ \frac{J_\nu(y) \cos(\nu \pi)-J_{-\nu} (y)}{\sin(\nu \pi)}},
\\
\forall n \in \mathbb N, &\qquad 
 Y_n(y) := \lim_{\nu\to n} Y_\nu (y).
\end{cases}
\end{equation}
For any $\nu \in \mathbb R_+$, the two functions $J_\nu$ and $Y_\nu$ are always linearly independent, see  \cite[section 3.63, eq. (1), p. 76]{Watson}.
In particular, in the case $\nu=n\in\mathbb N$, the pair $(J_n,Y_n)$ forms a fundamental system of solutions of the Bessel's equation for functions of order $n$.  

In the case $\nu = n \in \mathbb N$,  it will be useful to expand $Y_n$ under the form of a series of ascending powers. This can be done using Hankel's formula, 
see \cite[section 3.52, eq. (3), p. 62]{Watson} 
or \cite[eq. (5.5.3), p. 107]{Lebedev}:
\begin{multline}
\label{*expand-Yn}
\forall n \in \mathbb N^\star, \qquad Y_n(y)=
\frac{2}{\pi} J_n(y) \log \left( \frac{y}{2}\right) 
- \frac{1}{\pi} \sum_{m=0}^{n-1} \frac{(n-m-1)!}{m!} \left(\frac{y}{2}\right) ^{2m-n}
\\
- \frac{1}{\pi} \sum_{m=0}^{+\infty} \frac{(-1)^m}{m!(n+m)!}\left(\frac{y}{2}\right) ^{n+2m} 
\left[ \frac{\Gamma '(m+1)}{\Gamma (m+1)}  + \frac{\Gamma ' (m+n+1)}{\Gamma (m+n+1)}  \right],
\end{multline}
where $\frac{\Gamma '}{\Gamma}$ is the logarithmic derivative of the Gamma function, and satisfies
$\frac{\Gamma '(1)}{\Gamma (1)} = - \gamma$ (here $\gamma$ denotes Euler's constant) and 
$$\frac{\Gamma '(m+1)}{\Gamma(m+1)} = 1+\frac{1}{2} + \ldots \frac{1}{m} -\gamma \ \text{ for all } m \in \mathbb N.$$
In the case $n=0$, the first sum in \eqref{*expand-Yn} should be set equal to zero.


\subsubsection{Zeros of Bessel functions of order $\nu$ of the first kind,
and of order $-\nu$} \hfill
\label{Bessel-zeros}

The function $J_\nu$ has an infinite number of real zeros  which are simple with the possible exception of $y=0$ (\cite[section 15.21, p. 478-479 applied to $C_\nu=J_\nu$]{Watson} or \cite[section 5.13, Theorem 2, p. 127]{Lebedev}). We denote by $(j_{\nu,n})_{n\geq 1}$ the strictly increasing sequence of the positive zeros of $J_{\nu}$:
$$ 0< j_{\nu,1} < j_{\nu,2} < \dots < \ j_{\nu,n} < \dots$$
and we recall that 
$$  j_{\nu,n} \to +\infty \text{ as } n \to +\infty.$$
and the following bounds on the zeros, proved in Lorch and Muldoon \cite{Lorch}:
\begin{equation}
\label{*eq-Lorch2}
\forall \nu \in [0, \frac{1}{2}], \forall n\geq 1, \quad 
\pi (n + \frac{\nu}{2}-\frac{1}{4}) \leq j_{\nu, n} \leq \pi (n + \frac{\nu}{4}-\frac{1}{8}) ,
\end{equation}
and \begin{equation}
\label{*eq-Lorch}
\forall \nu \geq \frac{1}{2}, \forall n\geq 1, \quad 
\pi (n + \frac{\nu}{4}-\frac{1}{8}) \leq j_{\nu, n} \leq\pi (n + \frac{\nu}{2}-\frac{1}{4}) .
\end{equation}

Assume that $\nu \notin \Bbb N$. Then, as an application of the classical Sturm theorem, the function $J_{-\nu}$
has at least one zero between two consecutive zeros of $J_\nu$, and 
at most one, otherwise $J_\nu$ would have at least another zero inside, which is not possible. Hence $J_{-\nu}$ has an increasing  sequence of positive zeros,
interlaced with the ones of $J_\nu$. We will denote them $(j_{-\nu,n})_{n\geq 1}$.




\subsubsection{The solutions of \eqref{*pbm-vp-d} when $\alpha \in [1,2)$}  \hfill

We are going to prove the following 

\begin{Lemma}
\label{lem-struct-d}
Assume that $\Phi$ solves \eqref{*pbm-vp}. Define 
$$ \nu _\alpha := \frac{\alpha -1}{2-\alpha}, \quad \kappa _\alpha := \frac{2-\alpha}{2} .$$
Then there exists some $C_+ \in \Bbb R$ (possibly equal to $0$), and some $m\geq 1$ such that 
\begin{equation}
\label{Phi-struct-d}
\forall x\in (0,1), \quad \Phi (x) = C_+ x^{\frac{1-\alpha}{2}} J_{\nu _\alpha} ( j_{\nu _\alpha, m} x^{\kappa_\alpha}) .
\end{equation}

\end{Lemma}

The proof of Lemma \ref{lem-struct-d} is based on the previous results on Bessel functions, and has to be divided in three, according to
$\nu _\alpha \notin \mathbb N$,
$\nu _\alpha \in \mathbb N^*$, $\nu _\alpha =0$. We study these three cases in the following.


\subsubsection{Proof of Lemma \ref{lem-struct-d} when $\nu _\alpha := \frac{\alpha -1}{2-\alpha} \notin \mathbb N$}  \hfill

Let us assume that $\nu _\alpha \not \in \mathbb N$. Then we have
$$ \Phi = C_+ \Phi _+ + C_- \Phi _-$$
where 
\begin{equation}
\label{*base-phi}
\forall x \in (0,1), \quad 
\begin{cases}
\Phi _+ (x) := x^{\frac{1-\alpha}{2}} J_{\nu _\alpha} (\frac{2}{2-\alpha} \sqrt{\lambda} x^{\frac{2-\alpha}{2}}), \\
\Phi _- (x) := x^{\frac{1-\alpha}{2}} J_{-\nu _\alpha} (\frac{2}{2-\alpha} \sqrt{\lambda} x^{\frac{2-\alpha}{2}}).
\end{cases}
\end{equation}
Then, using the series expansion of $J_{\nu _\alpha}$ and $J_{-\nu _\alpha}$, one obtains
\begin{equation}
\label{*serie-phi-cas2}
\forall x \in (0,1), \quad 
\begin{cases}
\Phi _+ (x) = \sum _{m=0} ^\infty \tilde{c} _{\nu _\alpha ,m} ^+  x ^{(2-\alpha) m}, \\
\Phi _- (x) = \sum _{m=0} ^\infty \tilde{c} _{\nu _\alpha ,m} ^- x ^{1-\alpha + (2-\alpha) m} ,
\end{cases}
\end{equation}
where the coefficients $\tilde{c} _{\nu _\alpha ,m} ^+$ and $\tilde{c} _{\nu _\alpha ,m} ^-$ are defined by
\begin{equation}
\label{*coeffs-phi}
 \tilde{c} _{\nu _\alpha ,m} ^+ := c_{\nu _\alpha ,m} ^+ \Bigl( \frac{2}{2-\alpha} \sqrt{\lambda} \Bigr) ^{2m+\nu _\alpha},
\quad 
\tilde{c} _{\nu _\alpha ,m} ^- := c_{\nu _\alpha ,m} ^- \Bigl( \frac{2}{2-\alpha} \sqrt{\lambda} \Bigr) ^{2m-\nu _\alpha} .
\end{equation}
Then
$$ \Phi _+ (x) \sim _{0} \tilde{c} _{\nu _\alpha ,0} ^+,
\quad x^{\alpha /2} \Phi _+ ' (x) \sim _{0}  (2-\alpha) \tilde{c} _{\nu _\alpha ,1} ^+  x ^{1-\alpha /2} ,$$
$$ \Phi _- (x) \sim _{0} \tilde{c} _{\nu _\alpha ,0} ^- x ^{1-\alpha},
\quad x^{\alpha /2} \Phi _- ' (x) \sim _{0}  (1-\alpha) \tilde{c} _{\nu _\alpha ,0} ^-  x ^{-\alpha /2} ,$$
hence $x^{\alpha /2} \Phi _+ ' \in L^2 (0,1)$, while $x^{\alpha /2} \Phi _- '
\notin L^2(0,1)$ because $\alpha \geq 1$.
Hence, since $\Phi \in H^1 _\alpha (-1,1)$, one has 
$x^{\alpha /2} \Phi  ' \in L^2 (0,1)$, hence
$x^{\alpha /2} (C_+ \Phi _+ ' + C_- \Phi _- ' )\in L^2 (0,1)$.
This implies that $C_- =0$, hence $\Phi = C_+ \Phi _+$. 

Finally, we look at the boundary condition $\Phi (1)=0$: 
if $C_+ =0$ (which is not forbidden), it is automatically satisfied,
and if $C_+ \neq 0$, then the boundary condition implies that there is some $m \in \Bbb N$, $m\geq 1$ such that
$$ \lambda = \kappa _\alpha ^2  j_{\nu _\alpha, m}^2 ;$$
hence in any case, 
 there is some $C_+$ and some $m\geq 1$ such that \eqref{Phi-struct-d} holds.

In the same way, any $\Phi (x) := C x^{\frac{1-\alpha}{2}} J_{\nu _\alpha} ( j_{\nu _\alpha, m} x^{\kappa_\alpha}) $ is solution of \eqref{*pbm-vp-d}.


\subsubsection{Proof of Lemma \ref{lem-struct-d} when $\nu _\alpha := \frac{\alpha -1}{2-\alpha} \in \mathbb N ^*$}  \hfill

 Let us assume that $\nu _\alpha =n _\alpha \in \mathbb N ^*$. In this case, we have recalled in subsection \ref{*sub-bessel-entier} that a fundamental system of the differential equation \eqref{*eq-bessel-ordre-nu} is given by $J_{ n_\alpha}$ and $Y_{n _\alpha}$. Hence, denoting 
\begin{equation}
\label{*def-phi+-}
\forall x\in (0,1), \quad \Phi _{+,-} (x) := x^{\frac{1-\alpha}{2}} Y_{n_\alpha}  (\frac{2}{2-\alpha} \sqrt{\lambda} x^{\frac{2-\alpha}{2}}),
\end{equation}
 we see that the restriction of $\Phi$ on $(0,1)$ is a linear combination of $\Phi _+$ and $\Phi _{+,-}$. As in the previous case, we are going to study the behavior of
$\Phi _{+,-}$ near $0$. It follows from \eqref{*expand-Yn} that

\begin{multline}
\label{*serie-+-}
\forall x\in (0,1), \quad \Phi _{+,-} (x)
= \frac{2}{\pi} \Phi _{+} (x) \log \left( \frac{1}{2-\alpha} \sqrt{\lambda} x^{\frac{2-\alpha}{2}}\right)  
\\
+  \sum_{m=0}^{n_\alpha-1} \hat a_m x ^{(1-\alpha) +(2-\alpha)m}
+  \sum_{m=0}^{+\infty} \hat b_m x ^{(2-\alpha)m},
\end{multline}
where
$$\hat a_m:= 
- \frac{1}{\pi}\frac{(n_\alpha -m-1)!}{m!}\left(\frac{\sqrt{\lambda}}{2\kappa _\alpha }\right) ^{2m-n_\alpha}$$
and 
$$\hat b_m:=- \frac{1}{\pi} \frac{(-1)^m}{m!(n_\alpha+m)!}\left(\frac{\sqrt{\lambda}}{2\kappa_\alpha}\right) ^{2m+ n_\alpha} \left[ \frac{\Gamma ' (m+1)}{\Gamma (m+1)} + \frac{\Gamma ' (m+n_\alpha+1)}{\Gamma (m+n_\alpha+1)}  \right].$$
We study these three functions that appear in the formula of $\Phi _{+,-}$.
First $$ \Phi _{+,-,1} (x) := \frac{2}{\pi} \Phi _{+} (x) \log \left( \frac{1}{2-\alpha} \sqrt{\lambda} x^{\frac{2-\alpha}{2}}\right)$$
satisfies
$$ \Phi _{+,-,1} (x) \sim _0 \frac{2\kappa_\alpha}{\pi} \tilde{c} _{n _\alpha ,0} ^+ \ln x, \quad x^{\alpha /2} \Phi _{+,-,1}' (x) \sim _0 \frac{2\kappa_\alpha}{\pi} \tilde{c} _{n _\alpha ,0} ^+ x^{-1+\alpha /2},$$
hence $x^{\alpha /2} \Phi _{+,-,1}' \in L^2 (0,1)$ since $\alpha >1$.
Next $$ \Phi _{+,-,2} (x) := \sum_{m=0}^{n_\alpha-1} \hat a_m x ^{(1-\alpha) +(2-\alpha)m} $$
satisfies
$$ \Phi _{+,-,2} (x) \sim _0  \hat{a}_0 x^{1-\alpha}, \quad x^{\alpha /2} \Phi _{+,-,2}' (x) \sim _0 (1-\alpha)  \hat{a}_0 x^{-\alpha /2},$$
hence $x^{\alpha /2} \Phi _{+,-,2}' \notin L^2 (0,1)$, since $\alpha >1$ and $ \hat{a}_0 \neq 0$.
Finally, $$ \Phi _{+,-,3} (x) := \sum_{m=0}^{+\infty} \hat b_m x ^{(2-\alpha)m} $$
satisfies
$$ \Phi _{+,-,3} (x) \sim _0 \hat b _0, \quad x^{\alpha /2} \Phi _{+,-,3}' (x) \sim _0 (2-\alpha) \hat b_1 x^{1-\alpha /2},$$
hence $x^{\alpha /2} \Phi _{+,-,3}' \in L^2 (0,1)$.
Thus $x^{\alpha /2}  \Phi _{+,-} ' \notin L^2(0,1)$, and since
$\Phi = C_+ \Phi _+ + C_{+,-} \Phi _{+,-}$ and $x^{\alpha /2} \Phi ' \in L^2 (0,1)$, then necessarily $C_{+,-}=0$, and 
$\Phi = C_+ \Phi _+$. Then we are in the same position as in the previous case and the conclusion is the same.


\subsubsection{Proof of Lemma \ref{lem-struct-d} when $\nu _\alpha =0$ (hence $\alpha =1$)}  \hfill

In this case, the first sum in the decomposition of $Y_0$ is equal to zero, hence we have $ \Phi _{+,-} = \Phi _{+,-,1} + \Phi _{+,-,3}$. Moreover,
$$ \Phi _{+,-,1} (x) \sim _0 \frac{2\kappa_1}{\pi} \tilde{c} _{0,0} ^+ \ln x, \quad x^{\alpha /2} \Phi _{+,-,1}' (x) \sim _0 \frac{2\kappa_1}{\pi} \tilde{c} _{0,0} ^+ x^{-1/2},$$
hence $x^{\alpha /2} \Phi _{+,-,1}' \notin L^2 (0,1)$. On the contrary $x^{\alpha /2} \Phi _{+,-,3} ' \in L^2(0,1))$,
hence once again $x^{\alpha /2} \Phi _{+,-}' = x^{\alpha /2} (\Phi _{+,-,1}' + \Phi _{+,-,3} ')\notin L^2 (0,1))$, and the conclusion is the same. 

This concludes the proof of Lemma \ref{lem-struct-d}. \qed 


\subsection{The study of \eqref{*pbm-vp-g} when $\alpha \in [1,2)$}  \hfill

We are going to prove the following 

\begin{Lemma}
\label{lem-struct-g}
Assume that $\Phi$ solves \eqref{*pbm-vp}. 
Then there exists some $C_- \in \Bbb R$ (possibly equal to $0$), and some $m'\geq 1$ such that 
\begin{equation}
\label{Phi-struct-g}
\forall x\in (-1,0), \quad \Phi (x) = C_- \vert x \vert ^{\frac{1-\alpha}{2}} J_{\nu _\alpha} ( j_{\nu _\alpha, m'} \vert x \vert ^{\kappa_\alpha}) .
\end{equation}

\end{Lemma}

\noindent {\it Proof of Lemma \ref{lem-struct-g}.} Define
$$ \forall x \in (0,1), \quad \tilde \Phi  (x) := \Phi (-x) .$$
Then $\tilde \Phi$ verifies \eqref{*pbm-vp}, and 
Lemma \ref{lem-struct-g} follows from Lemma \ref{lem-struct-d}. \qed 


\subsection{The study of \eqref{*pbm-vp} when $\alpha \in [1,2)$: proof of Proposition \ref{*prop-vp}}  \hfill

Now assume that $\Phi$ solves \eqref{*pbm-vp}. We derive from Lemmas \ref{lem-struct-d} and \ref{lem-struct-g} that there exist $C_+, C_-, m, m'$ such that
\eqref{Phi-struct-d} and \eqref{Phi-struct-g} hold. 
Then we have
$$\forall x \in (0,1), \quad 
\lambda \Phi = -(x^\alpha \Phi _x)_x = \kappa _\alpha ^2 j_{\nu_\alpha ,m} ^2 \Phi ,
$$
and
$$ \forall x \in (-1,0), \quad 
\lambda \Phi = -(\vert x \vert ^\alpha \Phi _x)_x = \kappa _\alpha ^2 j_{\nu_\alpha ,m'} ^2 \Phi .$$
This is where the discussion according to the values of $C_+$ and $C_-$ arises:


\subsubsection{The possible eigenvalues and eigenfunctions} \hfill
\label{eigen-CN}

\begin{itemize}
\item if $C_-= 0$, then $C_+ \neq 0$ since $\Phi \neq 0$, and therefore $\lambda = \kappa _\alpha ^2 j_{\nu_\alpha ,m} ^2$, and 
$$ \Phi (x) = \begin{cases} 
C_+ x ^{\frac{1-\alpha}{2}} J_{\nu _\alpha} ( j_{\nu _\alpha, m} x ^{\kappa_\alpha}) \quad & \text{ if } x\in (0,1) ,\\
0 \quad & \text{ if } x\in (-1,0)
\end{cases} ;$$
\item if $C_+= 0$, then $C_- \neq 0$ since $\Phi \neq 0$, and therefore $\lambda = \kappa _\alpha ^2 j_{\nu_\alpha ,m'} ^2$, and 
$$ \Phi (x) = \begin{cases} 
0 \quad & \text{ if } x\in (0,1) ,\\
C_- \vert x \vert ^{\frac{1-\alpha}{2}} J_{\nu _\alpha} ( j_{\nu _\alpha, m'} \vert x \vert ^{\kappa_\alpha}) \quad & \text{ if } x\in (-1,0)
\end{cases} ;$$
\item if $C_+ \neq 0 \neq C_-$, then $\kappa _\alpha ^2 j_{\nu_\alpha ,m} ^2
= \kappa _\alpha ^2 j_{\nu_\alpha ,m'} ^2$, hence $m=m'$, and therefore
$\lambda = \kappa _\alpha ^2 j_{\nu_\alpha ,m} ^2$, and 
$$ \Phi (x) = \begin{cases} 
C_+ x ^{\frac{1-\alpha}{2}} J_{\nu _\alpha} ( j_{\nu _\alpha, m} x ^{\kappa_\alpha}) \quad & \text{ if } x\in (0,1) ,\\
C_- \vert x \vert ^{\frac{1-\alpha}{2}} J_{\nu _\alpha} ( j_{\nu _\alpha, m} \vert x \vert ^{\kappa_\alpha}) \quad & \text{ if } x\in (-1,0)
\end{cases} .$$

\end{itemize}
These are necessary conditions, and we have to study if such functions
$\Phi$ are sufficiently smooth, that is if they belong to $H^2 _\alpha (-1,1)$,
of if additionnal conditions on the coefficients appear.


\subsubsection{The eigenvalues and eigenfunctions} \hfill

We first prove the following:

\begin{Lemma}
\label{lem-regPhi}
Given $m\geq 1$, the functions
\begin{equation}
\label{Phi-eigen-d}
\Phi ^{(r)} _{\alpha,m} (x) = \begin{cases} 
 x ^{\frac{1-\alpha}{2}} J_{\nu _\alpha} ( j_{\nu _\alpha, m} x ^{\kappa_\alpha}) \quad & \text{ if } x\in (0,1) ,\\
0 \quad & \text{ if } x\in (-1,0)
\end{cases} 
\end{equation}
and
\begin{equation}
\label{Phi-eigen-g}
\Phi ^{(l)} _{\alpha,m} (x) = \begin{cases} 
0 \quad & \text{ if } x\in (0,1) ,\\
\vert x \vert ^{\frac{1-\alpha}{2}} J_{\nu _\alpha} ( j_{\nu _\alpha, m} \vert x \vert ^{\kappa_\alpha}) \quad & \text{ if } x\in (-1,0)
\end{cases}
\end{equation}
belong to $H^2 _\alpha (-1,1)$.
\end{Lemma}

\noindent {\it Proof of Lemma \ref{lem-regPhi}.} It is clear that $\Phi ^{(r)} _{\alpha,m}$
is locally absolutely continuous on $(0,1]$ and on $[-1,0)$, 
that $\vert x \vert ^{\alpha /2} (\Phi ^{(r)} _{\alpha,m} )' \in L^2 (-1,1)$, that $\Phi ^{(r)} _{\alpha,m} (-1)=0=\Phi ^{(r)} _{\alpha,m} (1)$, and finally, we note that 
$\vert x \vert ^\alpha (\Phi ^{(r)} _{\alpha,m})' \in H^1 (-1,1)$: indeed, 
\begin{multline*}
 \forall x \in (0,1), \quad x^\alpha (\Phi ^{(r)} _{\alpha,m} )' (x) = \sum _{m=0} ^\infty \tilde{c} _{\nu _\alpha ,m} ^+  (2-\alpha) m x ^{(2-\alpha) m-1+\alpha }
\\
= (2-\alpha) \tilde{c} _{\nu _\alpha ,1} ^+  x 
+ 2 (2-\alpha)\tilde{c} _{\nu _\alpha ,2} ^+ x ^{3 -\alpha } + \cdots ,
\end{multline*}
hence $x^\alpha (\Phi ^{(r)} _{\alpha,m} )' \in C^1 ([0,1])$, and since $x^\alpha (\Phi ^{(r)} _{\alpha,m} )' \to 0$ as $x\to 0^+$, the function
$$ \begin{cases} 
x^\alpha (\Phi ^{(r)} _{\alpha,m})' \quad & \text{ if } x\in (0,1) ,\\
0 \quad & \text{ if } x\in (-1,0]
\end{cases}$$
belongs to $H^1 (-1,1)$. We obtain that $\Phi ^{(r)} _{\alpha,m} \in H^2 _\alpha (-1,1)$,
and by symmetry the same property holds for $\Phi ^{(l)} _{\alpha,m}$. \qed

Now we are in position to prove the main part of Proposition \ref{*prop-vp}:
\begin{itemize}
\item we derived from the discussion of subsection \ref{eigen-CN} that the only possible eigenvalues are of the form $\kappa _\alpha ^2 j_{\nu_\alpha ,m} ^2$, and the associated eigenfunctions are of the form $C_+ \Phi ^{(r)} _{\alpha,m}$,
$C_- \Phi ^{(l)} _{\alpha,m}$ or $C_+ \Phi ^{(r)} _{\alpha,m} + C_- \Phi ^{(l)} _{\alpha,m}$;
\item we derive from Lemma \ref{lem-regPhi} that, given $m\geq 1$, and $C_+, C_-$, then the function $\Phi := C_+ \Phi ^{(r)} _{\alpha,m} + C_- \Phi ^{(l)} _{\alpha,m}$ belongs to $ H^2 _\alpha (-1,1)$, and satisfies the eigenvalue problem \eqref{*pbm-vp} with $\lambda = \kappa _\alpha ^2 j_{\nu_\alpha ,m} ^2$.
\end{itemize}
We have thus a necessary and sufficient condition: 
\begin{itemize}
\item the eigenvalues are exactly all the $\kappa _\alpha ^2 j_{\nu_\alpha ,m} ^2$,
\item and the associated eigenspace is exactly the space generated by 
$\Phi ^{(r)} _{\alpha,m}$ and $\Phi ^{(l)} _{\alpha,m}$, which form a basis of the eigenspace.
\end{itemize}
Hence the eigenspace if of dimension 2 and is generated by $\Phi ^{(r)} _{\alpha,m}$ and $\Phi ^{(l)} _{\alpha,m}$.

To conclude the proof of Proposition \ref{*prop-vp}, it remains to study what spaces the differents families generate.


\subsubsection{The spaces generated by the $\Phi ^{(r)} _{\alpha,m}$ and $\Phi ^{(l)} _{\alpha,m}$} \hfill
\label{sec-complet}

We prove the following:

\begin{Lemma}
\label{lem-spaces}
Consider
\begin{equation}
\label{Phi-eigen-d+}
\forall x\in (0,1), \quad \Phi ^{(r,+)} _{\alpha,m} (x) := 
 x ^{\frac{1-\alpha}{2}} J_{\nu _\alpha} ( j_{\nu _\alpha, m} x ^{\kappa_\alpha}) ,
\end{equation}
and
\begin{equation}
\label{Phi-eigen-g-}
\forall x\in (-1,0), \quad \Phi ^{(l,-)} _{\alpha,m} (x) = 
\vert x \vert ^{\frac{1-\alpha}{2}} J_{\nu _\alpha} ( j_{\nu _\alpha, m} \vert x \vert ^{\kappa_\alpha}) .
\end{equation}
Then the following properties hold true:
\begin{itemize}
\item the family $\{ \Phi ^{(r,+)} _{\alpha,m}, m\geq 1 \}$ is orthogonal and complete in $L^2(0,1)$,
\item the family $\{ \Phi ^{(l,-)} _{\alpha,m}, m\geq 1 \}$ is orthogonal and complete in $L^2(-1,0)$,
\item the family $\{ \Phi ^{(r)} _{\alpha, m}, \Phi ^{(l)} _{\alpha, m}, m\geq 1 \}$
is orthogonal and complete in $L^2 (-1,1)$.
\end{itemize}

\end{Lemma}

\noindent {\it Proof of Lemma \ref{lem-spaces}.} The orthogonality comes immediately: indeed, 
\begin{multline*}
 \int _{-1} ^1  \Phi ^{(r)} _{\alpha, m}  \Phi ^{(r)} _{\alpha, n}
= \frac{1}{\lambda _m} \int _{-1} ^1  \lambda _m \Phi ^{(r)} _{\alpha, m}  \Phi ^{(r)} _{\alpha, n}
= \frac{1}{\lambda _m} \int _{-1} ^1 ( -A \Phi ^{(r)} _{\alpha, m} ) \Phi ^{(r)} _{\alpha, n}
\\
= \frac{1}{\lambda _m} \int _{-1} ^1 \Phi ^{(r)} _{\alpha, m}  (-A \Phi ^{(r)} _{\alpha, n})
= \frac{\lambda _n}{\lambda _m} \int _{-1} ^1 \Phi ^{(r)} _{\alpha, m}  \Phi ^{(r)} _{\alpha, n} ,
\end{multline*}
hence $$ \forall m \neq n, \quad \int _{-1} ^1  \Phi ^{(r)} _{\alpha, m}  \Phi ^{(r)} _{\alpha, n} = 0 ,$$
which implies that
$$ \forall m \neq n, \quad \int _{0} ^1  \Phi ^{(r,+)} _{\alpha, m}  \Phi ^{(r,+)} _{\alpha, n} = 0 = \int _{-1} ^0  \Phi ^{(l,-)} _{\alpha, m}  \Phi ^{(l,-)} _{\alpha, n} ;$$
and of course 
$$ \forall m, n \geq 1 \quad \int _{-1} ^1  \Phi ^{(r)} _{\alpha, m}  \Phi ^{(l)} _{\alpha, n} = 0  .$$
This proves all the orthogonality properties stated in Lemma \ref{lem-spaces}.

Now we turn to the properties of completeness. We proceed as in Appendix in \cite{fatiha}.
Consider
$$ X^{(r)} := L^2(0,1) $$
endowed with its natural scalar product
$$ \forall f,g \in X^{(r)}, \quad \langle f,g \rangle = \int _0 ^1 fg ,$$
and next
\begin{multline}
\label{*H^1_a-2-g}
 H^1_{\alpha} (0,1):=\{ u \in L^2 (0,1) \ \mid \  u \text{ locally absolutely continuous in } (0,1]
\\ 
\vert x \vert ^{\alpha/2} u_x \in  L^2 (0,1) \text{ and } x^\alpha u_x (0)=0 =u(1) \} ,
 \qquad \qquad
\end{multline}
endowed with the natural scalar product
$$ \forall f,g \in  H^1_{\alpha} (0,1), \quad (f,g) = 
\int _{0} ^1 \vert x \vert ^\alpha f'(x) g'(x) + f(x) g(x) \, dx .$$
Next, consider 
\begin{equation}
\label{H2alpha-s01}
H^2_{\alpha} (0,1):=   \{ u \in H^1_{\alpha }(0,1)  \ \mid \int _{0} ^1 \vert (x ^\alpha u'(x) ) ' \vert ^2 < \infty  \} ,
\end{equation}
and the operator $A^{(r)}:D(A^{(r)})\subset X^{(r)} \to X^{(r)}$ will be  defined by
\begin{equation}
\label{D(A)-2r}
D(A^{(r)}) :=  H^2_{\alpha} (0,1) \quad \text{ and } \quad
\forall u \in D(A^{(r)}), \quad    A^{(r)} u:= (x ^\alpha  u_x)_x .
\end{equation}
Then the following results hold (as in Proposition \ref{Prop-A-s}):
 $H^1_{\alpha} (0,1)$ is an Hilbert space,
and $A^{(r)}: D(A^{(r)}) \subset X^{(r)} \to X^{(r)}$ is a self-adjoint negative operator with dense domain. It follows from the proof of Lemma \ref{lem-struct-d} that the eigenvalues of $A^{(r)}$ are once again $\lambda _m = \kappa _\alpha ^2 j_{\nu_\alpha,m} ^2$, associated to 
$\Phi ^{(r,+)} _{\alpha,m}$. Then consider $$ T_\alpha: L^2(0,1) \to L^2 (0,1), \quad f \mapsto T_\alpha (f) :=u_f $$
where $u_f \in D(A^{(r)}) $ is the solution of the problem $-A^{(r)} u_f = f $: we claim that $T_\alpha$ is self-adjoint and compact, and then it will be sufficient to apply the diagonalization theorem to conclude that the family $\{ \Phi ^{(r,+)} _{\alpha,m}, m\geq 1 \}$, composed of the eigenfunctions of $T_\alpha$, is complete in $L^2(0,1)$.
This follows from a result of Brezis \cite{Brezis} and some classical steps (see Appendix in \cite{fatiha}):
\begin{itemize}
\item $H^1_{\alpha} (0,1)$ is included in $L^2 (0,1)$: if $u\in H^1_{\alpha} (0,1)$, then
$$ \forall x\in (0,1), \quad u(x) = \int _x ^1 u'(s) \, ds ,$$
hence
\begin{multline}
\label{eq-estim}
\vert u(x) \vert 
\leq \Bigl( \int _x ^1 s^\alpha u'(s) ^2 \, ds \Bigr)^{1/2} 
\Bigl( \int _x ^1 \frac{1}{s^\alpha} \, ds \Bigr)^{1/2}
\\
\leq 
\begin{cases} \Bigl( \int _0 ^1 s^\alpha u'(s) ^2 \, ds \Bigr)^{1/2}  (-\ln x) \quad & \text{ if } \alpha =1 \\
\Bigl( \int _0 ^1 s^\alpha u'(s) ^2 \, ds \Bigr)^{1/2}  (\frac{1}{\alpha -1} \frac{1}{x^{\alpha -1}}) \quad & \text{ if } \alpha \in (1,2) 
\end{cases} ;
\end{multline}
\item \eqref{eq-estim} implies that the injection of $H^1_{\alpha} (0,1)$ into $L^2 (0,1)$ is continuous: indeed,
$$ \int _0 ^1 u(x) ^2 \, dx \leq 
\begin{cases} \Bigl( \int _0 ^1 s^\alpha u'(s) ^2 \, ds \Bigr) \int _0 ^1 (-\ln x) \, dx = \int _0 ^1 s^\alpha u'(s) ^2 \, ds \quad & \text{ if } \alpha = 1 \\
\Bigl( \int _0 ^1 s^\alpha u'(s) ^2 \, ds \Bigr) \int _0 ^1  \frac{1}{\alpha -1} \frac{1}{x^{\alpha -1}} \, dx = \Bigl( \int _0 ^1 s^\alpha u'(s) ^2 \, ds \Bigr) \frac{1}{\alpha -1} \frac{1}{2-\alpha}
 \quad & \text{ if } \alpha \in (1,2) 
\end{cases} ;
$$
\item in the same way, \eqref{eq-estim} also implies that 
\begin{equation}
\label{eq-estim1}
\int _0 ^\delta u(x) ^2 \, dx
+ \int _{1-\delta} ^1 u(x) ^2 \, dx < \omega (\delta) \int _0 ^1 s^\alpha u'(s) ^2 \, ds  
\end{equation}
with 
\begin{equation}
\label{eq-estim1bis}
\omega (\delta ) \to 0 \quad \text{ as } \delta \to 0 ;
\end{equation}
\item and finally, consider $[a,b] \subset (0,1)$, and $\vert h \vert < \min (a, 1-b)$; then 
\begin{multline*}
 \int _a ^b \vert u(x+h) - u(x) \vert ^2 \, dx 
=  \int _a ^b \vert \int _x ^{x+h} u'(s) \, ds \vert ^2 \, dx  
\\
\leq \int _a ^b \Bigl( \int _x ^{x+h} s^\alpha  u'(s) ^2 \, ds \Bigr) \vert \int _x ^{x+h} \frac{1}{s^\alpha} \, ds \vert  \, dx ,
\end{multline*}
hence if $\alpha =1$ and $h>0$, we have
\begin{multline*}
 \int _a ^b \vert u(x+h) - u(x) \vert ^2 \, dx
\leq  \Bigl( \int _0 ^{1} s^\alpha  u'(s) ^2 \, ds \Bigr) \int _a ^b \ln (x+h) - \ln x \, dx
\\
\leq  \Bigl( \int _0 ^{1} s^\alpha  u'(s) ^2 \, ds \Bigr) \Bigl( \int _b ^{b+h} \ln x \, dx - \int _a ^{a+h} \ln x \, dx \Bigr),
\end{multline*}
and clearly
$$ \tilde \omega (h) := \int _b ^{b+h} \ln x \, dx - \int _a ^{a+h} \ln x \, dx \to 0 \quad \text{ as } h \to 0 ^+ ;$$
we can do the same computations when $h<0$ and when $\alpha \in (1,2)$; hence if $\alpha \in [1,2)$ and $\vert h \vert < \min (a, 1-b)$, there exists
an explicit $\tilde \omega _\alpha (h) >0$ (independent of $u$) such that
\begin{equation}
\label{eq-estim2}
\int _a ^b \vert u(x+h) - u(x) \vert ^2 \, dx
\leq \tilde \omega (h) \Bigl( \int _0 ^{1} s^\alpha  u'(s) ^2 \, ds \Bigr) 
\end{equation}
and
\begin{equation}
\label{eq-estim2bis}
\tilde \omega (h) \to  0 \quad \text{ as } \quad h \to 0 ;
\end{equation}
\item then we are in position to apply Corollary IV.26 of Brezis \cite{Brezis}:

\begin{Theorem} (\cite{Brezis})
\label{thm-Brezis}
Consider $(a,b)$ an open interval of $\Bbb R$, and $\mathcal F$ a bounded set of $L^2(a,b)$.
Assume that
\begin{itemize}
\item for all $\varepsilon >0$, for all $[a',b'] \subset I$, there exists $0 < \delta < \min (a'-a, b-b')$ such that
$$ 
\begin{cases} 
\int _{a'} ^{b'} \vert u(x+h) - u(x) \vert ^2 \, dx < \varepsilon 
\\
\forall u \in \mathcal F, \forall \vert h \vert < \delta 
\end{cases} ,
$$
\item and for all $\varepsilon >0$, there exists $0<a''<b''<b$ such that
$$
\begin{cases} 
\int _{a} ^{a''} \vert u(x) \vert ^2 \, dx + \int _{b''} ^{b} \vert u(x) \vert ^2 \, dx < \varepsilon 
\\
\forall u \in \mathcal F 
\end{cases} .
$$
\end{itemize}
Then $\mathcal F$ is relatively compact in $L^2(a,b)$.

\end{Theorem}
Using \eqref{eq-estim1}-\eqref{eq-estim1bis} and \eqref{eq-estim2}-\eqref{eq-estim2bis}, we are in position to apply Theorem \ref{thm-Brezis}, and we obtain that the injection of $H^1 _\alpha (0,1)$ into $L^2(0,1)$ is compact.
\item This implies that the operator $T_\alpha: L^2(0,1) \to L^2(0,1)$ is compact.

\end{itemize}

Hence 
\begin{itemize}
\item the family $\{ \Phi ^{(r,+)} _{\alpha,m}, m\geq 1 \}$ is orthogonal and complete in $L^2(0,1)$;
\item in the same way, the family $\{ \Phi ^{(l,-)} _{\alpha,m}, m\geq 1 \}$ is orthogonal and complete in $L^2(-1,0)$;
\item hence the family $\{ \Phi ^{(r)} _{\alpha,m}, \Phi ^{(l)} _{\alpha,m}, m\geq 1 \}$ is orthogonal and complete in $L^2(-1,1)$.
\end{itemize}
This concludes the proof of Lemma \ref{lem-spaces}. \qed 

Finally, we can norm these orthogonal families: using \cite{Lebedev} (p. 129), we have

$$
\int _0 ^1 x^{1-\alpha} J_{\nu _\alpha} ^2 ( j_{\nu _\alpha, n} x^{\kappa _\alpha }) \, dx
= \frac{1}{\kappa _\alpha} \int _ 0 ^1 z  J_{\nu _\alpha} ^2 ( j_{\nu _\alpha, n} z) \, dz
= \frac{[J' _{\nu _\alpha} (j_{\nu _\alpha,n} ) ]^2 }{2\kappa _\alpha} ,
$$
hence the family given by \eqref{*fp-d} and \eqref{*fp-g} forms an orthonormal basis of $L^2(-1,1)$. The proof of Proposition \ref{*prop-vp} is complete. \qed


\section{Null controllability when $\alpha \in [1,2)$: proof of Proposition \ref{prop-NCstr}}
\label{sec-NCstrong}

\subsection{Negative results} \hfill
\label{ss-neg}

In a classical way, null controllability of \eqref{eq-control} holds if and only if 
observability holds for the adjoint problem: there exists $C>0$ such that
any solution $w$ of the ajoint problem
$$ \begin{cases}
w_t + (\vert x \vert ^\alpha w_x)_x = 0 , \quad x\in (-1,1) \\
w(-1,t)=0=w(1,t) 
\end{cases}$$
has to satisfy
\begin{equation}
\label{eq-obs-s}
\int _{-1} ^1 w(x,0) ^2 \, dx \leq C \int _0 ^T \int _a ^b w(x,t) ^2 \, dx \, dt .
\end{equation}
Assume that $0 < a < b < 1$. Then the observability inequality \eqref{eq-obs-s} cannot be satisfied. Indeed, choose $n\geq 1$, and consider $e^{\lambda _n t} \tilde \Phi ^{(l)} _{\alpha ,n}$: it is a solution of the adjoint problem, but since it is supported on the part $(-1,0)$, it will never satisfy the observability inequality \eqref{eq-obs-s}, hence null controllability does not hold.

More precisely: assume that the initial condition $u_0$ is partly supported in $(-1,0)$. Choose
$n \geq 1$ such that
$$ \int _{-1} ^1 u_0 (x) \Phi ^{(l)} _{\alpha ,n} (x) \, dx 
= \int _{-1} ^1 u_0 (x) \Phi ^{(l,-)} _{\alpha ,n} (x) \, dx \, \neq 0 ;$$
then multiplying \eqref{eq-control} by $\Phi ^{(l)} _{\alpha ,n}$ and integrating by parts, we obtain that
$$ \frac{d}{dt} \langle u , \Phi ^{(l)} _{\alpha ,n}  \rangle 
+ \lambda _n \langle u , \Phi ^{(l)} _{\alpha ,n}  \rangle = 0 ,$$
hence
$$ \langle u (t) , \Phi ^{(l)} _{\alpha ,n}  \rangle
 = e^{-\lambda _n t} \langle u _0 , \Phi ^{(l)} _{\alpha ,n}  \rangle ,$$
 which clearly implies that the control $h$ (supported in $(a,b) \subset (0,1)$) has no influence on the solution $u$, and $u(T)=0$ is not satisfied.


\subsection{Positive result} \hfill

Now we assume that the initial condition is supported in $[0,1)$, hence that
$$ \forall n\geq 1, \quad \langle u _0 , \Phi ^{(l)} _{\alpha ,n}  \rangle = 0 .$$
Then we have
$$\forall n\geq 1, \quad \langle u (t)  , \Phi ^{(l)} _{\alpha ,n}  \rangle = 0 ,$$
hence
$$ u(t) = \sum _{n=1} ^\infty \langle u (t)  , \Phi ^{(r)} _{\alpha ,n}  \rangle \Phi ^{(r)} _{\alpha ,n} .$$
Moreover, multiplying \eqref{eq-control} by $\Phi ^{(r)} _{\alpha ,n}$ and integrating by parts, we obtain that
$$ \frac{d}{dt} \langle u , \Phi ^{(r)} _{\alpha ,n}  \rangle 
+ \lambda _n \langle u , \Phi ^{(r)} _{\alpha ,n}  \rangle = \langle h \chi _{(a,b)}  , \Phi ^{(r)} _{\alpha ,n}  \rangle  .$$
Consider $u^{(r)}$ the solution of
\begin{equation}
\label{eq-contr-r}
\begin{cases}
u^{(r)} _t - (x^\alpha u^{(r)} _x)_x = h(x,t) \chi _{(a,b)} (x), \quad & t >0, \\
(x^\alpha u^{(r)} _x) (0,t)=0=u^{(r)} (1,t) ,\quad & t>0 \\
u^{(r)}(x,0)= u_0(x) , \quad & x \in (0,1)
\end{cases} :
\end{equation}
in the same way, we have
$$ \frac{d}{dt} \langle u^{(r)} , \Phi ^{(r)} _{\alpha ,n}  \rangle 
+ \lambda _n \langle u ^{(r)} , \Phi ^{(r)} _{\alpha ,n}  \rangle = \langle h \chi _{(a,b)}  , \Phi ^{(r)} _{\alpha ,n}  \rangle  ,$$
and then
$$ \forall n\geq 1, \quad \langle u (t)  , \Phi ^{(r)} _{\alpha ,n}  \rangle
= \langle u ^{(r)} , \Phi ^{(r)} _{\alpha ,n}  \rangle ,$$
which implies that
$$ u(t) \chi _{(0,1)}   = u ^{(r)} .$$
Now, applying \cite{sicon2008}, we know that there exists a control $h$ supported in $(a,b)$ such that $u ^{(r)} (T)=0$ (and in \cite{CMV-cost-loc} we constructed such control, using the moment method).
Hence this control drives the solution $u$ of \eqref{eq-control} to $0$ in time $T$. This concludes the proof of Proposition \ref{prop-NCstr}. \qed


\section{The Sturm-Liouville problem in the weakly degenerate case}
\label{sec-SLweak}

The goal of this section is to prove Proposition \ref{*prop-vp-w}:
we study the spectral problem \eqref{*pbm-vp} and the properties of the eigenvalues and eigenfunctions when $\alpha \in [0,1)$.

Once again, one can observe that if $\lambda$ is an eigenvalue, then $\lambda >0$: indeed, multiplying \eqref{*pbm-vp} by $\Phi $ and integrating by parts, then
$$ \lambda \int _{-1} ^1 \Phi ^2 = \int _{-1} ^1 \vert x \vert ^\alpha \Phi _x ^2 ,$$
which implies first $\lambda \geq 0$, and next that $\Phi =0$ if $\lambda =0$.
And, once again, we make the following observation: if $(\lambda,\Phi )$ solves \eqref{*pbm-vp}, \eqref{*pbm-vp-d} and \eqref{*pbm-vp-g} hold true.
In the following, we study \eqref{*pbm-vp-d} and \eqref{*pbm-vp-g} when $\alpha \in [0,1)$, and then we will be able to solve \eqref{*pbm-vp}.


\subsection{The study of \eqref{*pbm-vp-d} when $\alpha \in [0,1)$} \hfill

When $\alpha \in [0,1)$, we note that
$$ \nu _\alpha =\frac{\vert \alpha -1 \vert }{2-\alpha}=  \frac{1-\alpha}{2-\alpha} \in (0,\frac{1}{2}] ,$$ hence
we are in the case where $\nu_\alpha \notin \Bbb N$, hence all the solutions $\Phi$ of \eqref{*pbm-vp-d} are of the form
$$ \Phi = C_+ \Phi _+ + C_- \Phi _- ,$$
where, once again,
\begin{equation}
\label{*base-phi-w}
\forall x \in (0,1), \quad 
\begin{cases}
\Phi _+ (x) := x^{\frac{1-\alpha}{2}} J_{\nu _\alpha} (\frac{2}{2-\alpha} \sqrt{\lambda} x^{\frac{2-\alpha}{2}}) \\
\Phi _- (x) := x^{\frac{1-\alpha}{2}} J_{-\nu _\alpha} (\frac{2}{2-\alpha} \sqrt{\lambda} x^{\frac{2-\alpha}{2}})
\end{cases} ,
\end{equation}
but this time, using the series expansion of $J_{\nu _\alpha}$ and $J_{-\nu _\alpha}$, one obtains
\begin{equation}
\label{*serie-phi-cas2-w}
\forall x \in (0,1), \quad 
\begin{cases} \Phi _+ (x) = \sum _{m=0} ^\infty \tilde{c} _{\nu _\alpha ,m} ^+  x ^{1-\alpha + (2-\alpha) m} \\
\Phi _- (x) = \sum _{m=0} ^\infty \tilde{c} _{\nu _\alpha ,m} ^- x ^{ (2-\alpha) m} 
\end{cases} ,
\end{equation}
where the coefficients $\tilde{c} _{\nu _\alpha ,m} ^+$ and $\tilde{c} _{\nu _\alpha ,m} ^-$ are still defined by \eqref{*coeffs-phi}.
Hence 
$$ \Phi _+ (x) \sim _{0} \tilde{c} _{\nu _\alpha ,0} ^+ x ^{1-\alpha},
\quad x^{\alpha /2} \Phi _+ ' (x) \sim _{0}  (1-\alpha) \tilde{c} _{\nu _\alpha ,0} ^+  x ^{-\alpha /2} ,$$
$$ \Phi _- (x) \sim _{0} \tilde{c} _{\nu _\alpha ,0} ^- ,
\quad x^{\alpha /2} \Phi _- ' (x) \sim _{0}  (2-\alpha) \tilde{c} _{\nu _\alpha ,1} ^-  x ^{1-\alpha /2} ,$$
hence $x^{\alpha /2} \Phi _+ ' \in L^2 (0,1)$, and $x^{\alpha /2} \Phi _- '
\in L^2(0,1)$.

Before using the boundary condition $\Phi (1)=0$, we are going to study 
\eqref{*pbm-vp-g}.


\subsection{The study of \eqref{*pbm-vp-g} when $\alpha \in [0,1)$} \hfill

Take $\Phi$ solution of \eqref{*pbm-vp-d}. Then
$$ \tilde \Phi (x) := \Phi (-x)$$
satisfies \eqref{*pbm-vp-d}, hence there exists $\tilde C_+$, $\tilde C_-$
such that
$$ \tilde\Phi = \tilde C_+ \Phi _+ + \tilde C_- \Phi _- ,$$
hence
$$ \forall x\in (-1,0), \quad 
\Phi (x) = \tilde C_+ \Phi _+ (-x) + \tilde C_- \Phi _-(-x) .$$

Hence, to sum up, if $\Phi$ solves \eqref{*pbm-vp}, then there exist
$C_+, C_-, \tilde C_+, \tilde C_-$ such that
$$\begin{cases}
\Phi (x) = C_+ \Phi _+ (x) + C_- \Phi _- (x) , \quad & \forall x\in (0,1) , \\
\Phi (x) = \tilde C_+ \Phi _+ (-x) + \tilde C_- \Phi _-(-x), \quad & \forall x\in (-1,0) 
\end{cases} .$$


\subsection{The informations given by the space $H^2 _\alpha (-1,1)$ } \hfill

First, as a first consequence of the definition of $H^2 _\alpha (-1,1)$,
$\Phi$ has to be continuous, hence
$$ \lim _{x\to 0^-} \Phi (x) = \lim _{x\to 0^+} \Phi (x) .$$
This implies that 
$$ \tilde C_- \tilde{c} _{\nu _\alpha ,0} ^- = C_- \tilde{c} _{\nu _\alpha ,0} ^- ,$$
hence 
$$  \tilde C_- = C_- .$$

In the same spirit, it follows from the definition of $H^2 _\alpha (-1,1)$
that $x\mapsto \vert x \vert ^\alpha \Phi ' (x)$ has a limit in $0$, hence
$$ \lim _{x\to 0^-} \vert x \vert ^\alpha \Phi ' (x)  = \lim _{x\to 0^+} x^\alpha \Phi ' (x) .$$
This implies that
$$ - (1-\alpha) \tilde C_+ \tilde{c} _{\nu _\alpha ,0} ^+ = (1-\alpha) C_+ \tilde{c} _{\nu _\alpha ,0} ^+  ,$$
hence
$$ \tilde C_+ = - C_+ .$$

Hence, to sum up, if $\Phi$ solves \eqref{*pbm-vp}, then there exist
$C_+, C_-$ such that
$$\begin{cases}
\Phi (x) = C_+ \Phi _+ (x) + C_- \Phi _-(x) , \quad & \forall x\in (0,1) , \\
\Phi (x) = - C_+ \Phi _+ (-x) +  C_- \Phi _-(-x), \quad & \forall x\in (-1,0) 
\end{cases} .$$
Therefore it is natural to introduce
\begin{equation}
\label{def-Phi-imp}
\Phi ^{(e)} :=
\begin{cases}
\Phi _- (x) \quad & \text{ if } x\in (0,1) \\
\Phi _- (-x) \quad & \text{ if } x\in (-1,0)
\end{cases} ,
\end{equation}
and

\begin{equation}
\label{def-Phi-p}
\Phi  ^{(o)} :=
\begin{cases}
\Phi _+ (x) \quad & \text{ if } x\in (0,1) \\
-\Phi _+ (-x) \quad & \text{ if } x\in (-1,0)
\end{cases} ,
\end{equation}
in such a way that $\Phi  ^{(e)}$ is even, 
$\Phi  ^{(o)}$ is odd, and if $\Phi$ solves \eqref{*pbm-vp}, then there exist
$C_+, C_-$ such that
\begin{equation}
\label{w-str}
\forall x\in (0,1), \quad \Phi (x) = C_+ \Phi ^{(o)} (x) + C_- \Phi ^{(e)} (x) .
\end{equation}


\subsection{The informations given by the boundary conditions} \hfill

We have now to use the informations on the boundary conditions:
$\Phi (-1)=0=\Phi (1)$.
Then, \eqref{w-str} implies that
$$ \begin{cases}
C_+ \Phi  ^{(o)} (1) + C_- \Phi ^{(e)} (1) =0, \\
- C_+ \Phi ^{(o)} (1) + C_- \Phi ^{(e)} (1) =0 
\end{cases} ,$$
hence
$$ C_- \Phi  ^{(e)} (1) = 0, \quad \text{ and } \quad C_+ \Phi ^{(o)} (1) = 0 .$$
Now
\begin{itemize}
\item if $C_- = 0$, then $C_+ \neq 0$, and $\Phi ^{(o)} (1) = 0$,
hence there exists $n\geq 1$ such that 
$$ \lambda = \kappa _\alpha ^2 j_{\nu_\alpha, n} ^2 ,$$
and $\Phi = C_+ \Phi  ^{(o)}$; 
\item in the same way, if $C_+ = 0$, then $C_- \neq 0$, and $\Phi ^{(e)} (1) = 0$,
hence there exists $n\geq 1$ such that 
$$ \lambda = \kappa _\alpha ^2 j_{-\nu_\alpha, n} ^2 ,$$
where now $j_{-\nu_\alpha, n}$ denotes the $n^{\text{th}}$ zero of $J_{-\nu_\alpha}$, and $\Phi = C_- \Phi  ^{(o)}$;
\item if $C_- \neq 0$, then $\Phi  ^{(e)} (1) = 0$, hence there exists $n\geq 1$ such that 
$ \lambda = \kappa _\alpha ^2 j_{-\nu_\alpha, n} ^2$; but then 
$$ \Phi  ^{(o)} (1) = \Phi  _+ (1) = J_{\nu_\alpha} (j_{-\nu_\alpha, n}) \neq 0 $$
since $J_{\nu_\alpha}$ and $J_{-\nu_\alpha}$ have different zeros; hence $C_+=0$;
\item in the same way, if $C_+ \neq 0$, then $\Phi  ^{(o)} (1) = 0$, hence there exists $n\geq 1$ such that 
$ \lambda = \kappa _\alpha ^2 j_{\nu_\alpha, n} ^2$; but then 
$$ \Phi  ^{(e)} (1) = \Phi  _- (1) = J_{-\nu_\alpha} (j_{\nu_\alpha, n}) \neq 0 $$
since $J_{\nu_\alpha}$ and $J_{-\nu_\alpha}$ have different zeros; hence $C_-=0$.
\end{itemize}

To sum, we obtain the following necessary conditions:
if $\Phi$ solves \eqref{*pbm-vp}, then 
\begin{itemize}
\item either there exists $n\geq 1$ such that
$$ \lambda = \kappa _\alpha ^2 j_{\nu_\alpha, n} ^2
\quad \text{ and } \quad
\Phi (x) = C_+ x^{\frac{1-\alpha}{2}} J_{\nu _\alpha} (j_{\nu_\alpha, n} x^{\kappa _\alpha}) \text{ on } (0,1) $$
and is odd,
\item or there exists $n\geq 1$ such that
$$ \lambda = \kappa _\alpha ^2 j_{-\nu_\alpha, n} ^2
\quad \text{ and }  \quad 
\Phi (x) = C_- x^{\frac{1-\alpha}{2}} J_{-\nu _\alpha} ( j_{-\nu_\alpha, n} x^{\kappa _\alpha}) \text{ on } (0,1)$$
and is even.
\end{itemize}


\subsection{The eigenvalues and eigenfunctions: proof of Proposition \ref{*prop-vp-w}} \hfill

Now we consider the odd function $\Phi _{\alpha,n} ^{(o)}$ defined in \eqref{Phi-alpha-n-w-i}
and the even function $\Phi _{\alpha,n} ^{(e)}$ defined in \eqref{Phi-alpha-n-w-p}.
From the previous necessary conditions, it is sufficient to prove that $\Phi _{\alpha,n} ^{(o)}, \Phi _{\alpha,n} ^{(e)} \in H^2 _\alpha (-1,1)$. And this follows immediately from the series given in \eqref{*serie-phi-cas2-w}. Indeed, $\Phi _+$ is continuous on $[0,1]$, equal to $0$ at $0$ and $1$, $\Phi ' \in L^1 (0,1)$ hence $\Phi _{\alpha,n} ^{(o)}$ is absolutely continuous on $[0,1]$ and on $[-1,1]$. And the other integrability conditions are also satisfied. The same properties hold for $\Phi _{\alpha,n} ^{(e)}$. \qed


\subsection{The eigenspaces generated by the eigenfunctions} \hfill

Let us denote 
\begin{itemize}
\item $ L^{2,odd} (-1,1)$ the subspace of $L^2(-1,1)$ composed by 
odd functions, and  
\item $ L^{2,even} (-1,1)$ the subspace of $L^2(-1,1)$ composed by 
even functions.
\end{itemize}

We claim the following

\begin{Lemma}
\label{lem-decomp}
When $\alpha \in [0,1)$, we have the following orthogonal decomposition:
$$ L^2 (-1,1) = L^{2,odd} (-1,1) \oplus L^{2,even} (-1,1) .$$
Moreover, 
\begin{itemize}
\item the family $\{ \Phi _{\alpha,n} ^{(o)} , n\geq 1 \}$ forms an orthogonal basis of $L^{2,odd} (-1,1)$, 
\item the family $\{ \Phi _{\alpha,n} ^{(e)} , n\geq 1 \}$ forms an orthogonal basis of $L^{2,even} (-1,1)$.
\end{itemize}
\end{Lemma}

\noindent {\it Proof of Lemma \ref{lem-decomp}.}
It will be an immediate consequence of the following

\begin{Lemma}
\label{lem-decomp-w}
We have the following properties:

\begin{itemize}
\item The restrictions on $(0,1)$ of $\Phi _{\alpha,n} ^{(o)}$ are the eigenfunctions of the eigenvalue problem
\begin{equation}
\label{eq-vp01-DD}
 \begin{cases}
- (x^\alpha \Phi ')'  = \lambda \Phi , \quad x\in (0,1), \\
\Phi (0)=0= \Phi(1) 
\end{cases} ,
\end{equation}
and they form an orthogonal basis of $L^2 (0,1)$.

\item The restrictions on $(0,1)$ of $\Phi _{\alpha,n} ^{(e)}$ are the eigenfunctions of the eigenvalue problem
\begin{equation}
\label{eq-vp01-ND}
\begin{cases}
- (x^\alpha \Phi')' = \lambda \Phi , \quad x\in (0,1), \\
x^\alpha \Phi' (0)=0=\Phi(1) 
\end{cases} ,
\end{equation}
and they form an orthogonal basis of $L^2 (0,1)$.
\end{itemize}
\end{Lemma}

\noindent {\it Proof of Lemma \ref{lem-decomp-w}.}
First we consider the eigenvalue problem \eqref{eq-vp01-DD}.
Then once again, all the solutions of the second order differential equation are of the form
$$ \Phi = C_+ \Phi _+ + C_- \Phi _- ,$$
where $\Phi _+$ and $\Phi _-$ are given in \eqref{*base-phi-w}.
The developments given in \eqref{*serie-phi-cas2-w} and the boundary condition $\Phi(0)=0$
gives $C_-=0$, and the boundary condition $\Phi(1)=0$ gives $\lambda = \kappa _\alpha ^2 j_{\nu_\alpha, n}^2$. Hence the restrictions on $(0,1)$ of $\Phi _{\alpha,n} ^{(o)}$ are the eigenfunctions of the eigenvalue problem \eqref{eq-vp01-DD}.

Now we turn to the properties of completeness. We proceed as in Appendix in \cite{fatiha}
and as in subsection \ref{sec-complet}:
consider
\begin{multline}
\label{*H^1_a-2-g-w}
 H^1_{\alpha} (0,1):=\{ u \in L^2 (0,1) \ \mid \  u \text{ absolutely continuous in } [0,1]
\\ 
\vert x \vert ^{\alpha/2} u_x \in  L^2 (0,1) \text{ and } u (0)=0 =u(1) \} ,
 \qquad \qquad
\end{multline}
endowed with the natural scalar product
$$ \forall f,g \in  H^1_{\alpha} (0,1), \quad (f,g) = 
\int _{0} ^1 \vert x \vert ^\alpha f'(x) g'(x) + f(x) g(x) \, dx .$$
Next, consider 
\begin{equation}
\label{H2alpha-w01}
H^2_{\alpha} (0,1):=   \{ u \in H^1_{\alpha }(0,1)  \ \mid \int _{0} ^1 \vert (x ^\alpha u'(x) ) ' \vert ^2 < \infty  \} ,
\end{equation}
and the operator $A^{(r)}:D(A^{(r)})\subset X^{(r)} \to X^{(r)}$ will be  defined by
\begin{equation}
\label{D(A)-2r'}
D(A^{(r)}) :=  H^2_{\alpha} (0,1) \quad \text{ and } \quad
\forall u \in D(A^{(r)}), \quad    A^{(r)} u:= (x ^\alpha  u_x)_x .
\end{equation}
Then the following results hold (as in Proposition \ref{Prop-A-w}):
 $H^1_{\alpha} (0,1)$ is an Hilbert space,
and $A^{(r)}: D(A^{(r)}) \subset X^{(r)} \to X^{(r)}$ is a self-adjoint negative operator with dense domain. The eigenvalues of $A^{(r)}$ are once again $\lambda _m = \kappa _\alpha ^2 j_{\nu_\alpha,m} ^2$, associated to 
the restrictions on $(0,1)$ of $\Phi ^{(o)} _{\alpha,n}$. Then consider $$ T_\alpha: L^2(0,1) \to L^2 (0,1), \quad f \mapsto T_\alpha (f) :=u_f $$
where $u_f \in D(A^{(r)}) $ is the solution of the problem $-A^{(r)} u_f = f $. Proceeding as in subsection \ref{sec-complet}, it can be seen that $T_\alpha$ is self-adjoint and compact, and then the diagonalization theorem allows us to conclude that the restrictions on $(0,1)$ of $\Phi ^{(o)} _{\alpha,n}$, i.e. the eigenfunctions of $T_\alpha$, form a complete family in $L^2(0,1)$.
This proves the first statement of Lemma \ref{lem-decomp-w}.

Next we consider the eigenvalue problem \eqref{eq-vp01-ND}.
Once again, all the solutions of the second order differential equation are of the form
$$ \Phi = C_+ \Phi _+ + C_- \Phi _- ,$$
where $\Phi _+$ and $\Phi _-$ are given in \eqref{*base-phi-w}.
The developments given in \eqref{*serie-phi-cas2-w} and the boundary condition $x^\alpha \Phi '(0)=0$
gives $C_+=0$, and the boundary condition $\Phi(1)=0$ gives $\lambda = \kappa _\alpha ^2 j_{-\nu_\alpha, n}^2$. Hence the restrictions on $(0,1)$ of $\Phi _{\alpha,n} ^{(e)}$ are the eigenfunctions of the eigenvalue problem \eqref{eq-vp01-ND}.

Concerning the property of completeness, consider now 
\begin{multline}
\label{*H^1_a-2-g-w'}
 \tilde H^1_{\alpha} (0,1):=\{ u \in L^2 (0,1) \ \mid \  u \text{ absolutely continuous in } [0,1]
\\ 
\vert x \vert ^{\alpha/2} u_x \in  L^2 (0,1) \text{ and } x^\alpha u_x (0)=0 =u(1) \} ,
 \qquad \qquad
\end{multline}
endowed with the natural scalar product
$$ \forall f,g \in  H^1_{\alpha} (0,1), \quad (f,g) = 
\int _{0} ^1 \vert x \vert ^\alpha f'(x) g'(x) + f(x) g(x) \, dx .$$
Next, consider 
\begin{equation}
\label{H2alpha-wND}
\tilde H^2_{\alpha} (0,1):=   \{ u \in \tilde H^1_{\alpha }(0,1)  \ \mid \int _{0} ^1 \vert (x ^\alpha u'(x) ) ' \vert ^2 < \infty  \} ,
\end{equation}
and the operator $A^{(r)}:D(A^{(r)})\subset X^{(r)} \to X^{(r)}$ will be  defined by
\begin{equation}
\label{D(A)-2r01}
D(A^{(r)}) :=  \tilde H^2_{\alpha} (0,1) \quad \text{ and } \quad
\forall u \in D(A^{(r)}), \quad    A^{(r)} u:= (x ^\alpha  u_x)_x .
\end{equation}
Then the following results hold (as in Proposition \ref{Prop-A-w}):
 $\tilde H^1_{\alpha} (0,1)$ is an Hilbert space,
and $A^{(r)}: D(A^{(r)}) \subset X^{(r)} \to X^{(r)}$ is a self-adjoint negative operator with dense domain. The eigenvalues of $A^{(r)}$ are once again $\lambda _m = \kappa _\alpha ^2 j_{-\nu_\alpha,m} ^2$, associated to 
the restrictions on $(0,1)$ of $\Phi ^{(e)} _{\alpha,n}$. Then consider $$ T_\alpha: L^2(0,1) \to L^2 (0,1), \quad f \mapsto T_\alpha (f) :=u_f $$
where $u_f \in D(A^{(r)}) $ is the solution of the problem $-A^{(r)} u_f = f $. Proceeding as in subsection \ref{sec-complet}, it can be seen that $T_\alpha$ is self-adjoint and compact, and then the diagonalization theorem allows us to conclude that the restrictions on $(0,1)$ of $\Phi ^{(o)} _{\alpha,n}$, i.e. the eigenfunctions of $T_\alpha$, form a complete family in $L^2(0,1)$.
This proves the second statement of Lemma \ref{lem-decomp-w}. \qed

Lemma \ref{lem-decomp-w} immediately implies Lemma \ref{lem-decomp}. \qed


\section{Null controllability when $\alpha \in [0,1)$: proof of Proposition \ref{prop-NCstr-w}}
\label{sec-NCweak}

We apply the moment method:

\subsection{The eigenvalues and eigenfunctions} \hfill

Given an initial condition $u_0$, we can decompose it on the orthogonal basis of eigenfunctions.
The first thing is to order them: since $J_{\nu_\alpha} (0)=0$ and the zeros of $J_{\nu_\alpha}$ and $J_{-\nu_\alpha}$ are interlaced (because of Sturm's theorems), we have
$$ 0< j_{-\nu_\alpha, 1} < j_{\nu_\alpha, 1} < j_{-\nu_\alpha, 2} < j_{\nu_\alpha, 2} < \cdots ,$$
hence it is natural to denote
$$ \forall n\geq 1, \quad  \lambda _{\alpha, 2n-1} := \kappa _\alpha ^2 j_{-\nu_\alpha, n} ^2, \quad \lambda _{\alpha, 2n} := \kappa _\alpha ^2 j_{\nu_\alpha, n} ^2 ,$$
hence in such a way that
$$ 0 < \lambda _{\alpha ,1} < \lambda _{\alpha ,2} < \lambda _{\alpha ,3 } < \lambda _{\alpha ,4} < \cdots ,$$
and the associated normalized eigenfunctions
\begin{equation}
\label{1310-8}
\forall n\geq 1, \quad  \tilde \Phi  _{\alpha, 2n-1} := \frac{\sqrt{\kappa _\alpha}}{\vert J' _{-\nu _\alpha} (j_{-\nu _\alpha,n} ) \vert } \Phi ^{(e)} _{\alpha,n} , \quad
\tilde \Phi  _{\alpha, 2n} := \frac{\sqrt{\kappa _\alpha}}{\vert J' _{\nu _\alpha} (j_{\nu _\alpha,n} ) \vert } \Phi ^{(o)} _{\alpha,n} 
\end{equation}
form an orthonormal basis of $L^2(-1,1)$.


\subsection{The moment problem satisfied by a control $h \in L^2((a,b)\times (0,T))$} \hfill

First we expand the initial condition $u_0 \in L^2(-1,1)$: 
there exists $(\mu_{\alpha ,n} ^0)_{n\geq 1} \in \ell ^2(\mathbb N^\star)$
such that
$$ u_0 (x) =\sum _{n\geq 1} \mu_{\alpha ,n} ^0 \tilde\Phi_{\alpha ,n} (x), 
\qquad x \in (-1,1 ).$$
Next we expand the solution $u$ of \eqref{eq-control}: 
$$u(x,t)= \sum_{n\geq 1} \beta_{\alpha ,n} (t) \tilde\Phi_{\alpha ,n} (x), \qquad x\in (-1,1), \ t \geq 0
\quad \text{ with } \sum_{n\geq 1} \beta_{\alpha ,n} (t)^2 <+\infty.$$
Multiplying \eqref{eq-control} by $w_{\alpha ,n} (x,t):= \tilde\Phi  _{\alpha ,n} (x) e^{\lambda_{\alpha ,n} (t-T)}$, which is solution of the adjoint problem, one gets:
\begin{multline*}
\int _0 ^T \int _{-1} ^1 h(x,t) \chi_{[a,b]}(x) w_{\alpha ,n} (x,t) \, dx\, dt
= \int _0 ^T \int _{-1} ^1  w_{\alpha ,n} (x,t) (u_t - (\vert x \vert ^\alpha u_x)_x)
\\
= \int _{-1} ^1  [w_{\alpha ,n} u ]_0 ^T - \int _0 ^T \int _0 ^\ell (w_{\alpha ,n})_t u
- \int _0 ^T [  w_{\alpha ,n} (\vert x \vert ^\alpha u_x)]_{-1} ^1 + \int _0 ^T \int _{-1} ^1 (w_{\alpha ,n})_x \vert x \vert ^\alpha u_x
\\
= \int _{-1} ^1  \tilde\Phi_{\alpha ,n} u(T) - e^{-\lambda_{\alpha ,n} T}\int _{-1} ^1  \Phi_{\alpha ,n} u_0
- \int _0 ^T [  w_{\alpha ,n} (\vert x \vert ^\alpha u_x)]_{-1} ^1  
\\
+ \int _0 ^T [ \vert x \vert ^\alpha (w_{\alpha ,n})_x u ]_{-1} ^1 - \int _0 ^T \int _{-1} ^1 \Bigl( (w_{\alpha ,n})_t + ( \vert x \vert ^\alpha (w_{\alpha ,n})_x )_x \Bigr) u
\\
=\int _{-1} ^1  \tilde\Phi_{\alpha ,n} u(T) - e^{-\lambda_{\alpha ,n} T}\int _{-1} ^1  \Phi_{\alpha ,n} u_0 .
\end{multline*}
Hence, if $h$ drives the solution $u$ to $0$ in time $T$, we obtain the following moment problem:
\begin{equation}
\label{*moment2}
\forall n \geq 1, \qquad 
\int_0^T \int_{-1} ^1 h(x,t) \chi_{[a,b]}(x) \tilde \Phi_{\alpha ,n} (x) e^{\lambda_{\alpha ,n} t} dx  dt
=-\mu_{\alpha ,n} ^0  .
\end{equation}
Consider 
\begin{equation}
\label{1310-1}
h_{\alpha ,n}(t) := \int_{-1} ^1 h(x,t) \chi_{[a,b]}(x) \tilde \Phi_{\alpha ,n} (x) dx,\qquad t>0.
\end{equation}
Then the moment problem can be written in the following way:
\begin{equation}
\label{1310-2}
\forall n \geq 1, \quad  
\int_0^T h_{\alpha ,n}(t) e^{\lambda_{\alpha ,n} t} dt = -\mu_{\alpha ,n} ^0.
\end{equation}


\subsection{A formal solution to the moment problem, using a biorthogonal family} \hfill
\label{*s5}

Assume for a moment that there exists a family $(\sigma_{\alpha, m}^+)_{m\geq 1}$ in $L^2(0,T)$ that satisfies 
\begin{equation}
\label{biortho}
\forall m, n \geq 1, \quad \int _0 ^T \sigma _{\alpha,m} ^+ (t) \, e^{\lambda _{\alpha,n}t} \, dt = \delta _{mn} .
\end{equation}
Then let us define
\begin{equation}
\label{*contr} 
h(x,t):= \sum_{m\geq 1} -\mu_{\alpha ,m} ^0 \sigma_{\alpha ,m}^+ (t) \frac{\tilde \Phi_{\alpha ,m}(x)}{\int_a^b \tilde \Phi_{\alpha ,m}^2 } .
\end{equation}
Let us  prove that, formally, $h$ is solution of the moment problem \eqref{*moment2}:
\begin{multline*}
\int_0^T \int_{-1} ^1 h(x,t) \chi_{[a,b]}(x) \tilde \Phi_{\alpha ,n}(x) e^{\lambda_{\alpha ,n} t} dxdt
\\
=
\int_a^b \int_0^T \left(  \sum_{m\geq 1} -\mu_{\alpha ,m} ^0  \sigma_{\alpha ,m}^+(t) \frac{\tilde \Phi_{\alpha ,m}(x)}{\int_a^b \tilde \Phi_{\alpha ,m}^2 }  \right) \tilde \Phi_{\alpha ,n}(x) e^{\lambda_{\alpha ,n} t } dtdx
\\
=\int_a^b \sum_{m\geq 1}-\mu_{\alpha ,m} ^0  \frac{\tilde \Phi_{\alpha ,m}(x) \tilde \Phi_{\alpha, n}(x) }{\int_a^b \tilde \Phi_{\alpha ,m}^2 }
\left( \int_0^T \sigma_{\alpha ,m}^+(t) e^{\lambda_{\alpha ,n} t } dt \right) dx
\\
= \sum_{m\geq 1} -\mu_{\alpha ,m} ^0  \delta_{mn} \frac{\int_a^b \tilde \Phi_{\alpha ,m}(x) \tilde \Phi_{\alpha ,n}(x)dx  }{\int_a^b \tilde \Phi_{\alpha ,m}^2 }  
= -\mu_{\alpha ,n} ^0 .
\end{multline*}
Hence, formally, $h$ defined by \eqref{*contr}  solves the moment problem. It remains to check that all this makes sense, in particular that $h \in L^2((-1,1)\times (0,T))$. 


 
\subsection{The existence of a biorthogonal family} \hfill
\label{s-biortho}

For the existence of a biorthogonal family satisfying \eqref{biortho}, we use the following

\begin{Theorem} (\cite{cost-weak})
\label{thm-biortho1-gen}
Assume that 
$$ \forall n\geq 0, \quad \lambda_n \geq 0, $$
and that there is some $\gamma _{\text{min}}>0$ such that
\begin{equation}
\label{gap-min}
\forall n \geq 0, \quad \sqrt{\lambda _{n+1}} - \sqrt{\lambda _{n}}  \geq \gamma _{\text{min}} .
\end{equation}
Then there exists a family $(\sigma _{m} ^+)_{m\geq 0}$ which is biorthogonal to the family $(e^{\lambda _{n}t})_{n\geq0}$ in $L^2(0,T)$:
\begin{equation}
\label{*famillebi_qm-gen}
\forall m,n \geq 0, \quad \int _0 ^T \sigma _{m} ^+ (t) \,  e^{\lambda _{n}t} \, dt = \delta _{mn} .
\end{equation}
Moreover, it satisfies: there is some universal constant $C_u$ independent of $T$, $\gamma _{\text{min}}$ and $m$ such that, for all $m\geq 0$, we have
\begin{equation}
\label{*famillebi_qm-norme-gen}
\Vert \sigma _{m} ^+ \Vert _{L^2(0,T)} ^2  
\leq C_u e^{-2\lambda _{m} T} e^{C_u \sqrt{\lambda _m} / \gamma _{min} } B(T,\gamma _{min}),
\end{equation}
with
\begin{equation}
\label{eq(B}
B(T,\gamma _{min}) =
\begin{cases} 
\Bigl( \frac{1}{T} + \frac{1}{T^2 \gamma _{min} ^2} \Bigr) 
\, e^{\frac{C_u}{\gamma _{min} ^2T}}  \quad & \text{ if } T \leq \frac{1}{\gamma _{min} ^2} , \\
C_u \gamma _{min} ^2 \quad & \text{ if } T \geq \frac{1}{\gamma _{min} ^2} .
\end{cases}
\end{equation}
\end{Theorem}

To apply it, we need to prove that there is some $\gamma _{min} >0$ such that
\begin{equation}
\label{eq-gap-cond-alpha}
\forall n \geq 1, \quad \sqrt{\lambda _{\alpha, n+1}} - \sqrt{\lambda _{\alpha, n}}  \geq \gamma _{\text{min}} ,
\end{equation}
and this derives from general results about cylinder functions:
\begin{itemize}
\item first the Mac Mahon formula (Watson \cite{Watson} p. 506) gives
$$ j_{\nu_\alpha ,n} = \pi (n+ \frac{1}{2} \nu _\alpha -\frac{1}{4}) + O(\frac{1}{n});$$
\item next, using the Bessel function of second kind $Y_{\nu_\alpha}$ defined in \eqref{def-Ynu},
we have
$$ J_{-\nu_\alpha} (y) = J_{\nu_\alpha} (y) \cos (\nu_\alpha \pi) 
-  Y_{\nu_\alpha} (y) \sin (\nu_\alpha \pi) 
,$$
\item and then, since $\nu_\alpha \pi \in [0,\pi]$, once again the Mac Mahon formula (Watson \cite{Watson} p. 506) gives that
$$ j_{-\nu_\alpha ,n} = \pi (n+ \frac{1}{2} \nu _\alpha -\frac{1}{4}) - \nu_\alpha \pi + O(\frac{1}{n}) 
= \pi (n - \frac{1}{2} \nu _\alpha -\frac{1}{4}) + O(\frac{1}{n});$$
\item now we are in position to conclude:
\begin{multline*}
 \sqrt{\lambda _{\alpha, 2n+1}} - \sqrt{\lambda _{\alpha, 2n}}
= \kappa _\alpha (j_{-\nu_\alpha ,n+1}-j_{\nu_\alpha ,n}) 
\\
=  \kappa _\alpha \pi \Bigl( (n+1 - \frac{1}{2} \nu _\alpha -\frac{1}{4})
- (n+ \frac{1}{2} \nu _\alpha -\frac{1}{4}) \Bigr) + O(\frac{1}{n}) 
\\
= \kappa _\alpha \pi (1 - \nu _\alpha) + O(\frac{1}{n}) 
\longrightarrow \kappa _\alpha \pi (1 - \nu _\alpha) \text{ as } n \to \infty ,
\end{multline*}
and in the same way
\begin{multline*}
 \sqrt{\lambda _{\alpha, 2n}} - \sqrt{\lambda _{\alpha, 2n-1}}
= \kappa _\alpha (j_{\nu_\alpha ,n}-j_{-\nu_\alpha ,n}) 
\\
=  \kappa _\alpha \pi \Bigl( (n + \frac{1}{2} \nu _\alpha -\frac{1}{4})
- (n- \frac{1}{2} \nu _\alpha -\frac{1}{4}) \Bigr) + O(\frac{1}{n}) 
\\
= \kappa _\alpha \pi \nu _\alpha + O(\frac{1}{n}) 
\longrightarrow \kappa _\alpha \pi  \nu _\alpha \text{ as } n \to \infty ;
\end{multline*}
\item since 
$$ \kappa _\alpha \pi (1 - \nu _\alpha) >0 \quad \text{ and } \quad \kappa _\alpha \pi  \nu _\alpha >0 ,$$
there exists $\gamma _{min}(\alpha) >0$ such that \eqref{eq-gap-cond-alpha} holds.
\end{itemize}
Then Theorem \ref{thm-biortho1-gen} gives the existence of a biorthogonal family $(\sigma _{\alpha ,m} ^+)_{m\geq 1}$ satisfying \eqref{biortho}, and we derive from \eqref{*famillebi_qm-norme-gen} that
\begin{equation}
\label{*famillebi_qm-norme-gen-alpha}
\Vert \sigma _{\alpha,m} ^+ \Vert _{L^2(0,T)} ^2  
\leq C_u e^{-2\lambda _{\alpha,m} T} e^{C_u \sqrt{\lambda _{\alpha,m}} / \gamma _{min} } B(T,\gamma _{min}).
\end{equation}


\subsection{A lower bound for the norm of the eigenfunctions on the control region} \hfill
\label{s-minfp}

The last thing is to obtain a lower bound for $ \int _a ^b \tilde \Phi _{\alpha,m} ^2$. We prove the following:

\begin{Lemma}
\label{lem-minfp}
Given $\alpha \in [0,1)$, there exists $\mu (\alpha,a,b) >0$ such that
\begin{equation}
\label{eq-minfp}
\forall m \geq 1, \quad \int _a ^b \tilde \Phi _{\alpha,m} ^2 \geq \mu (\alpha,a,b) .
\end{equation}
\end{Lemma}

\noindent {\it Proof of Lemma \ref{lem-minfp}.} Of course, since the Bessel functions are nonzero solutions of a second order ordinary differential equation, we have
$$ \forall n \geq 1, \quad \int _a ^b \tilde \Phi _{\alpha,m} ^2 >0 .$$
We are going to study the behavior of this quantity as $m\to \infty$. 

Consider that $m$ is even: $m=2n$. Then,
%
%
using the change of variables $z= j_{\nu_{\alpha,n}} x^{\kappa _\alpha}$, we have
\begin{multline*}
\int _a ^b \tilde \Phi _{\alpha,2n} ^2
= \frac{\kappa _\alpha}{ J' _{\nu_\alpha} (j_{\nu_{\alpha,n}}) ^2} \int _a ^b x^{1-\alpha} 
J_{\nu_\alpha} (j_{\nu_{\alpha,n}} x^{\kappa _\alpha}) ^2 \, dx
\\
= \frac{1}{ J' _{\nu_\alpha} (j_{\nu_{\alpha,n}}) ^2 \, j_{\nu_{\alpha,n}} ^2} \int _{j_{\nu_{\alpha,n}} a^{\kappa _\alpha}} ^{j_{\nu_{\alpha,n}} b^{\kappa _\alpha}} z J_{\nu_\alpha} (z) ^2 \, dz.
\end{multline*}
%
Now we use the classical asymptotic development (Lebedev \cite{Lebedev}, formula (5.11.6) p. 122):
\begin{equation}
\label{Leb-asymp1}
J_\nu (z) = \sqrt{\frac{2}{\pi z}} \Bigl[ \cos(z-\frac{\nu\pi}{2}-\frac{\pi}{4}) (1+ O(\frac{1}{z^2})) + O(\frac{1}{z}) \Bigr] \quad \text{as } z\to \infty ,
\end{equation}
and the fact that
$$ J' _{\nu_\alpha} (j_{\nu_{\alpha,n}}) ^2 = J _{\nu_\alpha +1} (j_{\nu_{\alpha,n}}) ^2.$$
Hence
\begin{equation}
\label{Leb-asymp2}
J_\nu (z) = \sqrt{\frac{2}{\pi z}} \cos(z-\frac{\nu\pi}{2}-\frac{\pi}{4})+ O(\frac{1}{z^{3/2}}) \quad \text{as } z\to \infty ,
\end{equation}
and
\begin{equation}
\label{Leb-asymp3}
J_\nu (z)^2  = \frac{2}{\pi z} \cos^2 (z-\frac{\nu\pi}{2}-\frac{\pi}{4})+ O(\frac{1}{z^{2}}) \quad \text{as } z\to \infty .
\end{equation}
Then
\begin{multline*}
\int _{j_{\nu_{\alpha,n}} a^{\kappa _\alpha}} ^{j_{\nu_{\alpha,n}} b^{\kappa _\alpha}} z J_{\nu_\alpha} (z) ^2 \, dz 
= \int _{j_{\nu_{\alpha,n}} a^{\kappa _\alpha}} ^{j_{\nu_{\alpha,n}} b^{\kappa _\alpha}} 
\frac{2}{\pi } \cos^2 (z-\frac{\nu_\alpha \pi}{2}-\frac{\pi}{4})+ O(\frac{1}{z}) \, dz 
\\
= \frac{1}{\pi } \int _{j_{\nu_{\alpha,n}} a^{\kappa _\alpha}} ^{j_{\nu_{\alpha,n}} b^{\kappa _\alpha}} 
1+ \cos (2z-\nu_\alpha \pi-\frac{\pi}{2})+ O(\frac{1}{z}) \, dz 
\\
= j_{\nu_{\alpha,n}} \frac{b^{\kappa _\alpha}-a^{\kappa _\alpha}}{\pi} + O(1).
\end{multline*}
Hence
$$ \int _a ^b \tilde \Phi _{\alpha,2n} ^2
= \frac{1}{j_{\nu_{\alpha,n}}\, J _{\nu_\alpha +1} (j_{\nu_{\alpha,n}}) ^2} \Bigl( \frac{b^{\kappa _\alpha}-a^{\kappa _\alpha}}{\pi} + O(\frac{1}{j_{\nu_{\alpha,n}}}) \Bigr).$$
And finally, once again with the asymptotic development \eqref{Leb-asymp3}, we have
$$  z J_{\nu +1} (z)^2  = \frac{2}{\pi} \cos^2 (z-\frac{(\nu+1)\pi}{2}-\frac{\pi}{4})+ O(\frac{1}{z})
= \frac{2}{\pi} \sin^2 (z-\frac{\nu\pi}{2}-\frac{\pi}{4})+ O(\frac{1}{z}) ,$$
hence
$$ z J_{\nu } (z)^2 + z J_{\nu +1} (z)^2 = \frac{2}{\pi} + O(\frac{1}{z}) ,$$
which gives that
$$ z J_{\nu } (z)^2 + z J_{\nu +1} (z)^2 \to \frac{2}{\pi} \quad \text{as } z \to +\infty ,$$
hence
$$ j_{\nu_{\alpha,n}}\, J _{\nu_\alpha +1} (j_{\nu_{\alpha,n}}) ^2 \to \frac{2}{\pi} \quad \text{as } n \to \infty.$$
Therefore we obtain
$$ \int _a ^b \tilde \Phi _{\alpha,2n} ^2 \to \frac{b^{\kappa _\alpha}-a^{\kappa _\alpha}}{2} \quad \text{as } n \to \infty.$$

This remains true when $m$ is odd ($m=2n+1$), since the main formulae we used concerning Bessel functions remain valid when the order $\nu \geq \frac{-1}{2}$. Hence, in the same way,
$$ \int _a ^b \tilde \Phi _{\alpha,2n+1} ^2 \to \frac{b^{\kappa _\alpha}-a^{\kappa _\alpha}}{2} \quad \text{as } n \to \infty,$$
and this concludes the proof of Lemma \ref{lem-minfp}. \qed


\subsection{Conclusion: proof of Proposition \ref{prop-NCstr-w}} \hfill

We have all the arguments to conclude: 
\begin{itemize}
\item there exists a biorthogonal family $(\sigma _{\alpha ,m} ^+)_{m\geq 1}$ satisfying the biorthogonal properties \eqref{biortho} and the $L^2$ bounds \eqref{*famillebi_qm-norme-gen-alpha} (see subsection \ref{s-biortho}); 

\item then consider the control $h$ given by \eqref{*contr}: using the $L^2$ bounds
\eqref{*famillebi_qm-norme-gen-alpha} and the bound from below of the eigenfunctions in the control region \eqref{eq-minfp}, we obtain that 
$h\in L^2 ((-1,1)\times (0,T))$, 

\item and using the biorthogonal properties \eqref{biortho}, we obtain that the relations of the moment problem are all satisfied: 
$$
\forall n \geq 1, \qquad 
\int_0^T \int_{-1} ^1 h(x,t) \chi_{[a,b]}(x) \tilde \Phi_{\alpha ,n} (x) e^{\lambda_{\alpha ,n} t} dx  dt
=-\mu_{\alpha ,n} ^0  ;
$$
\item this implies that the solution $u$ of the well-posed problem \eqref{eq-control} satisfies
$$ \forall n\geq 1, \quad \langle u(T), \tilde \Phi _{\alpha,n} \rangle = 0 ,$$
hence $u(T)=0$; hence the control $h$ drives the solution $u$ to $0$ in time $T$. \qed
\end{itemize}


\subsection{The cost of controllability as $\alpha \to 1^-$: Proof of Theorem \ref{prop-cost}} \hfill

The proof of Theorem \ref{prop-cost} follows from the same arguments used in \cite{cost-weak, CMV-cost-loc}. First we analyze the behavior of the eigenvalues when $\alpha \to 1^-$, and we use it to understand the behavior of associated biorthogonal families and of the null controllability cost.

The starting point is the following: since the eigenvalues are given by
$$ \forall n\geq 1, \quad \lambda _{\alpha, 2n-1} = \kappa _\alpha ^2 j_{-\nu _\alpha, n} ^2, \quad \text{ and } \quad \lambda _{\alpha, 2n} = \kappa _\alpha ^2 j_{\nu _\alpha, n} ^2 ,$$
and since $\nu_\alpha = \frac{1-\alpha}{2-\alpha}$, we have
$$ \forall n\geq 1, \quad \lambda _{\alpha, 2n} - \lambda _{\alpha, 2n-1} \to 0 \quad \text{ as } \alpha \to 1^- .$$
This will cause the blow-up of the null controllability cost. The first thing to do is to precise the behavior of the gap between successive eigenvalues.


\subsubsection{The behavior of the eigenvalues when $\alpha \to 1^-$} \hfill

It is classical that the function $\nu \mapsto j_{\nu,n}$ is $C^1$ and increasing, see, e.g., Watson \cite{Watson}, p. 508. We prove the following uniform result:

\begin{Lemma}
\label{lem-j_nu-croissant}
Given $\eta >0$, there exists $0<m<M$ such that
\begin{equation}
\label{1002-1}
\forall \nu \in [-\eta, \eta], \forall n \geq 1, \quad m \leq \frac{d j_{\nu,n} }{d \nu} \leq M .
\end{equation}
\end{Lemma}

\noindent {\it Proof of Lemma \ref{lem-j_nu-croissant}.} 
We get from Watson \cite{Watson}, p. 508 that
\begin{equation}
\label{1003-0}
\frac{d j_{\nu,n} }{d \nu} = 2 j_{\nu,n} \int _0 ^{+\infty} K_0 (2j_{\nu,n} \sinh t) e^{-2\nu t} \, dt ,
\end{equation}
with the following formula (\cite{Watson}, p. 181):
\begin{equation}
\label{1003-0bis}
\forall x >0, \quad K_0 (x) = \int _0 ^{+\infty} e^{-x \cosh v} \, dv .
\end{equation}
Using that $\frac{1}{2} e^v \leq \cosh v \leq e^v$,  we obtain that
$$ \forall x>0, \quad \int _0 ^{+\infty} e^{- e^{v + \ln x} } \, dv
\leq K_0 (x) \leq \int _0 ^{+\infty} e^{- e^{v + \ln \frac{x}{2}} } \, dv ,
$$
hence
\begin{equation}
\label{1003-1}
\forall x>0, \quad F(\ln x ) \leq K_0 (x) \leq F(\ln \frac{x}{2} ), \quad
\text{ where } \quad F(X) := \int _X ^{+\infty} e^{- e^{w} } \, dw  .
\end{equation}
We investigate the behavior of $F(X)$ as $X\to \pm \infty$:
\begin{itemize}
\item when $X \geq 0$: since 
$$ \frac{d}{dw} e^{-e^w} = -e^w \, e^{-e^w} ,$$
an integration by parts gives us that
$$ \int _X ^{+\infty} (1+e^{-w}) \, e^{- e^{w} } \, dw = e^{-X} \, e^{-e^X} ,$$
and taking into account that $1\leq 1+e^{-w} \leq 2$, we obtain that
\begin{equation}
\label{1003-2}
\forall X \geq 0, \quad \frac{1}{2} e^{-X} \, e^{-e^X} \leq F(X) \leq e^{-X} \, e^{-e^X} ;
\end{equation}
\item when $X \leq 0$: decomposing 
$$ F(X) = \int _X ^0 e^{- e^{w} } \, dw + F(0), $$
and noting that
$$ \forall w\leq 0, \quad \frac{1}{e} \leq e^{- e^{w} } \leq 1 ,$$
we obtain that
\begin{equation}
\label{1003-3}
\forall X \leq 0, \quad -\frac{1}{e} X + \frac{1}{2e} \leq F(X) \leq -X + \frac{1}{e} .
\end{equation}
\end{itemize}
Then we derive from \eqref{1003-1}-\eqref{1003-3} the following 
\begin{itemize}
\item bounds from below:
\begin{equation}
\label{1003-4}
\forall x \in (0,1], \quad K_0 (x) \geq -\frac{1}{e} \ln x + \frac{1}{2e},
\quad \text{ and } \quad 
\forall x \geq 1, \quad K_0 (x) \geq \frac{1}{2x} e^{-x} ,
\end{equation}
\item and bounds from above:
\begin{equation}
\label{1003-5}
\forall x \in (0,2], \quad K_0 (x) \leq -\ln \frac{x}{2} + \frac{1}{e},
\quad \text{ and } \quad 
\forall x \geq 2, \quad K_0 (x) \leq \frac{2}{x} e^{-x/2} .
\end{equation}
\end{itemize}
This enables us to estimate $\frac{d j_{\nu,n} }{d \nu}$:
\begin{itemize}
\item bound from below: consider $\tau_{\nu,n}$ such that
$$ 2 j_{\nu, n} \sinh \tau_{\nu,n} = 1 ;$$
then we derive from \eqref{1003-0}, \eqref{1003-0bis} and \eqref{1003-4} that
\begin{multline*}
\frac{d j_{\nu,n} }{d \nu} \geq 2 j_{\nu,n} \int _0 ^{\tau_{\nu,n}} K_0 (2j_{\nu,n} \sinh t) e^{-2\nu t} \, dt 
\\
\geq 2 j_{\nu,n} \int _0 ^{\tau_{\nu,n}}  \Bigl( -\frac{1}{e} \ln (2j_{\nu,n} \sinh t) + \frac{1}{2e} \Bigr) e^{-2\nu t} \, dt 
\geq \frac{2 j_{\nu,n} }{2e} \tau_{\nu,n} \, e^{-2\vert \nu \vert \tau_{\nu,n}} ;
\end{multline*}
since $j_{\nu,n}\tau_{\nu,n} \to \frac{1}{2}$ as $n\to \infty$, we obtain the bound from below of \eqref{1002-1};
\item bound from above: consider $\tau ' _{\nu,n}$ such that
$$ 2 j_{\nu, n} \sinh \tau ' _{\nu,n} = 2 ;$$
then we derive from \eqref{1003-0}, \eqref{1003-0bis} and \eqref{1003-5} that
\begin{multline*}
\frac{d j_{\nu,n} }{d \nu} = 2 j_{\nu,n} \int _0 ^{\tau ' _{\nu,n}} K_0 (2j_{\nu,n} \sinh t) e^{-2\nu t} \, dt + 2 j_{\nu,n} \int _{\tau ' _{\nu,n}} ^{+\infty} K_0 (2j_{\nu,n} \sinh t) e^{-2\nu t} \, dt
\\
\leq 2 j_{\nu,n} \int _0 ^{\tau ' _{\nu,n}} \Bigl( -\ln \frac{2j_{\nu,n} \sinh t}{2} + \frac{1}{e} \Bigr) e^{-2\nu t} \, dt + 2 j_{\nu,n} \int _{\tau ' _{\nu,n}} ^{+\infty} \frac{2}{2j_{\nu,n} \sinh t} e^{-j_{\nu,n} \sinh t} e^{-2\nu t} \, dt 
\\
= I_1+I_2+I_3 ;
\end{multline*}
it remains to look to the three integrals of the righ-hand side; take $\nu \in [-\eta, \eta]$; then
\begin{itemize}
\item first, since $j_{\nu,n}\tau_{\nu,n} \to 1$ as $n\to \infty$, we have that
$$ I_2 = 2 j_{\nu,n} \int _0 ^{\tau ' _{\nu,n}} \frac{1}{e} e^{-2\nu t} \, dt
\leq \frac{2}{e} j_{\nu,n} \tau ' _{\nu,n} e^{2 \vert \nu \vert \tau ' _{\nu,n}} \leq M_2 ;$$
\item next, using the change of variables $y= j_{\nu,n} \sinh t$, we have
\begin{multline*}
I_1= 2 j_{\nu,n} \int _0 ^{\tau ' _{\nu,n}} -\ln \frac{2j_{\nu,n} \sinh t}{2} e^{-2\nu t} \, dt
\leq 2 j_{\nu,n} e^{2 \vert \nu \vert \tau ' _{\nu,n}} \int _0 ^{\tau ' _{\nu,n}} -\ln ( j_{\nu,n} \sinh t ) \, dt
\\
= 2 j_{\nu,n} e^{2 \vert \nu \vert \tau ' _{\nu,n}} \int _0 ^1 -\ln y \, \frac{1}{j_{\nu,n}} \, \frac{1}{\sqrt{1+\frac{y^2}{j_{\nu,n} ^2}}} \, dy
\leq 2  e^{2 \vert \nu \vert \tau ' _{\nu,n}} \int _0 ^1 -\ln y \, dy \leq M_1 ;
\end{multline*}
\item finally, in the same way with $y= j_{\nu,n} \sinh t$, we have
\begin{multline*}
I_3= 2 j_{\nu,n} \int _{\tau ' _{\nu,n}} ^{+\infty} \frac{2}{2j_{\nu,n} \sinh t} e^{-j_{\nu,n} \sinh t} e^{-2\nu t} \, dt
\\
= 2 \int _1 ^{+\infty} \frac{1}{\frac{y}{j_{\nu,n}}} e^{-y} e^{-2\nu \text{ argsinh } \frac{y}{j_{\nu,n}}} \frac{1}{j_{\nu,n}} \frac{1}{\sqrt{1+\frac{y^2}{j_{\nu,n} ^2}}} \, dy
\\
\leq 2 \int _1 ^{+\infty} e^{-y} e^{2 \vert \nu \vert \text{ argsinh } \frac{y}{j_{\nu,n}}}  \, dy \leq M_3 ;
\end{multline*}
\end{itemize}
using these three estimates, we obtain the bound from above of \eqref{1002-1}. 
\end{itemize}
This concludes the proof of Lemma \ref{lem-j_nu-croissant}.\qed

We immediately deduce from Lemma \ref{lem-j_nu-croissant} the following

\begin{Lemma}
\label{lem-appl-j_nu-croissant}
There exists $0<m_* \leq M_*$ and $\alpha _* \in (0,1)$ such that
\begin{equation}
\label{1006-1}
\forall \alpha \in [\alpha _*,1), \forall n \geq 1, 
\quad m_* (1-\alpha) \leq j_{\nu _\alpha,n} - j_{-\nu _\alpha,n} \leq M_* (1-\alpha)  .
\end{equation}
\end{Lemma}
It follows from Lemma \ref{lem-appl-j_nu-croissant} that we are in the following situation: there are two subfamilies of eigenvalues: $(\lambda _{\alpha, 2n-1})_{n\geq 1}$ and $(\lambda _{\alpha, 2n})_{n\geq 1}$, these two subfamilies satisfies the following uniform gap conditions: there exists some $\gamma _{min} >0$ independent of $\alpha \in [0,1)$ such that
$$ \forall n \geq 1, \quad \sqrt{\lambda _{\alpha, 2n+2}} -  \sqrt{\lambda _{\alpha, 2n}} \geq \gamma_{min} ,
\quad
\sqrt{\lambda _{\alpha, 2n+1}} -  \sqrt{\lambda _{\alpha, 2n-1}} \geq \gamma_{min} ,$$
and on the other hand
$$ \forall n\geq 1, \quad \sqrt{\lambda _{\alpha, 2n-1}} \rightarrow \frac{1}{2} j_{0,n} \leftarrow \sqrt{\lambda _{\alpha, 2n}} \quad \text{ as } \alpha \to 1 ^- ,$$
with
$$ \forall n \geq 1, \quad m_* (1-\alpha) \leq \sqrt{\lambda _{\alpha, 2n}} -  \sqrt{\lambda _{\alpha, 2n-1}} \leq 2 M_* (1-\alpha) .$$
In the following we estimate the null controllability cost when $\alpha \to 1^-$.


\subsubsection{Upper bound of the null controllability cost: proof of \eqref{1310-12}} \hfill

In the situation described in the previous subsection concerning the eigenvalues, Theorem \ref{thm-biortho1-gen} can be applied with some $\gamma _{min} (\alpha)$ of the order $1-\alpha$,
which gives an upper bound of the constructed biorthogonal family of the order $e^{1/(1-\alpha)}$. But in fact {\it the proof} of Theorem \ref{thm-biortho1-gen} gives a much better upper bound: here we do not write everything in detail, but we indicate how to adapt the proof of Theorem \ref{thm-biortho1-gen} in \cite{cost-weak} (adapted from \cite{Seid-Avdon}) to our present case:

\begin{itemize}
\item the counting function: consider
$$ N_{\alpha, n} (\rho) := \text{ card } \{k, 0< \vert \lambda _{\alpha ,n} - \lambda _{\alpha ,k} \vert \leq \rho \};$$ 
then in our context, one easily sees that there exists $c_u, C_u >0 $ independent of $n$ and of $\alpha$ such that
\begin{equation}
\label{1008-1}
\begin{cases}
\forall \rho \leq c_u (1-\alpha), \quad N_{\alpha, n} (\rho) =0 , \\
\forall \rho \geq c_u (1-\alpha), \quad N_{\alpha, n} (\rho) \leq 1+ C_u \sqrt{\rho} ;
\end{cases} 
\end{equation}
indeed, there is at most only one eigenvalue close to $\lambda _{\alpha ,n}$; then, if $\rho$ is small enough (of the order $1-\alpha$), $N_{\alpha, n} (\rho) =0$, and if $\rho$ exceeds this threshold value, then $N_{\alpha, n} (\rho)$ takes into account the eigenvalue close to $\lambda _{\alpha ,n}$, and the others, whose number is $O(\sqrt{\rho})$ (see \cite{cost-weak});

\item the growth of the Weierstrass product: we use the counting function $N_{\alpha, m}$ to estimate the growth of 
$$ F_{\alpha, m} (z) := \prod _{k=1, k\neq m} ^\infty \Bigl( 1 - \Bigl( \frac{iz - \lambda _{\alpha,m}}{\lambda _{\alpha,k} - \lambda _{\alpha,m}}\Bigr) ^2 \Bigr) , $$
and combining the proof of \cite{cost-weak} with \eqref{1008-1}, we obtain
\begin{multline*}
\ln \vert F_{\alpha,m} (z - i \lambda _{\alpha ,m} ) \vert
\leq 2 \int _0 ^\infty \frac{N_{\alpha, m} (\rho)}{\rho} \frac{\vert z \vert ^2 }{\vert z \vert ^2 + \rho ^2}  \, d\rho
\\
\leq 2 \int _{c_u (1-\alpha)} ^\infty \frac{1+ C_u \sqrt{\rho}}{\rho} \frac{\vert z \vert ^2 }{\vert z \vert ^2 + \rho ^2}  \, d\rho
= \ln \Bigl( \frac{\vert z \vert ^2}{c_u ^2 (1-\alpha)^2} +1 \Bigr) + O(1) \sqrt{\vert z \vert } ,
\end{multline*}
which gives that 
$$ \vert F_{\alpha,m} (z) \vert 
\leq \Bigl( 2 \frac{\vert z \vert ^2 + \lambda _{\alpha,m} ^2}{c_u ^2 (1-\alpha)^2}  +1 \Bigr) \, e^{O(1) ( \sqrt{\vert z\vert}  + \sqrt{\lambda _{\alpha,m} } )} ,
$$
hence there exists some $C_u '$ independent of $\alpha \in [0,1)$ and of $m\geq 1$ such that
\begin{equation}
\label{1008-2}
\forall m\geq 1, \forall z \in \Bbb C, \quad \vert F_{\alpha,m} (z) \vert 
\leq \frac{C_u ' }{(1-\alpha)^2} \, e^{C_u ' ( \sqrt{\vert z\vert}  + \sqrt{\lambda _{\alpha,m} } )} ;
\end{equation}

\item then (following \cite{cost-weak}), we obtain the existence of a family
$(\sigma _{m} ^+)_{m\geq 1}$ biorthogonal to the family $(e^{\lambda _{\alpha,n}t})_{n\geq 1}$
satisfying
\begin{equation}
\label{1008-3}
\Vert \sigma _{m} ^+ \Vert _{L^2(0,T)} 
\leq \frac{C'' _u }{(1-\alpha)^2} e^{-\lambda _{\alpha , m} T} e^{C'' _u \sqrt{\lambda _{\alpha,m}} } e^{C''_u /T} ,
\end{equation}
where the (possible) blow-up of $\Vert \sigma _{m} ^+ \Vert _{L^2(0,T)} $ with respect to the parameter $\alpha$ is controlled by the term $\frac{1}{(1-\alpha)^2}$;

\item as noted in subsection \ref{s-minfp}, we need to have an estimate from below for the $L^2$ norm of the eigenfunctions in the control region, but here it is sufficient to note from \eqref{*base-phi-w}-\eqref{*serie-phi-cas2-w} that
the variation of the eigenfunctions with respect to $\alpha$ is bounded, and their norm is uniformy bounded from below by a positive constant when $\alpha =1$, hence the same property holds for all $\alpha \in [0,1)$: there exists $m^* (a,b) >0$ such that
\begin{equation}
\label{1008-4}
\forall \alpha \in [0,1), \forall n\geq 1, \quad \int _a ^b \Phi _{\alpha,n} ^2 
\geq m^* ;
\end{equation}

\item finally, this allows us to bound from above the null controllability cost: consider the control $h$ given by \eqref{*contr}: using the $L^2$ bounds \eqref{1008-3}, we obtain the existence of a control driving the initial condition $u_0$ to $0$ in time $T$, satisfying, with some $C$ independent of $\alpha \in [0,1)$ and of $u_0$:
$$ \Vert h \Vert _{L^2(a,b)\times(0,T))} \leq \frac{C}{(1-\alpha)^2}  e^{-T/C} e^{C/T} \Vert u_0 \Vert _{L^2(0,1)} ,$$
which gives the estimate \eqref{1310-12} given in Theorem \ref{prop-cost}. \qed

\end{itemize}


\subsubsection{Lower bound of the null controllability cost: proof of \eqref{1310-11}} \hfill

We are going to take advantage on the relations given by the moment method. Consider  $u_0 = - \tilde \Phi _{\alpha,2}$, and let $h$ be any control that drives the initial condition $u_0$ to $0$ in time $T$. Then we deduce from \eqref{1310-1}-\eqref{1310-2} that
\begin{equation}
\label{1310-3}
\int_0^T h_{\alpha ,1}(t) e^{\lambda_{\alpha ,1} t} dt = 0, 
\quad \text{ and } \quad
\int_0^T h_{\alpha ,2}(t) e^{\lambda_{\alpha ,2} t} dt = 1 ,
\end{equation}
where we recall that
\begin{equation}
\label{1310-4}
 h_{\alpha ,n}(t) := \int_{-1} ^1 h(x,t) \chi_{[a,b]}(x) \tilde \Phi_{\alpha ,n} (x) dx 
 = \int_{a} ^b h(x,t) \tilde \Phi_{\alpha ,n} (x) dx .
\end{equation}
The blow-up of the null controllability cost when $\alpha \to 1^-$ will come from \eqref{1310-3} and the fact that
$$ \lambda_{\alpha ,2}- \lambda_{\alpha ,1} \to 0
\quad \text{ and } \quad  \tilde \Phi_{\alpha ,2} - \tilde \Phi_{\alpha ,1} \to 0 \quad \text{ as } \alpha \to 1^-.$$
Indeed, first
\begin{multline*}
\int_0^T h_{\alpha ,2}(t) e^{\lambda_{\alpha ,1} t} \, dt
= \int_0^T ( h_{\alpha ,2}(t) - h_{\alpha ,1}(t)) e^{\lambda_{\alpha ,1} t} \, dt
\\
= \int_0^T \Bigl( \int_{a} ^b h(x,t) (\tilde \Phi_{\alpha ,2} (x)- \tilde \Phi_{\alpha ,1} (x)) dx \Bigr) e^{\lambda_{\alpha ,1} t} \, dt,
\end{multline*}
hence
\begin{multline*}
\Bigl\vert \int_0^T h_{\alpha ,2}(t) e^{\lambda_{\alpha ,1} t}\,  dt \Bigr\vert
\leq  \int_0^T \Bigl( \Vert h (\cdot,t) \Vert _{L^2(a,b)} \Vert \tilde \Phi_{\alpha ,2}-\tilde\Phi_{\alpha ,1} \Vert _{L^2(a,b)}  e^{\lambda_{\alpha ,1} t} \, dt
\\
\leq e^{\lambda_{\alpha ,1} T} \Vert \tilde \Phi_{\alpha ,2}-\tilde\Phi_{\alpha ,1} \Vert _{L^2(a,b)} \int_0^T \Vert h (\cdot,t) \Vert _{L^2(a,b)} \, dt
\\
\leq \sqrt{T} e^{\lambda_{\alpha ,1} T} \Vert \tilde \Phi_{\alpha ,2}-\tilde\Phi_{\alpha ,1} \Vert _{L^2(a,b)} \Vert h \Vert _{L^2((a,b)\times (0,T))} .
\end{multline*}
Then we deduce that
\begin{multline}
\label{1310-6}
\int_0^T h_{\alpha ,2}(t) (e^{\lambda_{\alpha ,2} t} - e^{\lambda_{\alpha ,1} t} ) \,  dt
= 1 - \int_0^T h_{\alpha ,2}(t) e^{\lambda_{\alpha ,1} t}\,  dt
\\
\geq 1- \sqrt{T} e^{\lambda_{\alpha ,1} T} \Vert \tilde \Phi_{\alpha ,2}-\tilde\Phi_{\alpha ,1} \Vert _{L^2(a,b)} \Vert h \Vert _{L^2((a,b)\times (0,T))} .
\end{multline}
On the other hand,
$$ 0 \leq  e^{\lambda_{\alpha ,2} t} - e^{\lambda_{\alpha ,1} t}  \leq e^{\lambda_{\alpha ,2} t} (\lambda_{\alpha ,2} - \lambda_{\alpha ,1}) t ,$$
hence
\begin{multline*}
\int_0^T h_{\alpha ,2}(t) (e^{\lambda_{\alpha ,2} t} - e^{\lambda_{\alpha ,1} t} ) \,  dt
\leq (\lambda_{\alpha ,2} - \lambda_{\alpha ,1}) \int_0^T \vert h_{\alpha ,2}(t) \vert 
e^{\lambda_{\alpha ,2} t} t \, dt 
\\
\leq e^{\lambda_{\alpha ,2} T} (\lambda_{\alpha ,2} - \lambda_{\alpha ,1})
\Bigl( \int _0 ^T h_{\alpha,2} (t) ^2 \, dt \Bigr) ^{1/2} \Bigl( \int _0 ^T t^2 \, dt \Bigr) ^{1/2} ,
\end{multline*}
and since
\begin{multline*}
\int _0 ^T h_{\alpha,2} (t) ^2 \, dt
=  \int _0 ^T \Bigl( \int_{a} ^b h(x,t) \tilde \Phi_{\alpha ,2} (x) dx \Bigr) ^2 \, dt 
\\
\leq  \int _0 ^T\Bigl( \int _a ^b h(x,t) ^2 \, dx \Bigr) \Bigl( \int _a ^b \tilde \Phi_{\alpha ,2} (x) ^2 dx \Bigl) \, dt 
\\
\leq \int _0 ^T\Bigl( \int _a ^b h(x,t) ^2 \, dx \Bigr) \, dt  =  \Vert h \Vert _{L^2((a,b)\times (0,T))} ^2 ,
\end{multline*}
hence
\begin{equation}
\label{1310-7}
\int_0^T h_{\alpha ,2}(t) (e^{\lambda_{\alpha ,2} t} - e^{\lambda_{\alpha ,1} t} ) \,  dt
\leq T^{3/2} e^{\lambda_{\alpha ,2} T} (\lambda_{\alpha ,2} - \lambda_{\alpha ,1})
\Vert h \Vert _{L^2((a,b)\times (0,T))}.
\end{equation}
We deduce from \eqref{1310-6}-\eqref{1310-7} that
$$ T^{3/2} e^{\lambda_{\alpha ,2} T} (\lambda_{\alpha ,2} - \lambda_{\alpha ,1})
\Vert h \Vert _{L^2(a,b)\times (0,T)}
\geq 1- \sqrt{T} e^{\lambda_{\alpha ,1} T} \Vert \tilde \Phi_{\alpha ,2}-\tilde\Phi_{\alpha ,1} \Vert _{L^2(a,b)} \Vert h \Vert _{L^2((a,b)\times (0,T))} , $$
hence
\begin{equation}
\label{1310-10}
\Vert h \Vert _{L^2((a,b)\times (0,T))}
\geq \frac{e^{-\lambda_{\alpha ,2} T}}{\sqrt{T}} \frac{1}{T (\lambda_{\alpha ,2} - \lambda_{\alpha ,1}) + \Vert \tilde \Phi_{\alpha ,2}-\tilde\Phi_{\alpha ,1} \Vert _{L^2(a,b)}} .
\end{equation}
Now we are in position to conclude:
\begin{itemize}
\item it follows from Lemma \ref{lem-appl-j_nu-croissant} that the difference $\lambda_{\alpha ,2} - \lambda_{\alpha ,1}$ behaves as $1-\alpha$:
\begin{multline}
\label{1310-70}
\lambda_{\alpha ,2} - \lambda_{\alpha ,1} = 
\kappa _\alpha ^2 (j_{\nu_\alpha ,1} ^2 - j_{-\nu_\alpha ,1} ^2)
= \kappa _\alpha ^2 (j_{\nu_\alpha ,1} - j_{-\nu_\alpha ,1})(j_{\nu_\alpha ,1} + j_{-\nu_\alpha ,1}) 
\\
\leq 2 j_{\nu_\alpha ,1} (j_{\nu_\alpha ,1} - j_{-\nu_\alpha ,1})
\leq  2 j_{\nu_\alpha ,1} M_* (1-\alpha) ;
\end{multline}

\item and from the expressions of $\tilde \Phi_{\alpha ,1}$ and $\tilde \Phi_{\alpha ,2}$
given in \eqref{Phi-alpha-n-w-i}, \eqref{Phi-alpha-n-w-p} and \eqref{1310-70}:
$$ \forall x \in (a,b), \quad 
\begin{cases}
\tilde \Phi  _{\alpha, 1} (x) 
= \frac{\sqrt{\kappa _\alpha}}{\vert J _{1-\nu _\alpha} (j_{-\nu _\alpha,1} ) \vert } x^{\frac{1-\alpha}{2}} J_{-\nu_\alpha} (j_{-\nu_\alpha,1} x^{\kappa _\alpha}) 
\\
\tilde \Phi  _{\alpha, 2} (x) = \frac{\sqrt{\kappa _\alpha}}{\vert J _{1+\nu _\alpha} (j_{\nu _\alpha,1} ) \vert } x^{\frac{1-\alpha}{2}} J_{\nu_\alpha} (j_{\nu_\alpha,1} x^{\kappa _\alpha})
\end{cases} ,$$
we deduce from the smoothness of $(\nu,x) \mapsto J_{\nu} (x)$, the fact that $\nu_\alpha =  O(1-\alpha)$ and the property that
$ j_{\nu_\alpha,1} - j_{-\nu_\alpha,1} = O(1-\alpha)$, we obtain that:
\begin{equation}
\label{1310-9}
\exists C^*, \forall x\in (a,b), \quad \vert \tilde \Phi  _{\alpha, 2}(x) - \tilde \Phi  _{\alpha, 1} (x) \vert \leq C^* (1-\alpha) . 
\end{equation}
\end{itemize}
Then we deduce from \eqref{1310-10}-\eqref{1310-9} that
$$ \Vert h \Vert _{L^2((a,b)\times (0,T))} \geq \frac{e^{-\lambda_{\alpha ,2} T}}{\sqrt{T}} \frac{1}{ 2 j_{\nu_\alpha ,1} T M_* (1-\alpha) + C^* (1-\alpha)} .$$
This gives \eqref{1310-11}, and concludes the proof of Theorem \ref{prop-cost}. \qed


\section{Numerical approximation}
\label{sec-num}

We provide in this section a finite element framework for the approximation of the controllability problem, present a direct minimization algorithm for the optimal control problem associated to controllability, and show in the numerical examples that it provides results consistent with the theory.

\subsection{Finite element scheme for the state equation}

The framework for the finite element approximation of the state equation is quite classical, but we summarize it here for completeness, following \cite{T06}. We write the weak formulation of \eqref{eq-control}, using test functions from a finite dimensional space $V_\delta$, indexed by a space discretization parameter $\delta$. In its basic form, this space will be chosen as the space of piecewise linear ($\mathbb{P}_1$) functions on $[-1,1]$, with a uniform space grid of step $\delta$ and attaining a zero value at the boundary. The weak formulation, for any $v_\delta\in V_\delta$ and $t\in (0,T)$, reads
\begin{equation}\label{eq:galerkin}
\left\langle\frac{\partial u_\delta(t)}{\partial t},v_\delta\right\rangle - a(u_\delta(t),v_\delta) = \left\langle h_\delta(t)\chi_{(a,b)},v_\delta\right\rangle
\end{equation}
where
\begin{equation}\label{eq:bilinear}
a(u_\delta(t),v_\delta) = \int_{-1}^1 |x|^\alpha \>\frac{\partial u_\delta(x,t)}{\partial x} v_\delta'(x) dx.
\end{equation}
We look for a approximations $u_\delta(x,t)$ and $h_\delta(x,t)$ of the form
\[
u_\delta(x,t) = \sum_{j=1}^N u_j(t) \phi_j(x), \quad h_\delta(x,t) = \sum_{j=j_a}^{j_b} h_j(t) \phi_j(x),
\]
with $j_a$, $j_b$ respectively the first and the last index of nodes in $[a,b]$. We also set $N_h=j_b-j_a+1$ to be the dimension of the discretized control space. \\
It clearly suffices to enforce \eqref{eq:galerkin} for all the functions $\phi_i$ which generate the space $V_\delta$. Then, taking $v_\delta=\phi_i$, the various terms of \eqref{eq:galerkin} may be written as
\begin{equation*}
\left\langle\frac{\partial u_\delta(t)}{\partial t},v_\delta\right\rangle = \sum_j \dot u_j(t) \int_{-1}^1 \phi_j(x) \phi_i(x) dx,
\end{equation*}
\begin{equation*}
a(u_\delta(t),v_\delta) = \sum_j u_j(t) \int_{-1}^1 |x|^\alpha \> \phi_j'(x) \phi_i'(x) dx.
\end{equation*}
\begin{equation*}
\left\langle h_\delta(t)\chi_{(a,b)},v_\delta\right\rangle = \sum_j h_j(t) \int_a^b \phi_j(x) \phi_i(x) dx.
\end{equation*}
Once defined the mass and stiffness matrices as
\[
M = (m_{ij}), \quad m_{ij} = \int_{-1}^1 \phi_j(x) \phi_i(x) dx,
\]
\[
A = (a_{ij}), \quad a_{ij} = \int_{-1}^1 \phi_j'(x) \phi_i'(x) dx,
\]
we obtain the semi-discrete approximation
\begin{equation}\label{eq:galerkin2}
\sum_{j=1}^N m_{ij} \dot u_j(t) - \sum_{j=1}^N a_{ij} u_j(t) = \sum_{j=j_a}^{j_b} m_{ij} h_j(t) \qquad (i=1,\ldots,N)
\end{equation}
or, in matrix form,
\begin{equation}\label{eq:matrix1}
M\dot{\bm u}(t) = A\bm u(t) + B \bm h(t),
\end{equation}
where the matrix $B\in\mathbb{R}^{N\times N_h}$ is obtained by selecting the columns $j_a$ to $j_b$ of the mass matrix. In \eqref{eq:matrix1}, $\bm u(t)=(u_1(t) \cdots u_N(t))^\top \in\mathbb{R}^N$ and $\bm h(t)=(h_{j_a}(t) \cdots h_{j_b}(t))^\top \in\mathbb{R}^{N_h}$ are the vectors of respectively the semi-discretized state and the semi-discretized control (note that, here and in what follows, $(\cdot)^\top$ will denote the transpose of a vector or matrix). The semi-discrete approximation \eqref{eq:matrix1} can be formally rewritten as
\[
\dot{\bm u}(t) = M^{-1}A\bm u(t) + M^{-1}B \bm h(t),
\]
and further discretized with respect to $t$ via a backward Euler scheme, in the form
\begin{equation}\label{eq:fully_discrete1}
{\bm u}^{n+1} = {\bm u}^n + \Delta t M^{-1}A\bm u^{n+1} + \Delta t M^{-1}B \bm h^n,
\end{equation}
for $n\in[0,N_T-1]$ and $N_T=T/\Delta t$.

Note that we have denoted by $\bm h^n$ the discrete control on the interval $(t_n,t_{n+1})$, although in the implicit setting this would be typically identified with $\bm h(t_{n+1})$. In fact, this is irrelevant since we will later use this variables to minimize the norm of the final state.

In order to obtain a more explicit form, we can write \eqref{eq:fully_discrete1} as
\begin{equation*}
(I-\Delta t M^{-1}A){\bm u}^{n+1} = {\bm u}^n + \Delta t M^{-1}B \bm h^n,
\end{equation*}
that is,
\begin{eqnarray}\label{eq:fully_discrete2}
{\bm u}^{n+1} & = & (I-\Delta t M^{-1}A)^{-1}{\bm u}^n + \Delta t (I-\Delta t M^{-1}A)^{-1}M^{-1}B \bm h^n \nonumber \\
& = & S{\bm u}^n + R \bm h^n,
\end{eqnarray}
where $\bm u^n=(u_1^n \cdots u_N^n)^\top$ and $\bm h^n=(h_{j_a}^n \cdots h_{j_b}^n)^\top$, and
\[
S = (I-\Delta t M^{-1}A)^{-1},\qquad R = \Delta t S M^{-1}B,
\]
with an initial condition $\bm u^0$ defined as a suitable projection of $u_0(x)$ on the space $V_\delta$.

Convergence of the scheme is proved by standard techniques for nondegenerate diffusions. The degenerate case may be treated along the guidelines in \cite[Chap. 18]{T06} to obtain convergence in the space $H^1_\alpha$.

\subsection{Direct minimization algorithm}

The zero-controllability condition has been numerically implemented by a discrete optimal control problem of Mayer type, in which we minimize the norm of the fully discrete solution in the space $V_\delta$, i.e.,
\begin{equation}\label{eq:mayer}
\min_{\bm h^0, \ldots, \bm h^{N_T-1}} J_\delta(\bm h^0, \ldots, \bm h^{N_T-1}) = \min_{\bm h^0, \ldots, \bm h^{N_T-1}} \frac 1 2 {\bm u^{N_T}}^\top M \bm u^{N_T}
\end{equation}
in which minimization has been carried out in a discretize-then-optimize strategy. We are not interested here in computing minimal-norm controls, as obtained for example via the Hilbert Uniqueness Method (see \cite{Lions}).

The gradient of the discretized functional $J_\delta$ is computed using the discrete adjoint problem, that is, the backward difference equation
\begin{equation}\label{eq:adj}
\begin{cases}
\bm p^n = S^\top \bm p^{n+1} \\
\bm p^{N_T} = M \bm u^{N_T}
\end{cases}
\end{equation}
and computing the derivatives with respect to the discrete controls at time $t_k$ as
\begin{equation}\label{eq:grad}
\frac{\partial J}{\partial \bm h^k} = R^\top \bm p^k.
\end{equation}
Note that the form \eqref{eq:adj}--\eqref{eq:grad} results from a completely standard computation on the discrete problem, but might also be obtained by using the same discretization of the state equation on the continuous adjoint equation.

Once computed the gradient, algorithms of Conjugate Gradient or Quasi-Newton type may be applied. In our case, we have used the BFGS Quasi-Newton algorithm. A possible lack of uniqueness of the minimum is known to be correctly handled by this kind of methods.

\subsection{Numerical examples}

We provide two numerical tests with $a=1/2$, $b=3/4$, i.e., with the state equation
\begin{equation}
\begin{cases}
u_t - (|x|^\alpha u_x)_x = \chi_{(1/2,3/4)} h(x,t) \\
u(x,0) = u_0(x),
\end{cases}
\end{equation}
for $T=0.5$ and different values of $\alpha$, with $u_0$ having its support respectively on the left and on the right of the degeneracy. The test has been carried out with $\delta=0.02$ and $\Delta t= 2.5\cdot 10^{-3}$.

\subsubsection*{Example 1.}

\begin{figure}[htbp]
\begin{center}
\includegraphics[width=5 cm, height=4 cm]{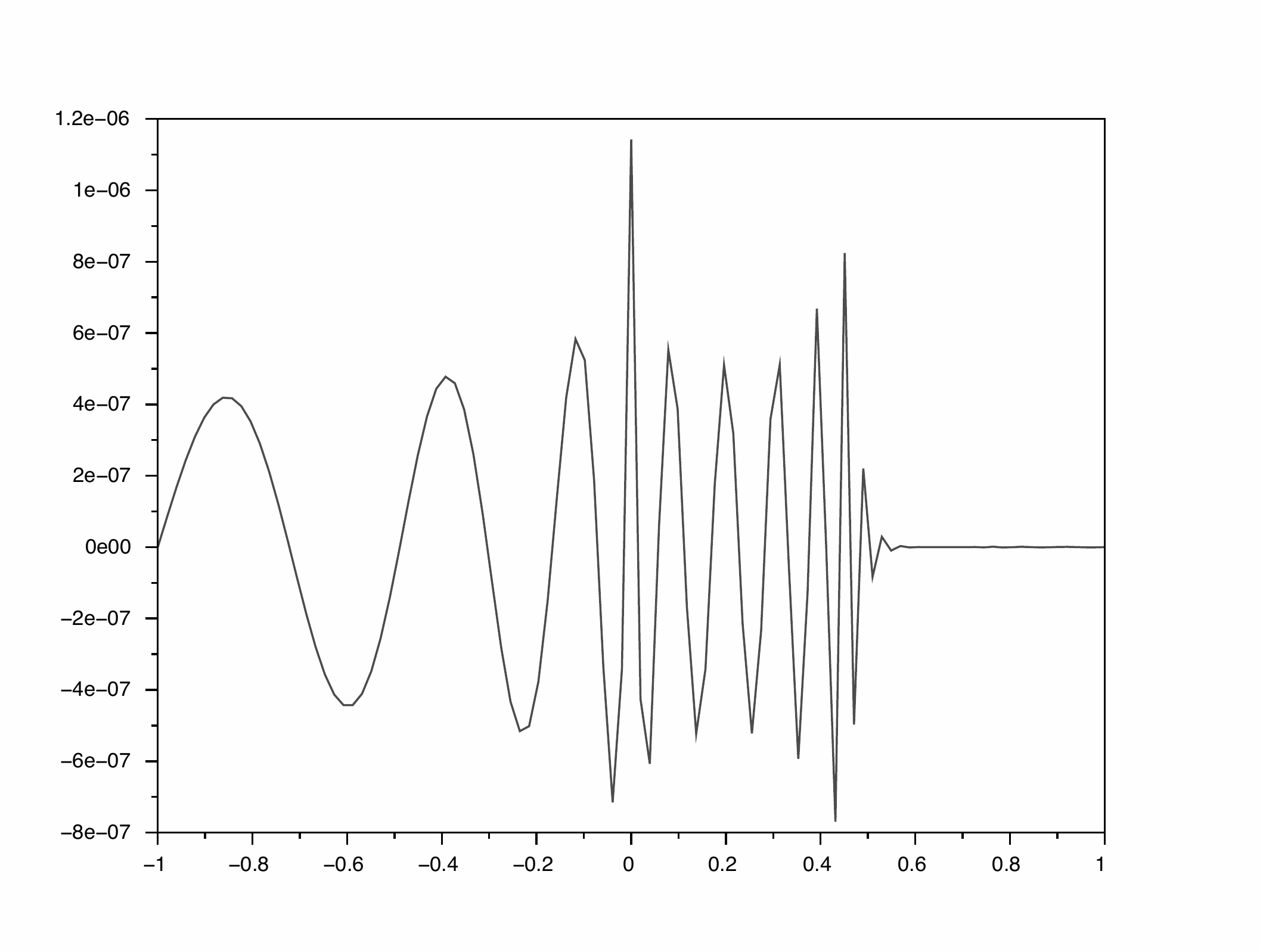}
\includegraphics[width=5 cm, height=4 cm]{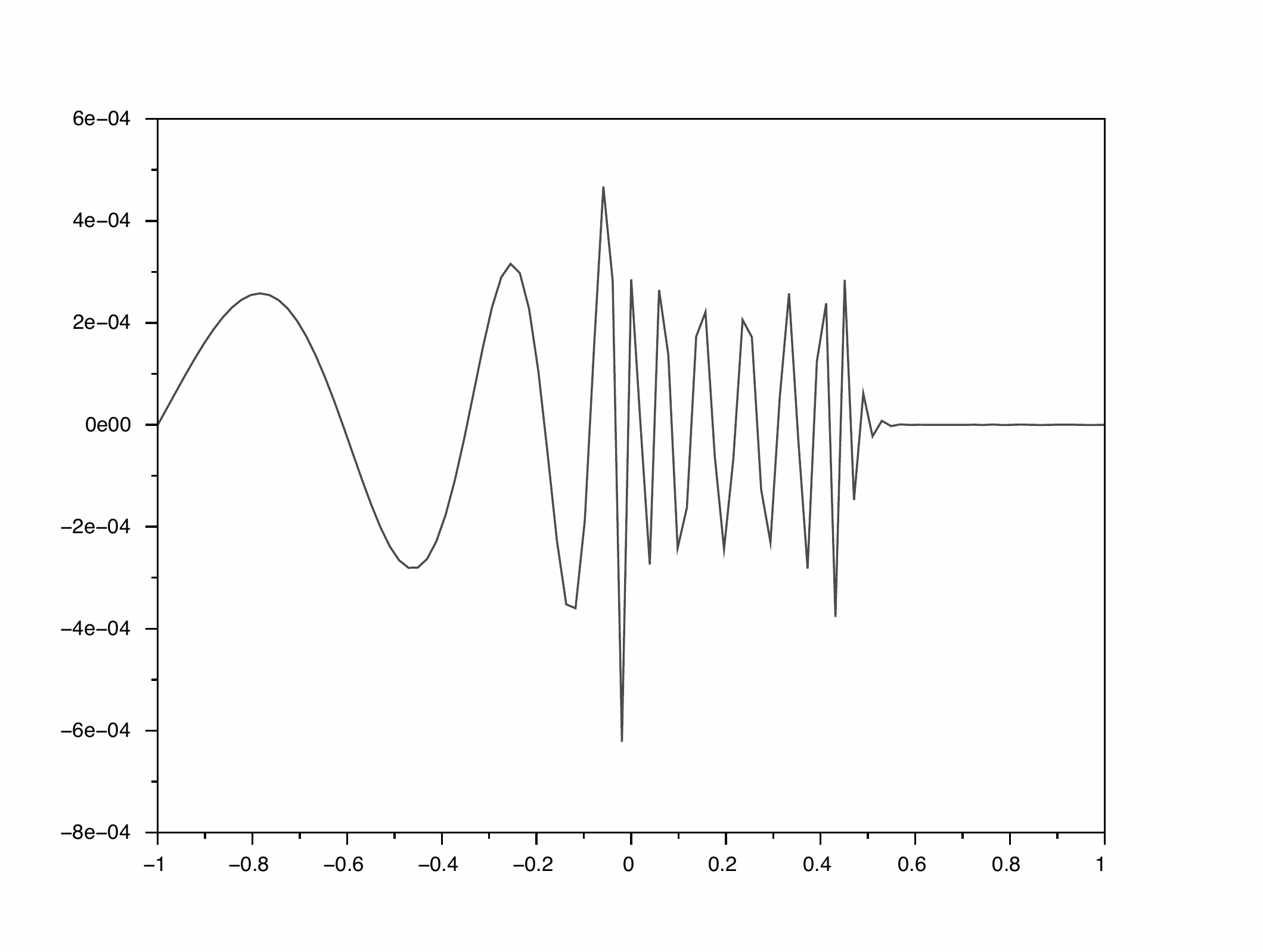}
\caption{Test 1, final state for $\alpha=0.8$ (left) and $\alpha=1.2$ (right).}
\label{fig:final_state_test1}
\end{center}
\end{figure}

\begin{figure}[htbp]
\begin{center}
\includegraphics[width=5 cm, height=4 cm]{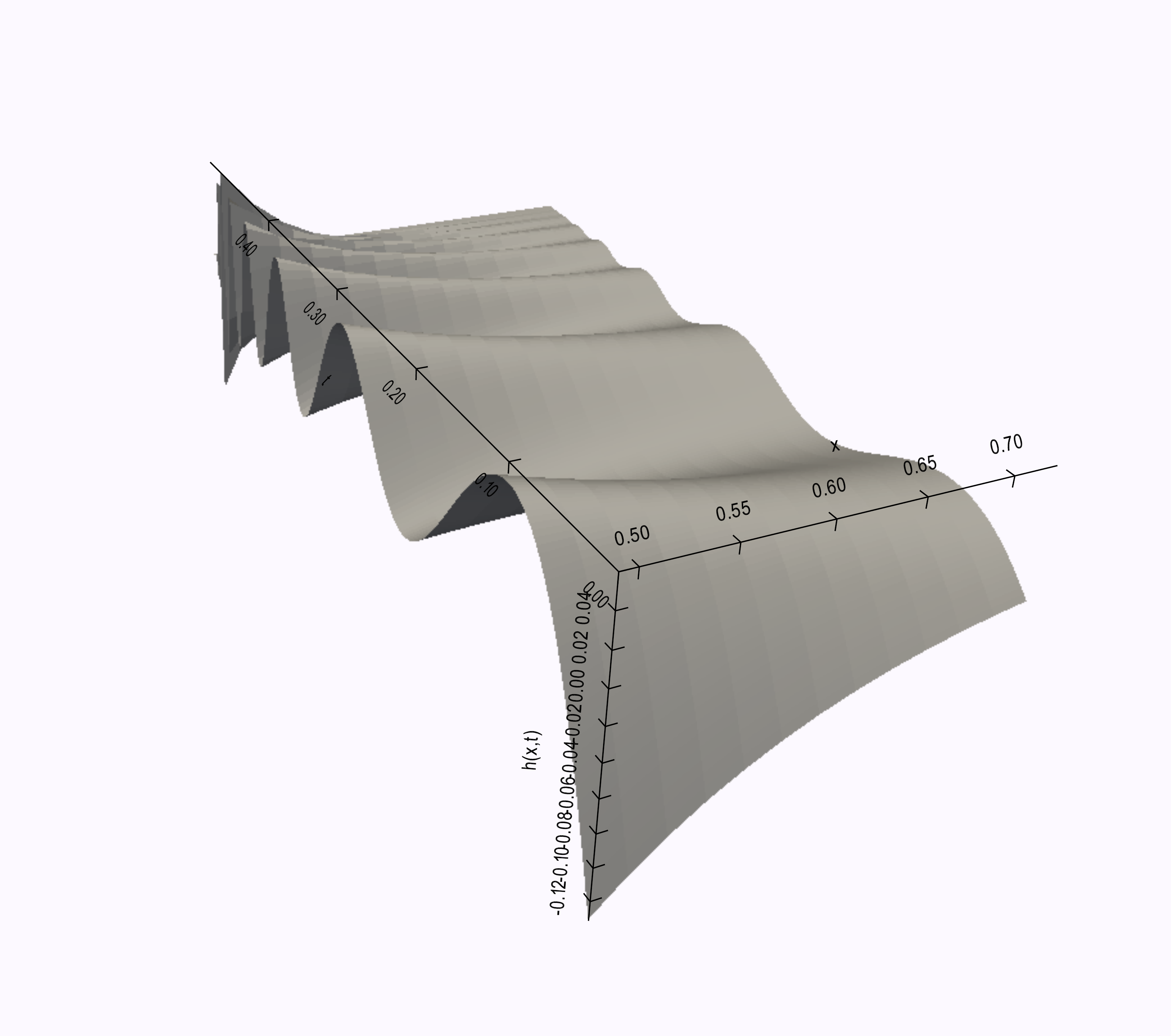}
\includegraphics[width=5 cm, height=4 cm]{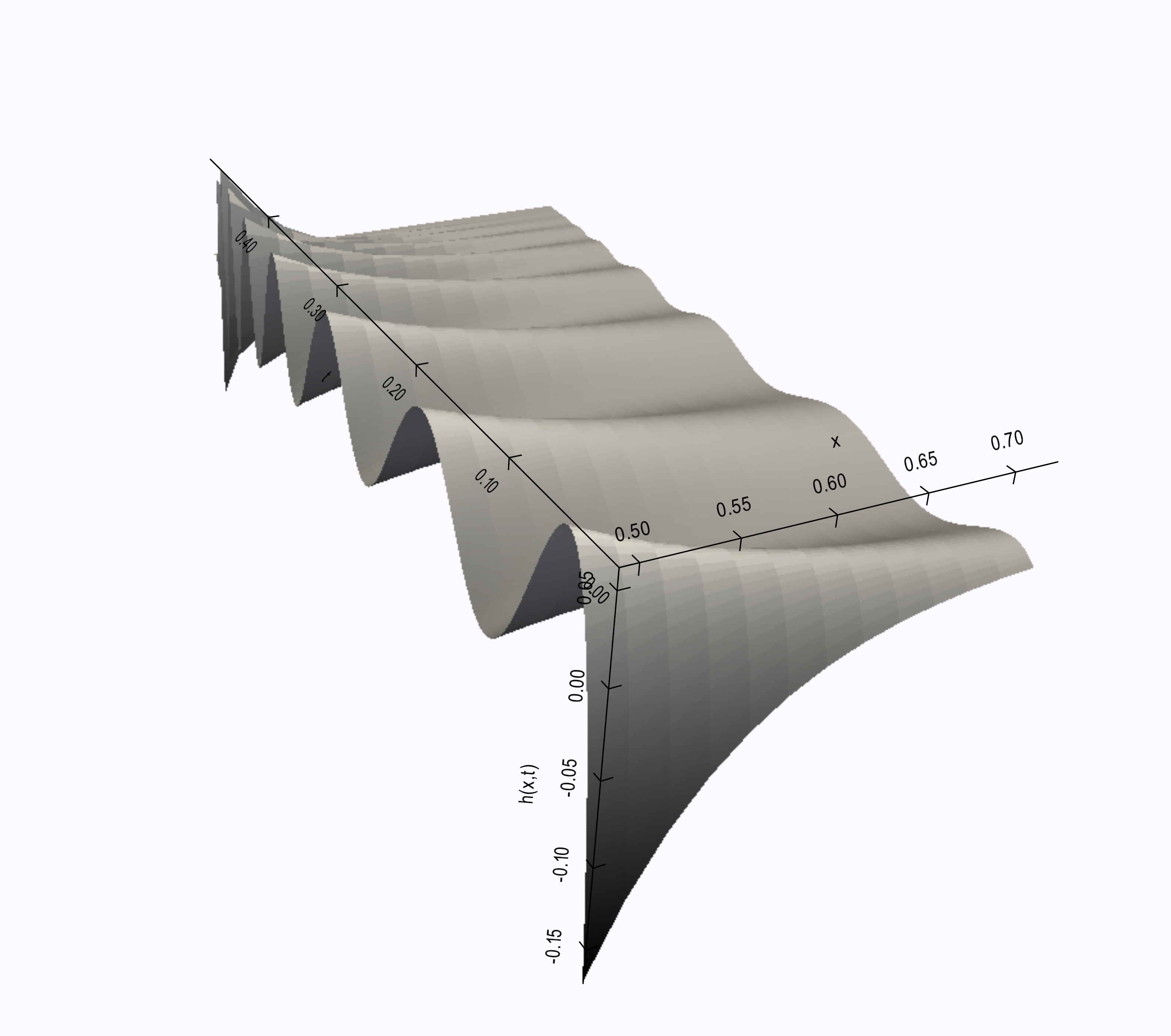}
\caption{Test 1, control $h(x,t)$ for $\alpha=0.8$ (left) and $\alpha=1.2$ (right). The scaling factor is respectively $10^{-5}$ and $2\cdot 10^{-8}$.}
\label{fig:h_test1}
\end{center}
\end{figure}

In this test, we take an initial state supported in ${\mathbb R}_-$:
\begin{equation}
u_0(x) = \chi_{(-1/2,-1/4)}.
\end{equation}
The two typical situations may be seen in Figg. \ref{fig:final_state_test1}--\ref{fig:h_test1}, where the former shows the final state at $T=0.5$ for $\alpha=0.8$ (left) and $\alpha=1.2$ (right), while the latter shows the corresponding control $h$ for $(x,t)\in [a,b]\times [0,T]$ (the plots of $h$ have been scaled by a factor of respectively $10^{-5}$ for $\alpha=0.8$ and $2\cdot 10^{-8}$ for $\alpha=1.2$). Note that, although the qualitative behaviour in the two cases is similar, the scales are very different: the final state has an $L^2$ norm of about $4.3\cdot 10^{-7}$ for $\alpha=0.8$ versus a norm of $2.4\cdot 10^{-4}$ for $\alpha=1.2$. The optimal control maintains the known feature of being highly oscillating, but its norm is considerably higher in the second case -- this reflects the transition between a zero-controllable situation and a non-controllable one.

%

\subsubsection*{Example 2.}

\begin{figure}[htbp]
\begin{center}
\includegraphics[width=5 cm, height=4 cm]{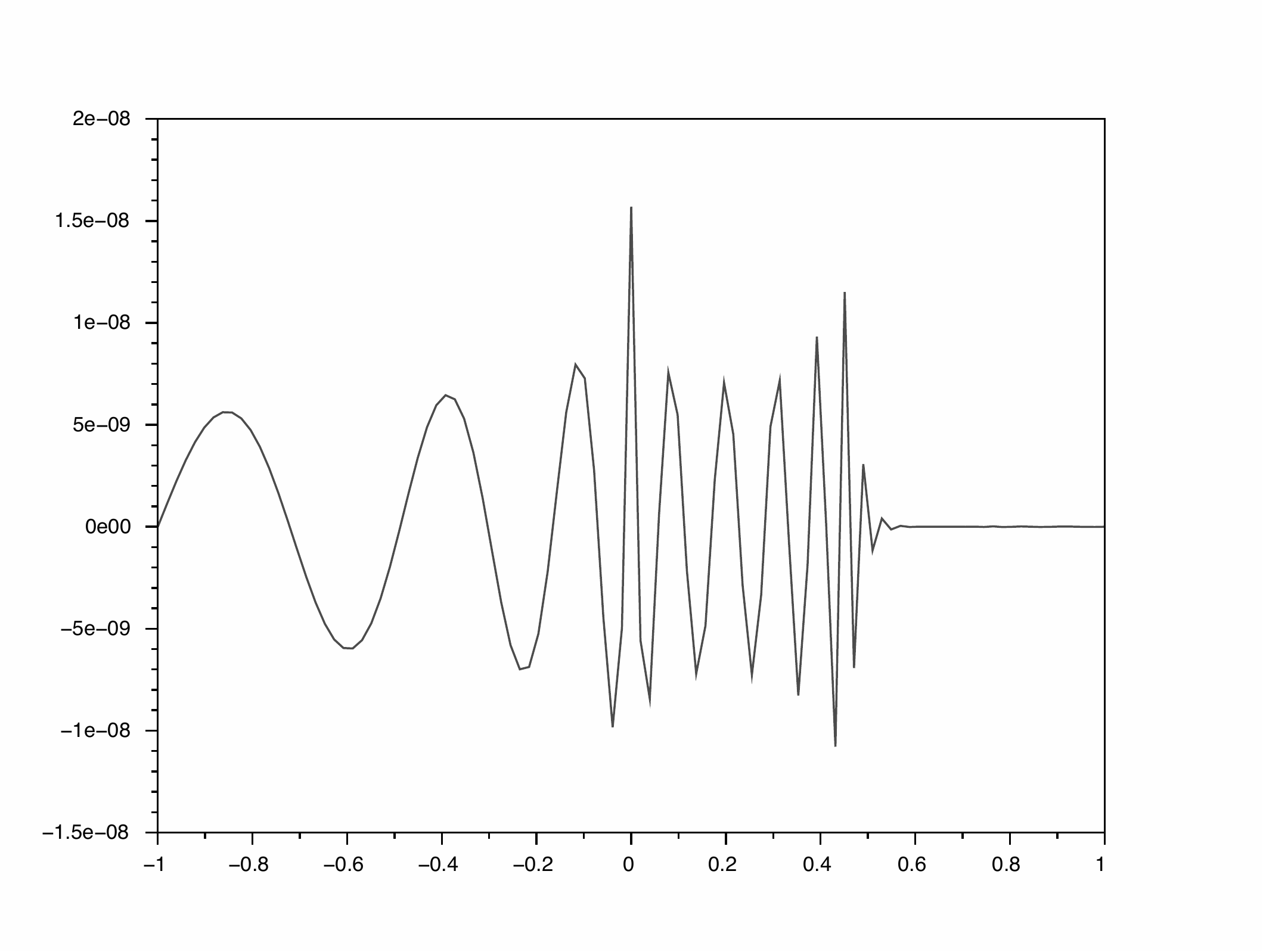}
\includegraphics[width=5 cm, height=4 cm]{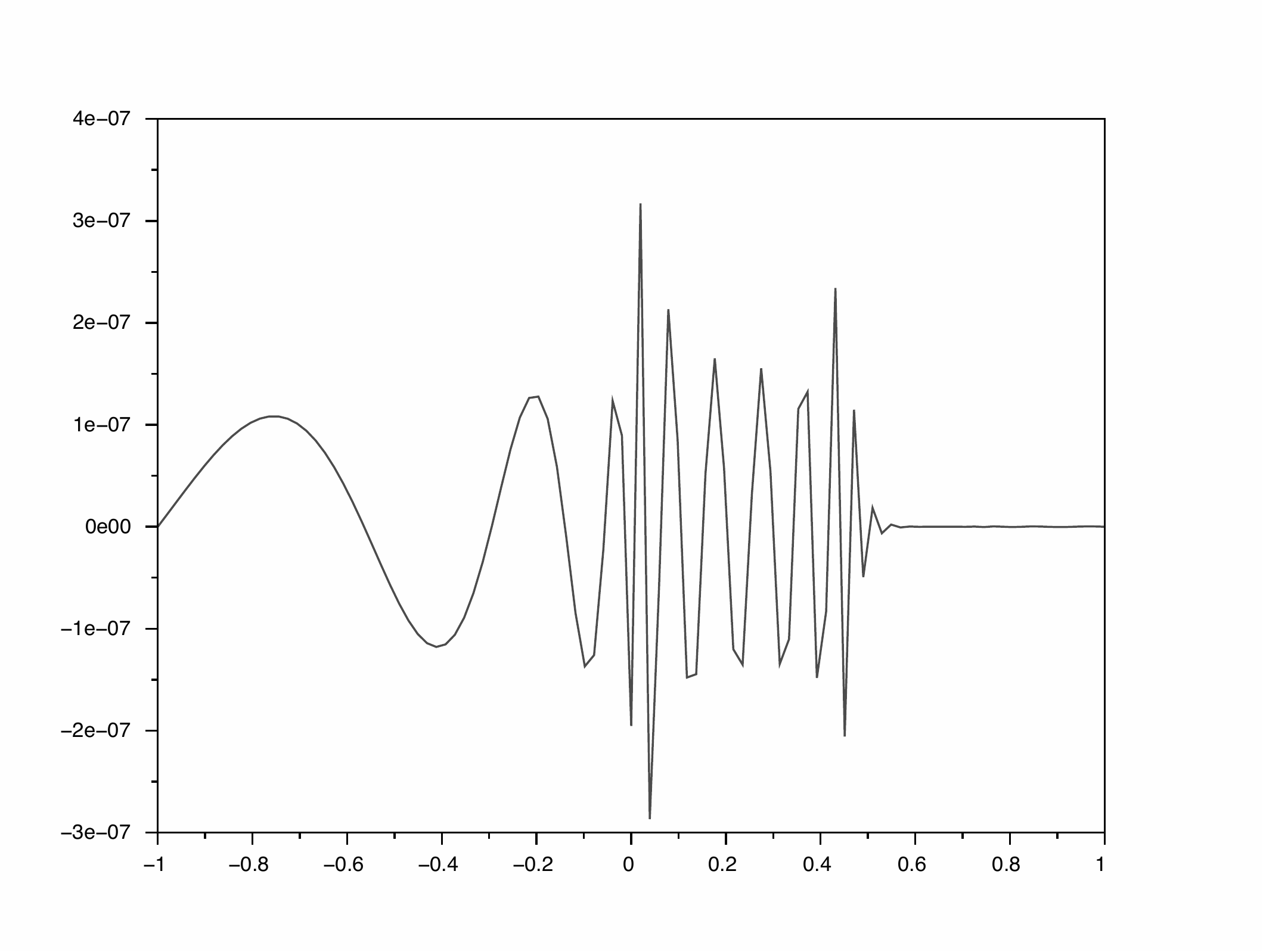}
\caption{Test 2, final state for $\alpha=0.8$ (left) and $\alpha=1.2$ (right).}
\label{fig:final_state_test2}
\end{center}
\end{figure}

\begin{figure}[htbp]
\begin{center}
\includegraphics[width=5 cm, height=4 cm]{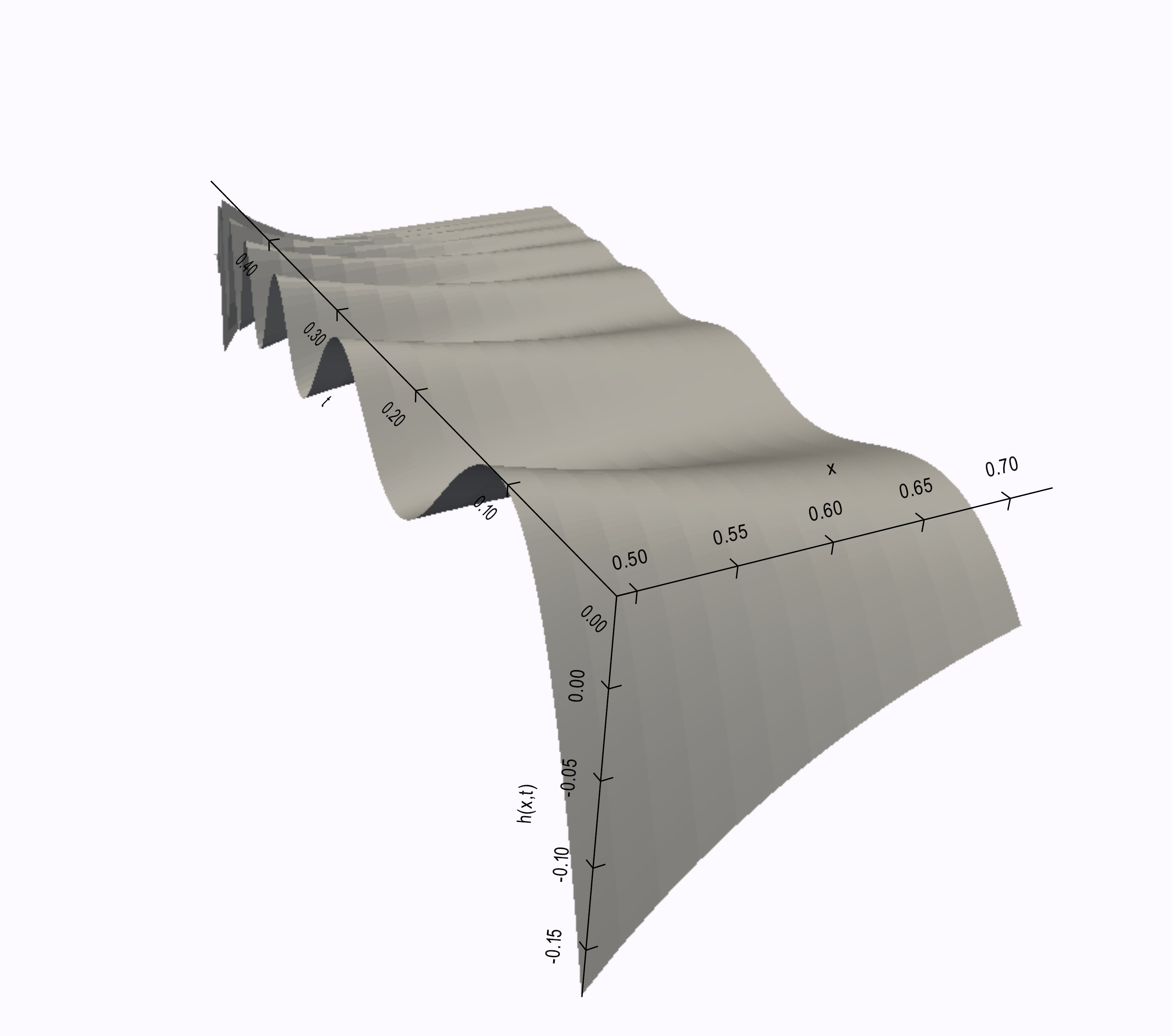}
\includegraphics[width=5 cm, height=4 cm]{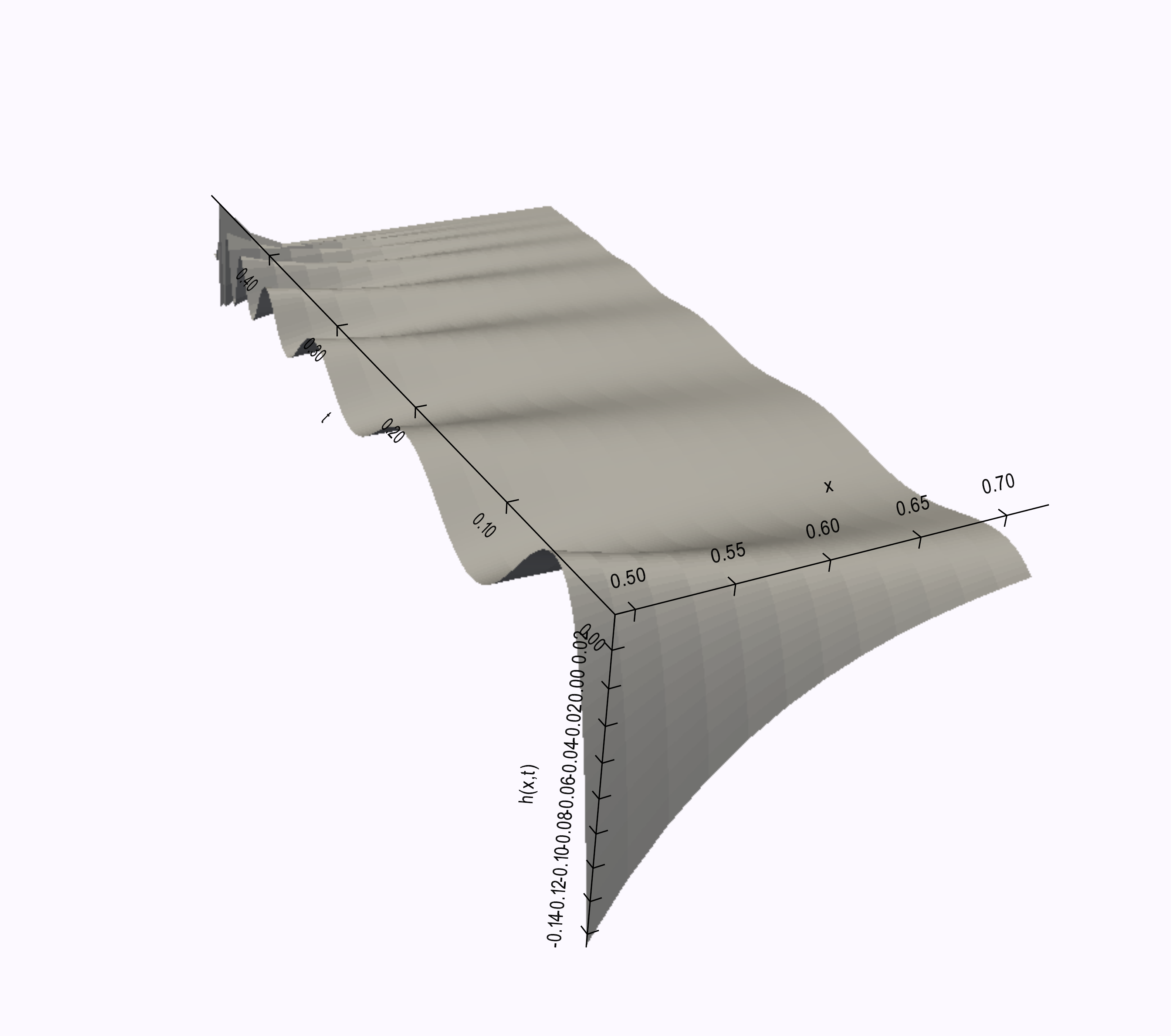}
\caption{Test 2, control $h(x,t)$ for $\alpha=0.8$ (left) and $\alpha=1.2$ (right). The scaling factor is respectively $5\cdot 10^{-4}$ and $5\cdot 10^{-5}$.}
\label{fig:h_test2}
\end{center}
\end{figure}

In the second numerical example, the initial state is supported in ${\mathbb R}_+$, and more precisely
\begin{equation}
u_0(x) = \chi_{(1/4,1/2)},
\end{equation}
which results in a zero-controllable state for all values of $\alpha$. The results are shown in Figg. \ref{fig:final_state_test2}--\ref{fig:h_test2}, for both $\alpha=0.8$ (left) and $\alpha=1.2$ (right), where the plots of $h$ have been scaled by a factor of respectively $5\cdot 10^{-4}$ for $\alpha=0.8$ and $5\cdot 10^{-5}$ for $\alpha=1.2$.

In this case, the initial state is controllable and the scales are similar: the final state has an $L^2$ norm of about $6\cdot 10^{-9}$ for $\alpha=0.8$ and about $10^{-7}$ for $\alpha=1.2$, according to the theoretical prediction. Note that, since the initial state is symmetric to that of Example 1, in both cases we would have the same final norm for the uncontrolled state equation.

%

As a final remark, we note that the agreement between theory and numerics must be understood in a qualitative form: numerical tests show that the value $\alpha=1$ works as a threshold between controllable and non-controllable situations. However, in the numerical setting under consideration the transition between the two cases appears to be continuous, and the influence of degeneracy on the null controllability of the discrete system is definitely worth a deeper investigation.


\end{document}